\documentclass[11pt]{article}
\usepackage{amsmath,epsfig,pstricks,pst-node}

\newenvironment{fullfigure}[2]
    {\begin{figure}[htb]\begin{center}\def\ffa{#1}\def\ffb{#2}}
    {\vspace{\baselineskip}\caption{\ffb.}\label{\ffa}\end{center}\end{figure}}

\newcommand{\eatline}{\vspace{-\baselineskip}}

\newcommand{\verta}[1]{\pspicture[.6](-.7,-1.9)(.7,.7)
\pcline(0,0)(0,.7)\mto    \pcline(0,0)(.7,0)\mfro
\pcline(0,0)(-.7,0)\mfro  \pcline(0,0)(0,-.7)\mto
\rput(0,-1.2){$#1$}\endpspicture}
\newcommand{\vertb}[1]{\pspicture[.6](-.7,-1.9)(.7,.7)
\pcline(0,0)(0,.7)\mfro   \pcline(0,0)(.7,0)\mto
\pcline(0,0)(-.7,0)\mto   \pcline(0,0)(0,-.7)\mfro
\rput(0,-1.2){$#1$}\endpspicture}
\newcommand{\vertc}[1]{\pspicture[.6](-.7,-1.9)(.7,.7)
\pcline(0,0)(0,.7)\mto    \pcline(0,0)(.7,0)\mto
\pcline(0,0)(-.7,0)\mfro  \pcline(0,0)(0,-.7)\mfro
\rput(0,-1.2){$#1$}\endpspicture}
\newcommand{\vertd}[1]{\pspicture[.6](-.7,-1.9)(.7,.7)
\pcline(0,0)(0,.7)\mfro   \pcline(0,0)(.7,0)\mfro
\pcline(0,0)(-.7,0)\mto   \pcline(0,0)(0,-.7)\mto
\rput(0,-1.2){$#1$}\endpspicture}
\newcommand{\verte}[1]{\pspicture[.6](-.7,-1.9)(.7,.7)
\pcline(0,0)(0,.7)\mto    \pcline(0,0)(.7,0)\mfro
\pcline(0,0)(-.7,0)\mto   \pcline(0,0)(0,-.7)\mfro
\rput(0,-1.2){$#1$}\endpspicture}
\newcommand{\vertf}[1]{\pspicture[.6](-.7,-1.9)(.7,.7)
\pcline(0,0)(0,.7)\mfro   \pcline(0,0)(.7,0)\mto
\pcline(0,0)(-.7,0)\mfro  \pcline(0,0)(0,-.7)\mto
\rput(0,-1.2){$#1$}\endpspicture}

\newcommand{\mto}{\lput{:U}{\pspicture(0,0)(0,0)
\psline[arrows=->](2.3pt,0)(2.4pt,0)\endpspicture}}
\newcommand{\mfro}{\lput{:U}{\pspicture(0,0)(0,0)
\psline[arrows=->](-2.3pt,0)(-2.4pt,0)\endpspicture}}

\newtheorem{Theorem}{Theorem}[section]
\newtheorem{Corollary}[Theorem]{Corollary}

\newtheorem{Proposition}[Theorem]{Proposition}

\newcommand{\beq}{\begin{equation}}
\newcommand{\eeq}{\end{equation}}

\def\ssc{\scriptscriptstyle}

\def\ov{\overline}

\def\ph1{{\phantom1}}
\def\ph-{{\phantom-}}

\def\x{{\bf x}}
\def\y{{\bf y}}
\def\z{{\bf z}}
\def\u{{\bf u}}
\def\v{{\bf v}}
\def\w{{\bf w}}
\def\p{{\bf p}}
\def\q{{\bf q}}
\def\a{{\bf a}}
\def\b{{\bf b}}

\def\ra{\rightarrow}

\def\con{{\rm con}}

\def\wgt{{\rm wgt}}
\def\hgt{{\rm hgt}}
\def\str{{\rm str}}

\def\bar{{\rm bar}}

\def\row{{\rm row}}

\def\={\!\!=\!\!}
\def\+{\!\!+\!\!}
\def\-{\!\!-\!\!}
\def\0{\phantom{0}}

\def\mystrut{\hbox{\vrule height10.6pt depth2pt width0pt}}
\def\mybox{\hbox to 8.9pt}

\def\mywidebox{\hbox to 18.9pt}
\def\norulefill{\leaders\hrule height0pt\hfill}
\def\nr#1{\multispan{#1}\norulefill}
\def\hr#1{\multispan{#1}\hrulefill}

\def\ov{\overline}

\def\ph1{{\phantom1}}
\def\ph-{{\phantom-}}

\def\con{{\rm con}}

\def\wgt{{\rm wgt}}
\def\hgt{{\rm hgt}}
\def\str{{\rm str}}

\def\bar{{\rm bar}}

\def\row{{\rm row}}

\def\var{{\rm var}}
\def\col{{\rm col}}

\def\0{\phantom{0}}

\def\ra{\ \longrightarrow\ }

\def\mystrut{\hbox{\vrule height10.6pt depth2pt width0pt}}
\def\mybox{\hbox to 8.9pt}

\def\mywidebox{\hbox to 18.9pt}
\def\norulefill{\leaders\hrule height0pt\hfill}
\def\nr#1{\multispan{#1}\norulefill}
\def\hr#1{\multispan{#1}\hrulefill}

\begin{document}

\title{Bijective proofs of shifted tableau 
and alternating sign matrix identities}
\author{A M Hamel\thanks{e-mail: ahamel@wlu.ca} \\
Department of Physics and Computer Science,\\
Wilfrid Laurier University, 
Waterloo, Ontario N2L 3C5, Canada\\
\\
and \\
\\
R C King\thanks{e-mail: R.C.King@maths.soton.ac.uk} \\
School of Mathematics, University of Southampton, \\
Southampton SO17 1BJ, England \\
}

\maketitle


\begin{abstract}
We give a bijective proof of an identity relating primed shifted
$gl(n)$-standard tableaux to the product of a $gl(n)$ character in the
form of a Schur function and $\prod_{1\leq i < j \leq n} (x_i + y_j)$.
This result generalises a number of well--known results due to 
Robbins and Rumsey, Chapman, Tokuyama, Okada and Macdonald. 
An analogous result is then obtained in the case of
primed shifted $sp(2n)$-standard tableaux which are bijectively
related to the product of a $t$-deformed $sp(2n)$ character and 
$\prod_{1\leq i < j \leq n}(x_i+t^2x_i^{-1}+y_j+t^2y_j^{-1})$. 
All results are also interpreted in terms of alternating sign matrix 
(ASM) identities,  including a result regarding subsets of ASMs
specified by conditions on certain restricted column sums.
\end{abstract}

\section{Introduction}
The expression
\begin{equation}
\prod_{1\leq i< j\leq n} (x_i+y_j)
\label{Eq-xy}
\end{equation}
appears in a number of contexts in symmetric function theory. Given
$\y=(y_1,y_2,\ldots, y_n)$ and $\x=(x_1,x_2,\ldots,x_n)$,
when ${\bf y} = -{\bf x}$, the expression (\ref{Eq-xy}) is just
the Vandermonde determinant that appears in Weyl's denominator formula 
\begin{equation}
\det(x_{i}^{n-j}) = \prod_{1\leq i< j\leq n} (x_i -x_j).
\label{Eq-Weyl}
\end{equation}

For ${\bf y} = \lambda {\bf x}$, the expression (\ref{Eq-xy}) becomes the
subject of the $\lambda$--determinant 
formula of Robbins
and Rumsey~\cite{RR86}:
\begin{equation}
\prod_{1\leq i< j\leq n} (x_i + \lambda x_j) = \sum_{A\in {\cal A}_{n}} 
\lambda^{SE(A)} (1+\lambda)^{NS(A)} \prod_{i=1}^{n} x_{i}^{NE_{i}(A)+ SE_{i}(A)
+ NS_{i}(A)},
\label{Eq-RR}
\end{equation}
where the exponents are various parameters associated with alternating sign
matrices
and defined in Section \ref{sec:gl-bijection}.  
Robbins and Runsey use different notation but do include the square
ice concepts, 
although they use different terminology.
Bressoud \cite{B01} asked for
a combinatorial proof of (\ref{Eq-RR}). This was provided by Chapman \cite{C01}
who generalised it to:
\begin{equation}
\prod_{1\leq i < j \leq n} (x_i+y_j) = \sum_{A\in {\cal A}_{n}}
\prod_{i=1}^{n}
x_{i}^{NE_{i}(A)} y_{i}^{SE_{i}(A)} (x_i + y_i)^{NS_{i}(A)}.
\label{Eq-Chapman}
\end{equation}
For ${\bf y} = t {\bf x}$, there is also the $t$-deformation of a
Weyl denominator formula for $gl(n)$ due to Tokuyama \cite{T88}:
\begin{equation}
\prod_{i=1}^n x_i\ \prod_{1\leq i < j \leq n} (x_i+tx_j)\ s_{\lambda}(\x)=
\sum_{ST\in {\cal ST}^{\mu}(n)} t^{\hgt(ST)} (1+t)^{\str(ST)-n}\ \x^{\wgt(ST)},
\label{Eq-Tokuyama}
\end{equation}
where the sum is over semistandard shifted tableaux $ST$ of shape 
$\mu=\lambda+\delta$ with $\delta=(n,n-1,\ldots,1)$, and where
$\hgt,\; \str,\; $ and $\wgt$ are parameters associated with semistandard
shifted tableaux. They are defined in Section \ref{sec:background}.
Note also that $s_\lambda(\x)$, the Schur function specified by
the partition $\lambda$, with a suitable interpretation of the
indeterminates $x_i$ for $i=1,2,\ldots,n$, is the character
of an irreducible representation of $gl(n)$ whose highest weight
is specified by the partition $\lambda$. 

Here we present a general identity that unifies the results
(\ref{Eq-Weyl})-(\ref{Eq-Tokuyama}). This identity 
is our first main result and is expressed 
in terms of a certain generalisation of Schur $P$--functions and
also in terms of the corresponding generalisation of Schur $Q$--functions.
These $P$  and $Q$  functions are defined combinatorially in 
Section~\ref{sec:background}.
%
%

\begin{Proposition}
Let $\mu=\lambda+\delta$ be a strict partition of length $\ell(\mu)=n$,
with $\lambda$ a partition of length $\ell(\lambda)\leq n$ and 
$\delta=(n,n-1,\ldots,1)$. In addition, 
let $\x=(x_1,x_2,\ldots,x_n)$ and $\y=(y_1,y_2,\ldots,y_n)$. Then
\beq
\begin{array}{rcl}
P_\mu(\x/\y)&=& s_\lambda(\x)\ \prod_{i=1}^n\ x_i\ 
\prod_{1\leq i<j\leq n}\ (x_i+y_j),\\ \\
Q_\mu(\x/\y)&=& s_\lambda(\x)\ 
\prod_{1\leq i\leq j\leq n}\ (x_i+y_j),\\ 
\end{array}
\label{Eq-PQ}
\eeq
where $P_\mu(\x/\y)$ and $Q_\mu(\x/\y)$ are as defined in 
Section~\ref{sec:background}.
\label{Prop-PQ}
\end{Proposition}

A bijective proof of this Proposition, along with a number of corollaries, is
provided in Section~\ref{sec:gl-bijection}. The case $\x=\y$
is an example of Macdonald  \cite{M95} (Ex2, p259, 2nd Edition). The case
$\y=t\x=(tx_1,tx_2,\ldots,tx_n)$ is equivalent to 
a Weyl denominator deformation theorem due to Tokuyama~\cite{T88}
for the Lie algebra $gl(n)$ and given a combinatorial proof by Okada~\cite{O90}.
The case $\lambda=0$ 
is equivalent to an alternating sign matrix (ASM) identity attributed to 
Robbins and Rumsey~\cite{RR86} and proved combinatorially 
by Chapman \cite{C01}. The connection with ASMs is provided in
Section~\ref{sec:asm}, in which both (\ref{Eq-RR}) and (\ref{Eq-Chapman})
are shown to be simple corollaries of Proposition~\ref{Prop-PQ}. 

It should be pointed out that the above Proposition is restricted
to the case of a strict partition $\mu$ of length $\ell(\mu)=n$.
Although a similar result applying to the case $\ell(\mu)=n-1$
may be obtained from the above by dividing both sides
by $s_{1^n}(\x)=x_1x_2\cdots x_n$, there is no similar product formula
for either $P_\mu(\x/\y)$ or $Q_\mu(\x/\y)$ in the case $\ell(\mu)<n-1$.

On the other hand, the above results may all be generalised to the case
of certain symplectic tableaux. The analogue of (\ref{Eq-xy}) in this
setting turns out to be
\begin{equation}
     \prod_{1\leq i< j\leq n}\ (x_i+t^2x_i^{-1}+y_j+t^2y_j^{-1}).
\label{Eq-xybar}
\end{equation}
When $\y=-\x$ and $t=-1$ the expression (\ref{Eq-xybar}) is a factor 
of the determinant that appears in Weyl's denominator formula for $sp(2n)$,
\begin{equation}
     \det(x_i^{n-j+1}-x_i^{-n+j-1}) =  
   \prod_{i=1}^n\ (x_i-x^{-1}_i)\  \prod_{1\leq i<j\leq n}\ (x_i+x_i^{-1}-x_j-x_j^{-1}).
\label{Eq-Weyl-sp}
\end{equation} 
More generally, for $\y=t\x$ we have~\cite{HK02}
\begin{equation}
\begin{array}{l} 
\prod_{i=1}^n (x_i+tx^{-1}_i)\  
  \prod_{1\leq i<j\leq n}\ (x_i+t^2x_i^{-1}+tx_j+tx_j^{-1})\ 
   sp_\lambda(\x;t) \\ \\
 \qquad\qquad =\  \sum_{ST\in{\cal ST}^\mu(n,{\ov n})} \
   t^{\var(ST)+\bar(ST)}\ (1+t)^{\str(ST)-n}\ \x^{\wgt(ST)},\\
\end{array}
\label{Eq-HK}
\end{equation}
where the sum is over semistandard shifted symplectic tableaux of shape
$\mu=\lambda+\delta$ with $\delta=(n,n-1,\ldots,1)$, and where $\var$, $\bar$,
$\str$ and $\wgt$ are defined in Section~\ref{sec:background}. Here 
$sp_\lambda(\x;t)$, once again
with a suitable interpretation of the indeterminates $x_i$ for $i=1,2,\ldots,n$,
is a $t$-deformation of the character $sp_\lambda(\x)$ of the
irreducible representation of the Lie algebra $sp(2n)$ whose highest weight is
specified by the partition $\lambda$.

Our second main result then takes the form

\begin{Proposition}
Let $\mu=\lambda+\delta$ be a strict partition of length $\ell(\mu)=n$,
with $\lambda$ a partition of length $\ell(\lambda)\leq n$ and 
$\delta=(n,n-1,\ldots,1)$. In addition, 
let $\x=(x_1,x_2,\ldots,x_n)$, $\y=(y_1,y_2,\ldots,y_n)$, 
${\ov\x}=({\ov x}_1,{\ov x}_2,\ldots,{\ov x}_n)$ and  
${\ov\y}=({\ov y}_1,{\ov y}_2,\ldots,{\ov y}_n)$, with ${\ov x}_i=x_i^{-1}$
and ${\ov y}_i=y_i^{-1}$ for $i=1,2,\ldots,n$. Then
\beq
Q_\mu(\x/\y;t)= sp_\lambda(\x;t)\ 
\prod_{1\leq i\leq j\leq n}\ 
(x_i+t^2{\ov x}_i+y_j+t^2{\ov y}_j),
\label{Eq-MainResult-sp}
\eeq
where $Q_\mu(\x/\y;t)$ is defined in 
Section~\ref{sec:background}.
\label{Prop-PQsp}
\end{Proposition}

Here $Q(\x/\y;t)$ is a generalisation of $Q(\x/\y)$ that associates
factors of $t^2$ with the barred components of $\ov\x$ and $\ov\y$.
Although a similar generalisation $P(\x/\y;t)$ of $P(\x/\y)$
exists, as we shall see, there does not exist a corresponding identity for
$P(\x/\y;t)$ that is analogous to the identity (\ref{Eq-MainResult-sp})
for $Q(\x/\y;t)$.

Our paper is arranged as follows.  In Section~\ref{sec:background} the necessary
background is introduced regarding both the relevant semistandard, shifted 
and primed tableaux, and the various $P$ and $Q$ functions and characters 
of $gl(n)$ and $sp(2n)$. For the $gl(n)$ case, Section~\ref{sec:gl-bijection} 
opens in Section \ref{Main Result} with a formal statement of the combinatorial identity
upon which the first main result, Proposition~\ref{Prop-PQ}, is based. 
A bijective proof of this identity is then provided. A detailed example appears in Section \ref{subsec:gl-example}. 
In Section~\ref{subsec:gl-corollaries} a number of corollaries are gathered together.
 
Turning to the $sp(2n)$ case, the combinatorial identity necessary to establish
the second main result, Proposition~\ref{Prop-PQsp}, is stated, bijectively 
proved and exemplified in Section~\ref{sec:sp-bijection}. Once again two
corollaries are supplied in Section~\ref{subsec:sp-corollaries}, including 
a proof of Proposition~\ref{Prop-PQsp}.  

Finally, in Section~\ref{sec:asm} the connection is made with alternating sign matrices
and U-turn alternating sign matrices in the $gl(n)$ and $sp(2n)$ cases, respectively.

\section{Background}
\label{sec:background}

\subsection{$gl(n)$ tableaux}\ \
\label{subsec:gl-tableaux}
Let $\lambda=(\lambda_1,\lambda_2,\ldots,\lambda_p)$ 
with $\lambda_1\geq\lambda_2\geq\cdots\geq\lambda_p>0$
be a partition of weight 
$|\lambda|=\lambda_1+\lambda_2+\cdots+\lambda_p$
and length $\ell(\lambda)=p$, where each $\lambda_i$
is a positive integer for all $i=1,2,\ldots,p$. Then 
$\lambda$ defines a Young diagram $F^\lambda$
consisting of $p$ rows of boxes of lengths
$\lambda_1,\lambda_2\ldots,\lambda_p$ left-adjusted
to a vertical line. 

A partition 
$\mu=(\mu_1,\mu_2,\ldots,\mu_q)$ of length $\ell(\mu)=q$
is said to be a strict partition if all the parts 
of $\mu$ are distinct; that is,
$\mu_1>\mu_2>\cdots>\mu_q>0$. A strict partition
$\mu$ defines a shifted Young diagram $SF^\mu$
consisting of $q$ rows of boxes
of lengths $\mu_1,\mu_2,\ldots,\mu_q$
left-adjusted this time to a diagonal line.

For any partition $\lambda$ of length
$\ell(\lambda)\leq n$
let ${\cal T}^\lambda(n)$
be the set of all semistandard tableaux $T$ 
obtained by numbering all the boxes of $F^\lambda$ 
with entries taken from the set $\{1,2,\ldots,n\}$, 
subject to the usual total ordering $1<2<\cdots<n$.
The numbering must be such that the entries are:

\begin{tabular}{rl}
T1&weakly increasing across each row from left to right;\\
T2&strictly increasing down each column from top to bottom.\\
\end{tabular}

\noindent
The weight of the tableau $T$ is given by
$\wgt(T)=(w_1,w_2,\ldots,w_n)$, where 
$w_k$ is the number of times $k$ appears
in $T$ for $k=1,2,\ldots,n$. For example in the case
$n=6$, $\lambda=(3,3,2,1,1)$ we have 
\beq
T=\ 
{\vcenter
{\offinterlineskip
\halign{&\mystrut\vrule#&\mybox{\hss$\scriptstyle#$\hss}\cr
\hr{7}\cr 
&1&&2&&3&\cr 
\hr{7}\cr 
&3&&5&&5&\cr 
\hr{7}\cr 
&4&&6&\cr 
\hr{5}\cr 
&5&\cr
\hr{3}\cr 
&6&\cr 
\hr{3}\cr
\omit& &\omit\cr
\hr{0}\cr
 }}}
\ \in{\cal T}^{33211}(6)
\quad\hbox{with}\quad \wgt(T)=(1,1,2,1,3,2).
\label{Eq-T}
\eeq

By the same token, for any strict partition $\mu$
of length $\ell(\mu)\leq n$,  
let ${\cal ST}^\mu(n)$ be the set of 
all semistandard shifted tableaux $ST$ obtained
by numbering all the boxes of $SF^\mu$ with entries 
taken from the set $\{1,2,\ldots,n\}$,
subject to the total ordering $1<2<\cdots<n$.
The numbering must be such that the entries are:

\begin{tabular}{rl}
ST1&weakly increasing across each row from left to right;\\
ST2&weakly increasing down each column from top to bottom;\\
ST3&strictly increasing down each diagonal from top-left to bottom-right.\\
\end{tabular}

\noindent
The weight of the tableau $ST$ is again given by  
$\wgt(ST)=(w_1,w_2,\ldots,w_n)$, where 
$w_k$ is the number of times $k$ appears
in $ST$ for $k=1,2,\ldots,n$. 

The rules
ST1-ST3 serve to exclude any $2\times2$ blocks
of boxes all containing the same entry, and as a result,
each $ST\in{\cal ST}^\mu(n)$ consists of a
sequence of ribbon strips of boxes containing
identical entries. Any given ribbon strip
may consist of a number of disjoint 
connected components. Let $\str(ST)$
denote the total number of disjoint connected components
of all the ribbon strips. Let $\hgt(ST)$ be the height of the tableaux,
defined $\hgt(ST) = \sum_{k=1}^{n} (\row_k(ST)-\con_k(ST))$, where
$\row_k(ST)$ is the number of rows of $S$ containing an entry $k$,
and $\con_k(ST)$ is the number of connected components of the ribbon
strip of $ST$ consisting of all the entries $k$.

By way of illustration, consider the case $n=6$, $\mu=(9,8,6,4,3,1)$
and the semistandard shifted tableau:
\beq 
ST=\ {\vcenter
{\offinterlineskip
\halign{&\mystrut\vrule#&\mybox{\hss$\scriptstyle#$\hss}\cr
\hr{19}\cr 
&1&&1&&1&&2&&2&&2&&3&&3&&5&\cr
\hr{19}\cr 
\omit&
&&2&&2&&3&&3&&4&&5&&5&&6&\cr
\nr{2}&\hr{17}\cr 
\omit& &\omit&
&&3&&3&&4&&4&&5&&6&\cr 
\nr{4}&\hr{13}\cr 
\omit& &\omit& &\omit&
&&4&&5&&5&&5&\cr
\nr{6}&\hr{9}\cr 
\omit& &\omit& &\omit& &\omit&
&&5&&6&&6&\cr 
\nr{8}&\hr{7}\cr
\omit& &\omit& &\omit& &\omit& &\omit&
&&6&\cr 
\nr{10}&\hr{3}\cr
 }}}\ \in {\cal ST}^{986431}(6)
\quad\hbox{with}\quad
\begin{array}{l}
\wgt(ST)=(3,5,6,4,8,5)\\ \\
\str(ST)=12,\ \ \hgt(ST)=6.\\
\end{array}
\label{Eq-ST}
\eeq

Refining this construct, for any strict partition $\mu$
with $\ell(\mu)\leq n$,
let ${\cal PST}^\mu(n)$ be the set of 
all primed, or marked, semistandard shifted 
tableaux $PST$ obtained
by numbering all the boxes of $SF^\mu$ with entries 
taken from the set $\{1',1,2',2,\ldots,n',n\}$,
subject to the total ordering $1'<1<2'<2<\cdots<n'<n$.
The numbering must be such that the entries are:

\begin{tabular}{rl}
PST1&weakly increasing across each row from left to right;\\
PST2&weakly increasing down each column from top to bottom;\\
PST3&with no two identical unprimed entries in any column;\\
PST4&with no two identical primed entries in any row;\\
PST5&with no primed entries on the main diagonal.\\\
\end{tabular}

\noindent
The weight of the tableau $PST$ is then defined to
be $\wgt(PST)=(u_1,u_2,\ldots,u_n/v_1,v_2,\ldots,v_n)$,
where $u_k$ and $v_k$ are the number of times $k$ and $k'$,
respectively, appear in $PST$ for $k=1,2,\ldots,n$.

The passage from ${\cal ST}^\mu(n)$ to ${\cal PST}^\mu(n)$
is effected merely by adding primes  to the entries of each 
$ST\in{\cal ST}^\mu(n)$ in all possible ways that are
consistent with PST1-5 to give some $PST\in{\cal PST}^\mu(n)$. 
The only entries for which any choice is possible are
those in the lower left hand box at the
beginning of each connected component of a ribbon strip. 
Thereafter, in that connected component of the ribbon
strip, entries in the boxes of its horizontal portions are 
unprimed and those in the boxes of its vertical 
portions are primed. It should be noted that all the 
boxes on the main diagonal are necessarily at the lower
left hand end of a connected component of a ribbon strip, 
but their entries remain unprimed by virtue of PST5. 

To illustrate this let us assign primes to those entries 
of $ST$ in (\ref{Eq-ST}) for which it is essential (that is, 
for every entry lying immediately above the same entry)
and some of those for which it is optional (those entries off the 
main diagonal that are at the start of any continuous strip of equal
entries). This gives, for example,
\beq
PST=\ {\vcenter
{\offinterlineskip
\halign{&\mystrut\vrule#&\mybox{\hss$\scriptstyle#$\hss}\cr
\hr{19}\cr 
&1&&1&&1&&2'&&2&&2&&3&&3&&5&\cr
\hr{19}\cr 
\omit&
&&2&&2&&3'&&3&&4'&&5'&&5&&6'&\cr
\nr{2}&\hr{17}\cr 
\omit& &\omit&
&&3&&3&&4'&&4&&5'&&6&\cr 
\nr{4}&\hr{13}\cr 
\omit& &\omit& &\omit&
&&4&&5'&&5&&5&\cr
\nr{6}&\hr{9}\cr 
\omit& &\omit& &\omit& &\omit&
&&5&&6'&&6&\cr 
\nr{8}&\hr{7}\cr
\omit& &\omit& &\omit& &\omit& &\omit&
&&6&\cr 
\nr{10}&\hr{3}\cr
 }}}\ \in {\cal PST}^{986431}(6)
\quad\hbox{with}\quad
\wgt(PST)=(3,4,5,2,5,3/0,1,1,2,3,2).
\label{Eq-PST}
\eeq

We may replace PST1-4 by identical conditions QST1-4,
but discard PST5. 
This serves to define corresponding primed shifted tableaux
$QST\in{\cal QST}^\mu(n)$ that now involve both primed
and unprimed entries on the main diagonal.

Finally, in this $gl(n)$ context, for fixed positive integer 
$n$, let $\delta=(n,n-1,\ldots,1)$
and let ${\cal PD}^\delta(n)$ be the set of all primed shifted  
tableaux, $PD$, of shape $\delta$,
obtained by numbering the boxes of $SF^\delta$
with entries taken from the set $\{1',1,2',2,\ldots,n',n\}$ 
in such a way that

\begin{tabular}{rl}
PD1&each unprimed entry $k$ appears only in the $k$th row;\\
PD2&each primed entry $k'$ appears only in the $k$th column;\\
PD3&there are no primed entries on the main diagonal.\\
\end{tabular}

\noindent  
The weight of the tableau $PD$ is defined by
$\wgt(PD)=(\u/\v)=(u_1,u_2,\ldots,u_n/v_1,v_2,\ldots,v_n)$,
where $u_k$ and $v_k$ are the numbers of times $k$ and $k'$,
respectively, appear in $PD$ for $k=1,2,\ldots,n$. Typically
for $n=6$ we have
\beq
PD=\ 
{\vcenter
{\offinterlineskip
\halign{&\mystrut\vrule#&\mybox{\hss$\scriptstyle#$\hss}\cr
\hr{13}\cr 
&1&&2'&&1&&4'&&5'&&6'&\cr 
\hr{13}\cr 
\omit& &&2&&3'&&2&&5'&&2&\cr 
\nr{2}&\hr{11}\cr 
\omit& &\omit&
&&3&&4'&&3&&3&\cr 
\nr{4}&\hr{9}\cr 
\omit& &\omit& &\omit&
&&4&&5'&&6'&\cr
\nr{6}&\hr{7}\cr 
\omit& &\omit& &\omit& &\omit&
&&5&&5&\cr 
\nr{8}&\hr{5}\cr
\omit& &\omit& &\omit& &\omit& &\omit&
&&6&\cr 
\nr{10}&\hr{3}\cr
 }}}\ \in{\cal PD}^{654321}(6)
\quad\hbox{with}\quad
\wgt(PD)=(2,3,3,1,2,1/0,1,1,2,3,2).
\label{Eq-PD}
\eeq

Since the $i$th entry on the main diagonal is always $i$
and for $i<j$ the entry in the $(i,j)$th position
is either $i$ or $j'$, it is clear that 
\beq
   \sum_{PD\in{\cal PD}^\delta(n)} (\x/\y)^{\wgt(PD)}
     =\prod_{i=1}^n\ x_i\ \prod_{1\leq i<j\leq n}\ (x_i+y_j).
\label{Eq-PDprod}
\eeq

By way of a small variation of the above, if we replace PD1-2
by identical conditions QD1-2 and discard the condition PD3, 
the corresponding set ${\cal QD}^\delta(n)$ of 
primed shifted tableaux $QD$ differs from ${\cal PD}^\mu(n)$
only in allowing primed entries on the main diagonal. 
It follows that
\beq
   \sum_{QD\in{\cal QD}^\delta(n)} (\x/\y)^{\wgt(QD)}
     =\prod_{1\leq i\leq j\leq n}\ (x_i+y_j).
\label{Eq-QDprod}
\eeq

These formulae (\ref{Eq-PDprod}) and (\ref{Eq-QDprod})
offer a combinatorial interpretation of factors 
appearing in the expansions (\ref{Eq-PQxy}) of
Proposition~\ref{Prop-PQ}. This will be exploited later in 
Section~\ref{sec:gl-bijection}.

\subsection{$sp(2n)$ tableaux}\ \
\label{subsec:sp-tableaux}
In order to establish a similar approach to Proposition~\ref{Prop-PQsp}
it is necessary to extend our already copious list of tableaux to 
encompass certain tableaux associated with the symplectic algebra $sp(2n)$. 
As before it is helpful to start with definitions of the various types of tableaux, both 
shifted and unshifted. 

For any partition $\lambda$ of length
$\ell(\lambda)\leq n$,
let ${\cal T}^\lambda(n,\ov{n})$
be the set of all semistandard symplectic tableaux $T$ 
obtained by numbering all the boxes of 
$F^\lambda$ with entries from the set $\{\ov1,1,\ov2,2,\ldots, \ov{n},n\}$, 
subject to the usual total ordering $\ov1<1<\ov2<2<\cdots\ov{n}<n$.
The entries are:

\begin{tabular}{rl}
T1&weakly increasing across each row from left to right;\\
T2&strictly increasing down each column from top to bottom.\\
T${\ov3}$&$k$ or $\ov{k}$ may appear no lower than the $k$th row.\\
\end{tabular}

\noindent
The weight of the symplectic tableau $T$ is given by
$\wgt(T)=(\w)=(w_1,w_2,\ldots,w_n)$, with $w_k=n_k-n_{\ov{k}}$ where 
$n_k$ and $n_{\ov{k}}$ are the number of times $k$ and $\ov{k}$,
respectively, appear in $T$ for $k=1,2,\ldots,n$. 
The parameter $bar(T)$ is equal to the number of barred entries in the tableau.
For example in the case
$n=5$, $\lambda=(4,3,3)$ we have 
\beq
T=\ 
{\vcenter
{\offinterlineskip
\halign{&\mystrut\vrule#&\mybox{\hss$\scriptstyle#$\hss}\cr
\hr{9}\cr 
&\ov1&&1&&\ov2&&4&\cr 
\hr{9}\cr 
&3&&\ov4&&\ov4&\cr 
\hr{7}\cr 
&\ov4&&4&&4&\cr 
\hr{7}\cr 
 }}}
\ \in{\cal T}^{433}(5,5)
\quad\hbox{with}\quad
\begin{array}{l}
\wgt(T)=(0,-1,1,0,0)\\ \\
\bar(T)= 5.\\
\end{array} 
\label{Eq-Tsp}
\eeq

For any strict partition $\mu$
of length $\ell(\mu)\leq n$, 
let ${\cal ST}^\mu(n,\ov{n})$ be the set of 
all semistandard shifted symplectic tableaux $ST$ obtained
by numbering all the boxes of $SF^\mu$ with entries 
taken from the set $\{\ov1,1,\ov2,2,\ldots,\ov{n},n\}$,
subject to the total ordering $\ov1<1<\ov2<2<\cdots<\ov{n}<n$.
The numbering must be such that the entries are:

\begin{tabular}{rl}
ST1&weakly increasing across each row from left to right;\\
ST2&weakly increasing down each column from top to bottom;\\
ST3&strictly increasing down each diagonal from top-left to bottom-right;\\
ST${\ov4}$&with $d_k\in\{k,\ov{k}\}$, 
where $d_k$ is the $k$th entry on the main diagonal.\\
\end{tabular}

\noindent
The weight of the shifted symplectic tableau $ST$ is given by
$\wgt(ST)=(w_1,w_2,\ldots,w_n)$, with $w_k=n_k-n_{\ov{k}}$ where 
$n_k$ and $n_{\ov{k}}$ are the number of times $k$ and $\ov{k}$,
respectively, appear in $ST$ for $k=1,2,\ldots,n$. Once again it
is convenient, following~\cite{HK03},
to introduce $\str(ST)$ as the total number of
disjoint connected components of all ribbon strips of $ST$, and
$\var(ST)=\sum_{k=1}^n(\row_k(ST)-\con_k(ST)+\col_{\ov{k}}(ST)-\con_{\ov{k}}(ST))$,
where $\row_k(ST)$ is the number of rows of $ST$ containing an entry $k$,
$\col_{\ov{k}}(ST)$ is the number of columns containing an entry $\ov{k}$,
while $\con_k(ST)$ and $\con_{\ov{k}}(ST)$ are the number
of connected components of the ribbon strips of $ST$ consisting of 
all the entries $k$ and $\ov{k}$, respectively, 
and $\bar(ST)$ is equal to the total number of  barred  entries. 

Typically, for $n=5$ and $\mu=(9,7,6,2,1)$ we have
\begin{equation}
ST=\
{\vcenter
 {\offinterlineskip
 \halign{&\mystrut\vrule#&\mybox{\hss$\scriptstyle#$\hss}\cr
  \hr{19}\cr
  &\ov1&&1&&\ov2&&2&&\ov3&&\ov3&&\ov4&&4&&5&\cr
  \hr{19}\cr
  \omit& &&\ov2&&\ov2&&2&&3&&\ov4&&\ov4&&4&\cr
  \nr{2}&\hr{15}\cr
  \omit& &\omit& &&3&&\ov4&&4&&4&&4&&4&\cr
  \nr{4}&\hr{13}\cr
  \omit& &\omit& &\omit& &&4&&4&\cr
  \nr{6}&\hr{5}\cr
  \omit& &\omit& &\omit& &\omit& &&\ov5&\cr
  \nr{8}&\hr{3}\cr
 }}}
\ \in {\cal ST}^{97621}(5,\ov5)
\quad\hbox{with}\quad
\begin{array}{l}
\wgt(ST)=(0,-1,0,4,0)\\ \\
\bar(ST)= 11,\ \ \str(ST)=12,\ \ \var(ST)=7.\\
\end{array}
\label{Eq-STsp}
\end{equation}

Refining this construct, for any strict partition $\mu$
with $\ell(\mu)\leq n$,
let ${\cal PST}^\mu(n,\ov{n})$ be the set of 
all primed semistandard shifted symplectic
tableaux $PST$ obtained
by numbering all the boxes of $SF^\mu$ with entries 
taken from the set 
$\{\ov{1}',\ov{1},1',1,{\ov 2}',{\ov 2},2',2,\ldots,\ov{n},\ov{n},n',n\}$,
subject to the total ordering 
\begin{equation}
\ov1'<\ov1<1'<1<\ov2'<\ov2<2'<2<\cdots<{\ov n}'<{\ov n}<{n}'<{n}.
\label{Eq-order}
\end{equation}
The numbering must be such that the entries are:

\begin{tabular}{rl}
PST1&weakly increasing across each row from left to right;\\
PST2&weakly increasing down each column from top to bottom;\\
PST3&with no two identical unprimed entries in any column;\\
PST4&with no two identical primed entries in any row;\\
PST${\ov5}$&with $d_k\in\{\ov{k},k\}$,  
where $d_k$ is the $k$th entry on the main diagonal.\\
\end{tabular}

\noindent
The weight of the tableau $PST$ is then defined to
be $\wgt(PST)=(\u/\v)$ with
$\u=(u_1,u_2,\ldots,u_n)$ and $\v=(v_1,v_2,\ldots,v_n)$,
where $u_k=n_k-n_{\ov{k}}$ and $v_k=n_{k'}-n_{\ov{k}'}$, 
with $n_k$, $n_{\ov{k}}$, $n_{k'}$ and $n_{\ov{k}'}$
are the number of times $k$, $\ov{k}$, $k'$ and $\ov{k}'$,
respectively, appear in $PST$ for $k=1,2,\ldots,n$.
In addition, let $\bar(PST)$ be the total number
of barred entries in $PST$.

If we now replace PST1-4 by identical conditions QST1-4 and
replace PST${\ov5}$ by:

\begin{tabular}{rl}
QST${\ov5}$&with $d_k\in\{\ov{k}',\ov{k},k',k\}$,  
where $d_k$ is the $k$th entry on the main diagonal.\\
\end{tabular} 

\noindent
Then once again the corresponding primed shifted tableaux 
$QST\in{\cal QST}^\mu(n,{\ov n})$ now have
primes allowed on the main diagonal.

Typically, for $n=5$ and $\mu=(9,7,6,2,1)$ we have
\begin{equation}
 QST=\ 
{\vcenter
 {\offinterlineskip
 \halign{&\mystrut\vrule#&\mybox{\hss$\scriptstyle#$\hss}\cr
  \hr{19}\cr
  &\ov1&&1&&\ov2'&&2'&&\ov3'&&\ov3&&\ov4'&&4'&&5&\cr
  \hr{19}\cr
  \omit& &&\ov2'&&\ov2&&2&&3&&\ov4'&&\ov4&&4'&\cr
  \nr{2}&\hr{15}\cr
  \omit& &\omit& &&3'&&\ov4'&&4'&&4&&4&&4&\cr
  \nr{4}&\hr{13}\cr
  \omit& &\omit& &\omit& &&4'&&4&\cr
  \nr{6}&\hr{5}\cr
  \omit& &\omit& &\omit& &\omit& &&\ov5'&\cr
  \nr{8}&\hr{3}\cr
 }}}
\ \in {\cal QST}^{97621}(5,\ov5) 
\quad\hbox{with}\quad
\begin{array}{l}
\wgt(QST)=(0,0,0,3,1/0,-1,0,1,-1),\\ \\
\bar(QST)=11.\\
\end{array}
\label{Eq-QSTsp}
\end{equation}

To complete our set of $sp(2n)$ tableaux, for fixed positive integer 
$n$, let $\delta=(n,n-1,\ldots,1)$
and let ${\cal PD}^\delta(n,{\ov n})$ be the set of all primed shifted  
tableaux, $PD$, of shape $\delta$,
obtained by numbering the boxes of $SF^\delta$
with entries taken from the set 
$\{{\ov1}',\ov1,1',1,\ov2',\ov2,2',2,\ldots,{\ov n}',{\ov n},n',n\}$ 
in such a way that

\begin{tabular}{rl}
PD${\ov1}$&each unprimed entry $k$ or ${\ov k}$ appears only in the $k$th row;\\
PD${\ov2}$&each primed entry $k'$ or ${\ov k}'$ appears only in the $k$th column;\\
PD3&there are no primed entries on the main diagonal.\\
\end{tabular}

\noindent  
The weight of the tableau $PD$ is defined by
$\wgt(PD)=(\u/\v)$ with
$\u=(u_1,u_2,\ldots,u_n)$ and $\v=(v_1,v_2,\ldots,v_n)$,
where $u_k=n_k-n_{\ov{k}}$ and $v_k=n_{k'}-n_{\ov{k}'}$, 
with $n_k$, $n_{\ov{k}}$, $n_{k'}$ and $n_{\ov{k}'}$
are the number of times $k$, $\ov{k}$, $k'$ and $\ov{k}'$,
respectively, appear in $PD$ for $k=1,2,\ldots,n$.
In addition let $\bar(PD)$ be the total number
of barred entries in $PD$.

With this notation, 
since the entry in the $i$th position on the main diagonal
is either $i$ or ${\ov i}$
while for $i<j$ the entry in the $(i,j)$th position
is either $i$, ${\ov i}$, $j'$ or ${\ov j}'$, it is clear that 
\beq
   \sum_{PD\in{\cal PD}^\delta(n,{\ov n})} 
   t^{2\bar(PD)}\ (\x/\y)^{\wgt(PD)}
     =\prod_{i=1}^n\ (x_i+t^2{\ov x}_i)\ 
   \prod_{1\leq i<j\leq n}\ (x_i+t^2{\ov x}_i+y_j+t^2{\ov y}_j).
\label{Eq-PDprodbar}
\eeq

By way of a small variation of the above, if we replace PD1-2
by identical conditions QD1-2 and discard the condition PD3, 
the corresponding set ${\cal QD}^\delta(n)$ of 
primed shifted tableaux $QD$ differs from ${\cal PD}^\mu(n)$
only in allowing primed entries on the main diagonal. 

Typically
for $n=5$ we have
\beq
QD=\ 
{\vcenter
{\offinterlineskip
\halign{&\mystrut\vrule#&\mybox{\hss$\scriptstyle#$\hss}\cr
\hr{11}\cr 
&\ov1&&1&&\ov3'&&4'&&\ov1&\cr 
\hr{11}\cr 
\omit& &&2&&2&&2&&\ov2&\cr 
\nr{2}&\hr{9}\cr 
\omit& &\omit&
&&\ov3'&&3&&\ov3&\cr 
\nr{4}&\hr{7}\cr 
\omit& &\omit& &\omit&
&&4'&&4&\cr
\nr{6}&\hr{5}\cr 
\omit& &\omit& &\omit& &\omit&
&&\ov5'&\cr 
\nr{8}&\hr{3}\cr
 }}}\ \in{\cal QD}^{54321}(5,\ov5)
\quad\hbox{with}\quad
\begin{array}{l}
\wgt(QD)=(-1,2,0,1,0/0,0,-2,2,-1),\\ \\
\bar(QD)=7.\\
\end{array}
\label{Eq-QDbar}
\eeq

It follows from our definition of ${\cal QD}(n,{\ov n})$ that
\beq
   \sum_{QD\in{\cal QD}^\delta(n,{\ov n})}\
    t^{2\bar(QD)}\ (\x/\y)^{\wgt(QD)}
     =\prod_{1\leq i\leq j\leq n}\ (x_i+t^2{\ov x}_i+y_j+t^2{\ov y}_j).
\label{Eq-QDprodbar}
\eeq

These formulae (\ref{Eq-PDprodbar}) and (\ref{Eq-QDprodbar})
have been introduced so as to offer a combinatorial interpretation 
of factors appearing in the expansions (\ref{Eq-PQxybar}) of
Proposition~\ref{Prop-PQsp}. This will be exploited later in 
Section~\ref{sec:gl-bijection}.

\subsection{Schur's $P$ and $Q$ functions and their generalisations}
\label{subsec:PQ}\

Let $\x=(x_1,x_2,\ldots,x_n)$ be a vector of $n$ indeterminates
and let $\w=(w_1,w_2,\ldots,w_n)$ be a vector of $n$ non-negative 
integers. Then $\x^\w=x_1^{w_1}x_2^{w_2}\cdots x_n^{w_n}$.
With this notation it is well known that each partition 
$\lambda$ of length $\ell(\lambda)\leq n$ specifies a Schur 
function $s_\lambda(\x)$ with combinatorial definition:

\beq
     s_\lambda(\x) = \sum_{T\in{\cal T}^\lambda(n)}  \x^{\wgt(T)}
\label{Eq-Sf}
\eeq 

Similarly, each strict partition $\mu$ of length
$\ell(\mu)\leq n$ specifies a Schur $P$-function
and a Schur $Q$-function whose combinatorial definitions take the form:
\beq
\begin{array}{rcl}
P_\mu(\x) &=& \sum_{ST\in{\cal ST}^\mu(n)} 
       2^{\str(ST)-\ell(\mu)} \x^{\wgt(ST)};\\ \\
Q_\mu(\x) &=& \sum_{ST\in{\cal ST}^\mu(n)} 
 2^{\str(ST)} \x^{\wgt(ST)}. \\
\end{array}
\label{Eq-PQx}
\eeq

Now let $\z=(\x/\y)=(x_1,x_2,\ldots,x_n/y_1,y_2,\ldots,y_n)$,
where $\x$ and $\y$ are two vectors of $n$ indeterminates, and
let $\w=(\u/\v)=(u_1,u_2,\ldots,u_n/v_1,v_2,\ldots,v_n)$
where $\u$ and $\v$ are two vectors of $n$ non-negative integers.
Then let $\z^\w=(\x/\y)^{(\u/\v)}=\x^\u\,\y^\v
=x_1^{u_1}\cdots x_n^{u_n}\,y_1^{v_1}\cdots y_n^{v_n}$.
With this notation each strict partition $\mu$ of length
$\ell(\mu)\leq n$ serves to specify generalised Schur 
$P$ and $Q$-functions defined by:

\beq
\begin{array}{rcl}
   P_\mu(\x/\y)&=&\sum_{PST\in{\cal PST}^\mu(n)} (\x/\y)^{\wgt(PST)};\\ \\
   Q_\mu(\x/\y)&=&\sum_{QST\in{\cal QST}^\mu(n)} (\x/\y)^{\wgt(QST)}.\\
\end{array}
\label{Eq-PQxy}
\eeq

Since the maps back from $PST\in{\cal PST}^\mu(n)$ 
and from $QST\in{\cal QST}^\mu(n)$ to some 
$ST\in{\cal ST}^\mu(n)$ are effected merely by deleting primes,
and there are no primes on the main diagonal in the 
case of $PST$, it follows that

\beq
     Q_\mu(\x)=2^{\ell(\mu)}P_\mu(\x) \quad\hbox{with}\quad
            P_\mu(\x)=P_\mu(\x/\x) \quad\hbox{and}\quad
            Q_\mu(\x)=Q_\mu(\x/\x) 
\label{Eq-QxPx}
\eeq

It might be noted that $s_\lambda(\x)$, $P_\lambda(\x)$ and $Q_\lambda(\x)$
are nothing other than the specialisations $P_\lambda(\x;1)$,
$P_\mu(\x;-1)$ and $Q_\mu(\x;-1)$, respectively, of the
Hall-Littlewood functions $P_\mu(\x;t)$ and $Q_\mu(\x;t)$.

Turning to the symplectic case,
it is well known that each partition 
$\lambda$ of length $\ell(\lambda)\leq n$ specifies an
irreducible representation of $sp(2n)$ whose character
$sp_\lambda(\x)$ may be given a combinatorial definition:
\beq
     sp_\lambda(\x) = \sum_{T\in{\cal T}^\lambda(n,{\ov n})}  \x^{\wgt(T)}.
\label{Eq-Sp}
\eeq 
This may be $t$-deformed to give 
\beq
    sp_\lambda(\x;t) = \sum_{T\in{\cal T}^\lambda(n,{\ov n})}\ 
    t^{2\bar(T)}\ \x^{\wgt(T)}.
\label{Eq-Spt}
\eeq

In the case of a strict partition $\mu$ of length $\ell(\mu)=n$
the required generalisations of Schur $P$ and $Q$ functions
take the form:
\beq
\begin{array}{rcl}
     P_\mu(\x/\y;t) &=& \sum_{PST\in{\cal PST}^\mu(n,{\ov n})} \
     t^{2\bar(PST)}\ (\x/\y)^{\wgt(PST)}; \\ \\
     Q_\mu(\x/\y;t) &=& \sum_{QST\in{\cal QST}^\mu(n,{\ov n})} \
     t^{2\bar(QST)}\ (\x/\y)^{\wgt(QST)}. \\ 
\end{array}
\label{Eq-PQxybar}
\eeq

\section{The $gl(n)$ bijection}
\label{sec:gl-bijection}
\subsection{Main Result}

The generalisations of the combinatorial definitions of $P_\mu(\x)$
and $Q_\mu(\x)$ to $P_\mu(\x/\y)$ and $Q_\mu(\x/\y)$, respectively,
together with those of $s_\lambda(\x)$ and the product factors appearing 
in (\ref{Eq-PQxy}), allow us to establish the validity of 
Proposition~\ref{Prop-PQ} by first proving the following:

\begin{Theorem}
Let $\mu=\lambda+\delta$ be a strict partition of length $\ell(\mu)=n$,
with $\lambda$ a partition of length $\ell(\lambda)\leq n$ and 
$\delta=(n,n-1,\ldots,1)$. There exists a weight preserving, bijective map 
$\Theta$ from ${\cal PST}^\mu(n)$ to 
$\big({\cal PD}^\delta(n),{\cal T}^\lambda(n)\big)$ and
from ${\cal QST}^\mu(n)$ to 
$\big({\cal QD}^\delta(n),{\cal T}^\lambda(n)\big)$
such that for all $PST\in{\cal PST}^\mu(n)$ and 
for all $QST\in{\cal QST}^\mu(n)$
\beq
\begin{array}{rcl}
   \Theta:\ PST&\ \mapsto\ &(PD,T)  \quad\hbox{with}\quad
     \wgt(PST)=\wgt(PD)+\wgt(T). \\ \\
   \Theta:\ QST&\ \mapsto\ &(QD,T)  \quad\hbox{with}\quad
     \wgt(QST)=\wgt(QD)+\wgt(T).\\
\end{array}
\label{Eq-Theta}
\eeq
with $PD\in{\cal PD}^\delta(n)$, $QD\in{\cal QD}^\delta(n)$ 
and $T\in{\cal T}^\lambda(n)$.
\label{The-PQ}
\end{Theorem}

\noindent{\bf Proof}\ \

We choose to tackle the $PST$ case first with
the aim of describing a candidate map $\Theta$ and 
showing that it is both weight preserving and bijective.

The technique is to apply the jeu de taquin (\cite{S87},
\cite{W84}, \cite{SS90}, \cite{M95},\cite{HH92}) to the primed
entries $k'$ of $PST$ taken in turn starting with any $1'$s 
(actually there are none), then any $2'$s (at most one), then
any $3'$s (at most two) and so on. If for fixed $k$ there is more 
than one $k'$ in $PST$ then these are dealt with in turn from top to 
bottom before a final rearrangement is made of the entries in the
$k$th column. The map $\Theta$ is thus expressible in the form
$\Theta=\theta_{n'}\circ\cdots\circ\theta_{2'}\circ\theta_{1'}$.

We start by describing the map $\theta_{k'}$. This involves sliding 
each $k'$ in the north-west
direction by a sequence of interchanges with either its northern or
western neighbour until it reaches a position in the $k$th column
either in the topmost row, or immediately below another $k'$,
or immediately below some unprimed entry $i$ in the $i$th row.

This amounts to 
playing jeu de taquin, treating $k'$ to be strictly less than 
all the unprimed entries. At every stage all the unprimed entries 
must satisfy the semistandardness conditions T1 and T2; that 
is, they should be weakly increasing across rows and strictly 
increasing down columns. It is this that ensures that each move 
made by $k'$ is uniquely determined. Consider first the situation
illustrated by the tableau $T_0$ in (\ref{Eq-jeu-NW}) with
$k'$ not yet in the $k$th column, nor in the top row. This is to be 
thought of as the subtableau surrounding a particular $k'$ awaiting
its next move. For the time being we assume that $b,d$ are unprimed,
while $a,c,e,f,g,h$ may be primed, or unprimed, or even absent if
$k'$ is at or near the southern or eastern edge of the complete diagram. 
However, all the unprimed entries amongst $a,b,\ldots,h$ must, by hypothesis,
satisfy the semistandardness conditions T1 and T2.

Now for the jeu de taquin rules that define the map $\theta_{k'}$.
If $d\leq b$ then  $k'$ is to be interchanged
with $b$ and if $d>b$ then $k'$ is to be interchanged with
$d$. In the first case $k'$ moves north and the resulting tableau $T_N$
satisfies  both  $T1$  and $T2$ since $d\leq b\leq c<e$, while in the
second, $k'$ moves west and the resulting tableau $T_W$
satisfies  both  $T1$  and $T2$ since $b<d<f\leq g$.
\beq
\theta_{k'}:\ \
T_0=\ {\vcenter
{\offinterlineskip
\halign{&\mystrut\vrule#&\mybox{\hss$\scriptstyle#$\hss}\cr
\hr{7}\cr 
&a&&b&&c&\cr
\hr{7}\cr 
&d&&k'&&e&\cr
\hr{7}\cr
&f&&g&&h&\cr 
\hr{7}\cr
 }}}
\quad
\ \ \longrightarrow\ \
\quad
\left\{  
\begin{array}{cc}
T_N=\ {\vcenter
{\offinterlineskip
\halign{&\mystrut\vrule#&\mybox{\hss$\scriptstyle#$\hss}\cr
\hr{7}\cr 
&a&&k'&&c&\cr
\hr{7}\cr 
&d&&b&&e&\cr
\hr{7}\cr
&f&&g&&h&\cr 
\hr{7}\cr
 }}}
&\quad\hbox{if}\quad b\geq d;\\  \\
T_W=\ {\vcenter
{\offinterlineskip
\halign{&\mystrut\vrule#&\mybox{\hss$\scriptstyle#$\hss}\cr
\hr{7}\cr 
&a&&b&&c&\cr
\hr{7}\cr 
&k'&&d&&e&\cr
\hr{7}\cr
&f&&g&&h&\cr 
\hr{7}\cr
 }}}
&\quad\hbox{if}\quad b<d.\\
\end{array} 
\right.
\label{Eq-jeu-NW}
\eeq
If $k'$ is already in the topmost row, so that the row $a\,b\,c$ is
absent then $\theta_{k'}$
acts on $T_0$ as follows:
\beq
\theta_{k'}:\ \
T_0=\ {\vcenter
{\offinterlineskip
\halign{&\mystrut\vrule#&\mybox{\hss$\scriptstyle#$\hss}\cr
\hr{7}\cr 
&d&&k'&&e&\cr
\hr{7}\cr
&f&&g&&h&\cr 
\hr{7}\cr
 }}}
\quad
\ \ \longrightarrow \ \
\quad
T_W=\ {\vcenter
{\offinterlineskip
\halign{&\mystrut\vrule#&\mybox{\hss$\scriptstyle#$\hss}\cr
\hr{7}\cr 
&k'&&d&&e&\cr
\hr{7}\cr
&f&&g&&h&\cr 
\hr{7}\cr
 }}}
\label{Eq-jeu-W}
\eeq
Once again the unprimed entries in $T_W$ satisfy both T1 and T2,
since we still have $d<f\leq g$.

On the other hand, if $k'$ is already in the $k$th column,
so that the column $a\,d\,f$ is absent, the map $\theta_{k'}$ 
leaves $T_0$ unaltered, that is $k'$ has reached its final resting 
place, unless $k'$ lies in the $i$th row with an 
unprimed entry $b\geq i$ immediately above it. In such a case
$\theta_{k'}$ acts on $T_0$ as shown below
\beq
\theta_{k'}:\ \
\begin{array}{c}
T_0=\ {\vcenter
{\offinterlineskip
\halign{&\mystrut\vrule#&\mybox{\hss$\scriptstyle#$\hss}\cr
\hr{5}\cr 
&b&&c&\cr
\hr{5}\cr 
&k'&&e&\cr
\hr{5}\cr
&g&&h&\cr 
\hr{5}\cr
 }}}\\
\hbox{$k'$ in $i$th row and $k$th column}\\
b\geq i\\
\end{array}
\quad
\ \ \longrightarrow\ \
\quad
T_N=\ {\vcenter
{\offinterlineskip
\halign{&\mystrut\vrule#&\mybox{\hss$\scriptstyle#$\hss}\cr
\hr{5}\cr 
&k'&&c&\cr
\hr{5}\cr 
&b&&e&\cr
\hr{5}\cr
&g&&h&\cr 
\hr{5}\cr
 }}}
\label{Eq-jeu-N}
\eeq
 Yet again, the unprimed entries of $T_N$ satisfy both T1 and T2,
since we still have $b\leq c<e$.

 Now we return to the two possibilities that we had previously
set aside, namely $b=k'$ or $d=k'$. The first of these cannot
occur in the tableau $T_0$ of (\ref{Eq-jeu-NW}), since in such a case 
the uppermost $b=k'$ would have been moved either north or west before
attempting to move the central $k'$. In the case of the tableau 
$T_0$ of (\ref{Eq-jeu-N}), as we have alreay pointed out, if $b=k'$
then no further move of the $k'$ below $b=k'$ is required.

It follows that the only possible impediment to the movement of $k'$
in a north-westerly direction until it actually reaches the $k$th
column, is the existence of another $k'$ to its immediate left, that is 
in $T_0$ of (\ref{Eq-jeu-NW}) or (\ref{Eq-jeu-N}), we have $d=k'$. 
That this cannot occur is a 
corollary of the fact that the path followed by $k'$ always remains 
column by column below (that is strictly south of) the path followed 
by any preceding $k'$. To see this consider $k'$ arriving, 
 as shown below in the diagram on the left of (\ref{Eq-path-S}), at a
position due south of an entry $b$ which itself lies on the path
of the preceding $k'$. 
\beq
\theta_{k'}:\ \ 
{\vcenter
{\offinterlineskip
\halign{&\mystrut\vrule#&\mybox{\hss$\scriptstyle#$\hss}\cr
\hr{21}\cr 
     &\p&&      &\omit&\cdots&\omit&  &&\a&&\b&&  &\omit&\cdots&\omit& &&\q&\cr
\hr{21}\cr 
     &r&&      &\omit&\cdots&\omit&  &&c&&k'&&  \cr
\hr{13}\cr
}}}
\quad
\begin{array}{c}
b<c\\
\ \ \longrightarrow\ \ \\
\end{array}
\quad
{\vcenter
{\offinterlineskip
\halign{&\mystrut\vrule#&\mybox{\hss$\scriptstyle#$\hss}\cr
\hr{21}\cr 
     &\p&&      &\omit&\cdots&\omit&  &&\a&&\b&&  &\omit&\cdots&\omit& &&\q&\cr
\hr{21}\cr 
     &r&&      &\omit&\cdots&\omit&  &&k'&&c&&  \cr
\hr{13}\cr
}}}
\label{Eq-path-S}
\eeq
If this path of the preceding $k'$ moves north from $b$, then
there is no problem since the $k'$ can follow the same path north or
move west without violating the strictly south condition. On the other 
hand the path of the preceding $k'$ may move west along the indicated
boldface track from $q$ to $p$. In doing so, it must at one stage have
displaced $b$ from its original position at the site of $a$,
immediately above $c$ and satisfying the T2 condition $b<c$.
This condition then ensures that the $k'$ must itself move west 
as shown in (\ref{Eq-path-S}). It therefore stays south of the path of
the preceding $k'$ that passes through the position of $a$. This
implies that the path of $k'$ must always stay strictly south of the
path of the preceding $k'$, thereby excluding the possibility $d=k'$ in both 
(\ref{Eq-jeu-NW}) and (\ref{Eq-jeu-W}). This ensures that each $k'$ will
eventually reach the $k$th column by means of a sequence of moves
of type (\ref{Eq-jeu-NW})-(\ref{Eq-jeu-N}).

Following the action of $\theta_{k'}$ 
 the unprimed entry $k$ on the main diagonal of $PST$
remains fixed, and all $k'$s are in the $k$th column along with 
distinct unprimed entries $j$ with $1\leq j\leq k$.
If $k'$ appears in the $i$th row, then $i$ cannot appear above it,
since $k'$ would then move north  as in (\ref{Eq-jeu-N}),
and cannot appear below it, since it would then be to the right of
a diagonal entry greater than $i$ and thus violate the
weakly increasing condition T1. It follows that the unprimed 
entries $j$ in the $k$th column do not include the row 
numbers of $k'$. Since they are distinct and $1\leq j\leq k$, 
they must include all the other row numbers, and be arranged 
in strictly increasing order in accordance with T2.
This means that each unprimed entry 
in the $k$th column lies in its own row. Since the primed entries
in this column are all $k'$s all the
entries in the $k$th column satisfy PD1-3.

Iterating this procedure for all $k=1,2,\ldots,n$ results in all 
primed entries being moved to the first $k$ columns of $SF^\mu$
along with some unprimed entries, collectively satisfying
PD1-3 in this region of shape $SF^\delta$, and leaving only unprimed 
entries, all satisfying T1-2, in the right hand region of shape 
$F^\lambda$. That is, the result of applying $\Theta$ to 
$PST\in{\cal PST}^\mu$ is a semistandard tableau $T\in{\cal T}^\lambda(n)$
of shape $\lambda$ juxtaposed to a primed tableau $PD\in{\cal PD}^\delta(n)$.
This map is necessarily weight preserving since every individual step is a 
simple interchange which does not alter the number of $k$s or $k'$s for any $k$.

To show that this map $\Theta$ is bijective it should be noted that each step
may be reversed. 
One starts by juxtaposing an arbitrary pair of tableaux $PD\in{\cal PD}^\delta(n)$
and $T\in{\cal T}^\lambda(n)$ to create a tableaux of shape $SF^\mu$ with
$\mu=\lambda+\delta$. Then for each $k$ taken in turn from 
$n$ to $n-1$ down to $1$ one applies $\theta_{k'}^{-1}$
to all the primed entries $k'$; that is to say one reverses 
the action of $\theta_{k'}$ 
by playing jeu de taquin in the reverse direction with primed entries $k'$
treated in turn from bottom to top, moving each one in a south easterly direction 
with $k'$ now assumed to be larger than $i$ for $i=1,\ldots, k-1$ but less than 
$j$ for $j=k,k+1,\ldots,n$ 
with the semistandardness conditions T1 and T2 applying to all unprimed entries 
at all times. For example, in the following diagram this leads 
unambiguously from $T_0$ to $T_E$ if $e<g$ and from $T_0$ to $T_S$ 
if $e\geq g$, with all unprimed entries satisfying the semistandardness 
conditions T1 and T2 since in addition $b\leq c<e<g$ in $T_E$ and
$d<f\leq g\leq e$ in $T_S$.
\beq
\theta^{-1}_{k'}:\ \
T_0=\ {\vcenter
{\offinterlineskip
\halign{&\mystrut\vrule#&\mybox{\hss$\scriptstyle#$\hss}\cr
\hr{7}\cr 
&a&&b&&c&\cr
\hr{7}\cr 
&d&&k'&&e&\cr
\hr{7}\cr
&f&&g&&h&\cr 
\hr{7}\cr
 }}}
\quad
\ \ \longrightarrow\ \  
\quad
\left\{ 
\begin{array}{ll} 
T_E=\ {\vcenter
{\offinterlineskip
\halign{&\mystrut\vrule#&\mybox{\hss$\scriptstyle#$\hss}\cr
\hr{7}\cr 
&a&&b&&c&\cr
\hr{7}\cr 
&d&&e&&k'&\cr
\hr{7}\cr
&f&&g&&h&\cr 
\hr{7}\cr
 }}}
&\quad\hbox{if}\quad e<g;
\\ \\
T_S=\ {\vcenter
{\offinterlineskip
\halign{&\mystrut\vrule#&\mybox{\hss$\scriptstyle#$\hss}\cr
\hr{7}\cr 
&a&&b&&c&\cr
\hr{7}\cr 
&d&&g&&e&\cr
\hr{7}\cr
&f&&k'&&h&\cr 
\hr{7}\cr
 }}}
&\quad\hbox{if}\quad e\geq g.
\end{array}
\right.
\label{Eq-jeu-ES}
\eeq
If $k'$ is already in the lowest row or rightmost
column, the following diagrams illustrate the allowed moves of $k'$ 
east and south respectively.
\beq
\theta^{-1}_{k'}:\ \
T_0=\ {\vcenter
{\offinterlineskip
\halign{&\mystrut\vrule#&\mybox{\hss$\scriptstyle#$\hss}\cr
\hr{7}\cr 
&a&&b&&c&\cr
\hr{7}\cr
&d&&k'&&e&\cr 
\hr{7}\cr
 }}}
\ \ \ \longrightarrow\ \ \
T_E=\ {\vcenter
{\offinterlineskip
\halign{&\mystrut\vrule#&\mybox{\hss$\scriptstyle#$\hss}\cr
\hr{7}\cr 
&a&&b&&c&\cr
\hr{7}\cr
&d&&e&&k'&\cr 
\hr{7}\cr
 }}}
\label{Eq-jeu-E}
\eeq
and
\beq
\theta^{-1}_{k'}:\ \
T_0=\ {\vcenter
{\offinterlineskip
\halign{&\mystrut\vrule#&\mybox{\hss$\scriptstyle#$\hss}\cr
\hr{5}\cr 
&a&&b&\cr
\hr{5}\cr 
&d&&k'&\cr
\hr{5}\cr
&f&&g&\cr 
\hr{5}\cr
 }}}
\ \ \ \longrightarrow\ \
T_S=\ {\vcenter
{\offinterlineskip
\halign{&\mystrut\vrule#&\mybox{\hss$\scriptstyle#$\hss}\cr
\hr{5}\cr 
&a&&b&\cr
\hr{5}\cr 
&d&&g&\cr
\hr{5}\cr
&f&&k'&\cr 
\hr{5}\cr
 }}}
\label{Eq-jeu-S}
\eeq
Once again since the unprimed entries of $T_0$ satisfy both T1 and T2,
these rules are also satisfied by both $T_E$ and $T_S$ 
since we still have $b\leq c<e$ and $d<f\leq g$.

More importantly, returning to the original jeu de taquin moves
illustrated in (\ref{Eq-jeu-NW}), this reversed jeu de taquin is such that if
the conditions T1 and T2 are satisfied by the entries in $T_0$, then 
the reverse process leads directly from $T_W$ to $T_0$ since $d<f$
and from $T_N$ to $T_0$ since $b\leq c$. Thus the original steps 
along each of the $k'$ paths are retraced precisely. The same is true
of the maps (\ref{Eq-jeu-W}) and (\ref{Eq-jeu-N}).

The only task remaining is to show that the
endpoints of these retraced paths results in an element $PST$  of
${\cal PST}^{\mu}$. If $T_0$ is such that $k'$ has reached its
endpoint then neither $e$ nor $g$ is $<k'$. If either $e$ or $g$ 
is $\geq k$ then this poses no problem. If $g\geq  k'$ there is again
no problem since the rules PST1-5 allow two (or more) $k'$s in the same
column. It is only the case $e=k'$ that produces a violation, in  this
case of PST4.  Fortunately this case is excluded by  the following 
argument analogous
to that which led to the strictly south property of the original jeu
de taquin. Now we require a strictly north property.  The argument
goes as follows.  The fact that the strictly north property applies to the reverse jeu de
taquin follows from a consideration of the following diagram in which
a $k'$ south westerly  path meets a preceding west--east path,
indicated by means of boldface entries,  passing from
$p$ to $q$ through the positions of $a$ and $b$. The existence of the 
latter requires that $a$ must initially have been immediately
south of $c$, so that $c<a$. This in turn implies that $k'$ moves eastwards
staying strictly north of the preceding path, and ensuring that no 
two $k'$s can  appear in the same row.     
\beq
\theta_{k}^{-1}:\ \
{\vcenter
{\offinterlineskip
\halign{&\mystrut\vrule#&\mybox{\hss$\scriptstyle#$\hss}\cr
\nr{8}&\hr{13}\cr 
   \omit&  &\omit&  &\omit& &\omit&  &&k'&&c&&  &\omit&\cdots&\omit&  &&s&  \cr
\hr{21}\cr 
   &\p&     &  &\omit&\cdots&\omit&  &&\a&&\b&&  &\omit&\cdots&\omit&  &&\q&  \cr
\hr{21}\cr
}}}
\quad
\begin{array}{c}
a>c\\
\ \ \longrightarrow\ \ \\
\end{array}
\quad
{\vcenter
{\offinterlineskip
\halign{&\mystrut\vrule#&\mybox{\hss$\scriptstyle#$\hss}\cr
\nr{8}&\hr{13}\cr 
   \omit&  &\omit&  &\omit& &\omit&  &&c&&k'&&  &\omit&\cdots&\omit&  &&s&  \cr
\hr{21}\cr
     &\p&&      &\omit&\cdots&\omit&  &&\a&&\b&&  &\omit&\cdots&\omit& &&\q&\cr
\hr{21}\cr 
}}}
\label{Eq-path-N}
\eeq
Proceeding in this way, the process terminates when each $k'$ has moved
as far east and south as the jeu de taquin allows. The fact that no
two identical $k'$s may lie in the same row is sufficient, given T1 and T2, to
show that the resulting primed shifted tableau satisfies PST1-4. Since we
had already noted that the diagonal entries are always unprimed, PST5
is also satisfied.

It follows that $\Theta^{-1}$ is well defined and maps the juxtaposition
of any pair of tableaux $PD\in{\cal PD}^\delta(n)$ and $T\in{\cal T}^\lambda(n)$ 
to a unique $PST\in{\cal PST}^\mu(n)$. Thus the original
map $\Theta$ from ${\cal PST}^\mu(n)$ to $\big({\cal PD}^\delta(n),{\cal T}^\lambda(n)\big)$
is indeed bijective. Since it is also weight preserving, as argued earlier, this completes
the proof of the $PST$ case in Theorem~\ref{The-PQ}. 

The only difference between the $PST$ and
$QST$ cases is the fact that in the latter case primed entries 
are allowed on the main diagonal. This is reflected in the
same distinction between $PD$ and $QD$ on the right of the 
above formulae. In fact it is not difficult to see that the map 
$\Theta$ preserves the entries on the main diagonal in both cases;
that is, just as the main diagonal of $PST$ coincides with that of $PD$,
where $\Theta:\ PST\,\mapsto\,(PD,T)$,
so the main diagonal of $QST$, complete with any primes, 
coincides with that of $QD$, where $\Theta:\ QST\,\mapsto\,(QD,T)$.
This observation is sufficient to complete the proof of Theorem~\ref{The-PQ}. 
 
\subsection{Example}
\label{subsec:gl-example}

This bijection is illustrated by the map from $PST$ of (\ref{Eq-PST})
to $PD$ of (\ref{Eq-PD}) and $T$ of (\ref{Eq-T}); that is,
\beq
PST=\ {\vcenter
{\offinterlineskip
\halign{&\mystrut\vrule#&\mybox{\hss$\scriptstyle#$\hss}\cr
\hr{19}\cr 
&1&&1&&1&&2'&&2&&2&&3&&3&&5&\cr
\hr{19}\cr 
\omit&
&&2&&2&&3'&&3&&4'&&5'&&5&&6'&\cr
\nr{2}&\hr{17}\cr 
\omit& &\omit&
&&3&&3&&4'&&4&&5'&&6&\cr 
\nr{4}&\hr{13}\cr 
\omit& &\omit& &\omit&
&&4&&5'&&5&&5&\cr
\nr{6}&\hr{9}\cr 
\omit& &\omit& &\omit& &\omit&
&&5&&6'&&6&\cr 
\nr{8}&\hr{7}\cr
\omit& &\omit& &\omit& &\omit& &\omit&
&&6&\cr 
\nr{10}&\hr{3}\cr
 }}}
\ \ \longleftrightarrow\ \ 
PD={\vcenter
{\offinterlineskip
\halign{&\mystrut\vrule#&\mybox{\hss$\scriptstyle#$\hss}\cr
\hr{13}\cr 
&1&&2'&&1&&4'&&5'&&6'&\cr 
\hr{13}\cr 
\omit& &&2&&3'&&2&&5'&&2&\cr 
\nr{2}&\hr{11}\cr 
\omit& &\omit&
&&3&&4'&&3&&3&\cr 
\nr{4}&\hr{9}\cr 
\omit& &\omit& &\omit&
&&4&&5'&&6'&\cr
\nr{6}&\hr{7}\cr 
\omit& &\omit& &\omit& &\omit&
&&5&&5&\cr 
\nr{8}&\hr{5}\cr
\omit& &\omit& &\omit& &\omit& &\omit&
&&6&\cr 
\nr{10}&\hr{3}\cr
 }}}
\ \ , \ \ T=\ \
{\vcenter
{\offinterlineskip
\halign{&\mystrut\vrule#&\mybox{\hss$\scriptstyle#$\hss}\cr
\hr{7}\cr 
&1&&2&&3&\cr 
\hr{7}\cr 
&3&&5&&5&\cr 
\hr{7}\cr 
&4&&6&\cr 
\hr{5}\cr 
&5&\cr
\hr{3}\cr 
&6&\cr 
\hr{3}\cr
\omit& &\omit\cr
\hr{0}\cr
 }}}
\label{Eq-PST-PD-T}
\eeq

The paths traced out by the primed entries $k'$ of $PST$
as they move
northwest as far as but no further than the $k$th column
are illustrated by means of boldface entries in the tableaux
shown below:

First moving the single $2'$ under the map $\theta_{2'}$ gives:
\beq
{\vcenter
{\offinterlineskip
\halign{&\mystrut\vrule#&\mybox{\hss$\scriptstyle#$\hss}\cr
\hr{19}\cr 
&1&&\bf1&&\bf1&&\bf{2'}&&2&&2&&3&&3&&5&\cr
\hr{19}\cr 
\omit&
&&2&&2&&3'&&3&&4'&&5'&&5&&6'&\cr
\nr{2}&\hr{17}\cr 
\omit& &\omit&
&&3&&3&&4'&&4&&5'&&6&\cr 
\nr{4}&\hr{13}\cr 
\omit& &\omit& &\omit&
&&4&&5'&&5&&5&\cr
\nr{6}&\hr{9}\cr 
\omit& &\omit& &\omit& &\omit&
&&5&&6'&&6&\cr 
\nr{8}&\hr{7}\cr
\omit& &\omit& &\omit& &\omit& &\omit&
&&6&\cr 
\nr{10}&\hr{3}\cr
 }}}
\ra
{\vcenter
{\offinterlineskip
\halign{&\mystrut\vrule#&\mybox{\hss$\scriptstyle#$\hss}\cr
\hr{19}\cr 
&1&&\bf{2'}&&\bf1&&\bf1&&2&&2&&3&&3&&5&\cr
\hr{19}\cr 
\omit&
&&2&&2&&3'&&3&&4'&&5'&&5&&6'&\cr
\nr{2}&\hr{17}\cr 
\omit& &\omit&
&&3&&3&&4'&&4&&5'&&6&\cr 
\nr{4}&\hr{13}\cr 
\omit& &\omit& &\omit&
&&4&&5'&&5&&5&\cr
\nr{6}&\hr{9}\cr 
\omit& &\omit& &\omit& &\omit&
&&5&&6'&&6&\cr 
\nr{8}&\hr{7}\cr
\omit& &\omit& &\omit& &\omit& &\omit&
&&6&\cr 
\nr{10}&\hr{3}\cr
 }}}
\eeq

Under $\theta_{3'}$ the only $3'$ moves just one step west where it has, as required,
reached the $3$rd column. It does not move north because the entry $1$
immediately above already lies in its own row:
\beq
{\vcenter
{\offinterlineskip
\halign{&\mystrut\vrule#&\mybox{\hss$\scriptstyle#$\hss}\cr
\hr{19}\cr 
&1&&2'&&1&&1&&2&&2&&3&&3&&5&\cr
\hr{19}\cr 
\omit&
&&2&&\bf2&&\bf{3'}&&3&&4'&&5'&&5&&6'&\cr
\nr{2}&\hr{17}\cr 
\omit& &\omit&
&&3&&3&&4'&&4&&5'&&6&\cr 
\nr{4}&\hr{13}\cr 
\omit& &\omit& &\omit&
&&4&&5'&&5&&5&\cr
\nr{6}&\hr{9}\cr 
\omit& &\omit& &\omit& &\omit&
&&5&&6'&&6&\cr 
\nr{8}&\hr{7}\cr
\omit& &\omit& &\omit& &\omit& &\omit&
&&6&\cr 
\nr{10}&\hr{3}\cr
 }}}
\ra
{\vcenter
{\offinterlineskip
\halign{&\mystrut\vrule#&\mybox{\hss$\scriptstyle#$\hss}\cr
\hr{19}\cr 
&1&&2'&&1&&1&&2&&2&&3&&3&&5&\cr
\hr{19}\cr 
\omit&
&&2&&\bf{3'}&&\bf2&&3&&4'&&5'&&5&&6'&\cr
\nr{2}&\hr{17}\cr 
\omit& &\omit&
&&3&&3&&4'&&4&&5'&&6&\cr 
\nr{4}&\hr{13}\cr 
\omit& &\omit& &\omit&
&&4&&5'&&5&&5&\cr
\nr{6}&\hr{9}\cr 
\omit& &\omit& &\omit& &\omit&
&&5&&6'&&6&\cr 
\nr{8}&\hr{7}\cr
\omit& &\omit& &\omit& &\omit& &\omit&
&&6&\cr 
\nr{10}&\hr{3}\cr
 }}}
\eeq

There are two $4'$s. Under $\theta_{4'}$ the upper one must be moved first
and then the lower one:
\beq
{\vcenter
{\offinterlineskip
\halign{&\mystrut\vrule#&\mybox{\hss$\scriptstyle#$\hss}\cr
\hr{19}\cr 
&1&&2'&&1&&\bf1&&\bf2&&2&&3&&3&&5&\cr
\hr{19}\cr 
\omit&
&&2&&3'&&2&&\bf3&&\bf{4'}&&5'&&5&&6'&\cr
\nr{2}&\hr{17}\cr 
\omit& &\omit&
&&3&&3&&4'&&4&&5'&&6&\cr 
\nr{4}&\hr{13}\cr 
\omit& &\omit& &\omit&
&&4&&5'&&5&&5&\cr
\nr{6}&\hr{9}\cr 
\omit& &\omit& &\omit& &\omit&
&&5&&6'&&6&\cr 
\nr{8}&\hr{7}\cr
\omit& &\omit& &\omit& &\omit& &\omit&
&&6&\cr 
\nr{10}&\hr{3}\cr
}}}
\ra
{\vcenter
{\offinterlineskip
\halign{&\mystrut\vrule#&\mybox{\hss$\scriptstyle#$\hss}\cr
\hr{19}\cr 
&1&&2'&&1&&\bf{4'}&&\bf1&&2&&3&&3&&5&\cr
\hr{19}\cr 
\omit&
&&2&&3'&&2&&\bf2&&\bf3&&5'&&5&&6'&\cr
\nr{2}&\hr{17}\cr 
\omit& &\omit&
&&3&&3&&4'&&4&&5'&&6&\cr 
\nr{4}&\hr{13}\cr 
\omit& &\omit& &\omit&
&&4&&5'&&5&&5&\cr
\nr{6}&\hr{9}\cr 
\omit& &\omit& &\omit& &\omit&
&&5&&6'&&6&\cr 
\nr{8}&\hr{7}\cr
\omit& &\omit& &\omit& &\omit& &\omit&
&&6&\cr 
\nr{10}&\hr{3}\cr
 }}}
\ \ = \  \
{\vcenter
{\offinterlineskip
\halign{&\mystrut\vrule#&\mybox{\hss$\scriptstyle#$\hss}\cr
\hr{19}\cr 
&1&&2'&&1&&4'&&1&&2&&3&&3&&5&\cr
\hr{19}\cr 
\omit&
&&2&&3'&&2&&2&&3&&5'&&5&&6'&\cr
\nr{2}&\hr{17}\cr 
\omit& &\omit&
&&3&&\bf3&&\bf{4'}&&4&&5'&&6&\cr 
\nr{4}&\hr{13}\cr 
\omit& &\omit& &\omit&
&&4&&5'&&5&&5&\cr
\nr{6}&\hr{9}\cr 
\omit& &\omit& &\omit& &\omit&
&&5&&6'&&6&\cr 
\nr{8}&\hr{7}\cr
\omit& &\omit& &\omit& &\omit& &\omit&
&&6&\cr 
\nr{10}&\hr{3}\cr
 }}}
\ra
{\vcenter
{\offinterlineskip
\halign{&\mystrut\vrule#&\mybox{\hss$\scriptstyle#$\hss}\cr
\hr{19}\cr 
&1&&2'&&1&&4'&&1&&2&&3&&3&&5&\cr
\hr{19}\cr 
\omit&
&&2&&3'&&2&&2&&3&&5'&&5&&6'&\cr
\nr{2}&\hr{17}\cr 
\omit& &\omit&
&&3&&\bf{4'}&&\bf3&&4&&5'&&6&\cr 
\nr{4}&\hr{13}\cr 
\omit& &\omit& &\omit&
&&4&&5'&&5&&5&\cr
\nr{6}&\hr{9}\cr 
\omit& &\omit& &\omit& &\omit&
&&5&&6'&&6&\cr 
\nr{8}&\hr{7}\cr
\omit& &\omit& &\omit& &\omit& &\omit&
&&6&\cr 
\nr{10}&\hr{3}\cr
 }}}
\eeq

There are three $5'$s to deal with in turn from top to bottom using 
$\theta_{5'}$, but the last of these is already in the $3$rd column 
and directly below a $3$ in the $3$rd row, and so does not move:
\beq
{\vcenter
{\offinterlineskip
\halign{&\mystrut\vrule#&\mybox{\hss$\scriptstyle#$\hss}\cr
\hr{19}\cr 
&1&&2'&&1&&4'&&\bf1&&\bf2&&\bf3&&3&&5&\cr
\hr{19}\cr 
\omit&
&&2&&3'&&2&&2&&3&&\bf{5'}&&5&&6'&\cr
\nr{2}&\hr{17}\cr 
\omit& &\omit&
&&3&&4'&&3&&4&&5'&&6&\cr 
\nr{4}&\hr{13}\cr 
\omit& &\omit& &\omit&
&&4&&5'&&5&&5&\cr
\nr{6}&\hr{9}\cr 
\omit& &\omit& &\omit& &\omit&
&&5&&6'&&6&\cr 
\nr{8}&\hr{7}\cr
\omit& &\omit& &\omit& &\omit& &\omit&
&&6&\cr 
\nr{10}&\hr{3}\cr
 }}}
\ \ra\ 
{\vcenter
{\offinterlineskip
\halign{&\mystrut\vrule#&\mybox{\hss$\scriptstyle#$\hss}\cr
\hr{19}\cr 
&1&&2'&&1&&4'&&\bf{5'}&&\bf1&&\bf2&&3&&5&\cr
\hr{19}\cr 
\omit&
&&2&&3'&&2&&2&&3&&\bf3&&5&&6'&\cr
\nr{2}&\hr{17}\cr 
\omit& &\omit&
&&3&&4'&&3&&4&&5'&&6&\cr 
\nr{4}&\hr{13}\cr 
\omit& &\omit& &\omit&
&&4&&5'&&5&&5&\cr
\nr{6}&\hr{9}\cr 
\omit& &\omit& &\omit& &\omit&
&&5&&6'&&6&\cr 
\nr{8}&\hr{7}\cr
\omit& &\omit& &\omit& &\omit& &\omit&
&&6&\cr 
\nr{10}&\hr{3}\cr
 }}}
\ =\  
{\vcenter
{\offinterlineskip
\halign{&\mystrut\vrule#&\mybox{\hss$\scriptstyle#$\hss}\cr
\hr{19}\cr 
&1&&2'&&1&&4'&&5'&&1&&2&&3&&5&\cr
\hr{19}\cr 
\omit&
&&2&&3'&&2&&\bf2&&\bf3&&3&&5&&6'&\cr
\nr{2}&\hr{17}\cr 
\omit& &\omit&
&&3&&4'&&3&&\bf4&&\bf{5'}&&6&\cr 
\nr{4}&\hr{13}\cr 
\omit& &\omit& &\omit&
&&4&&5'&&5&&5&\cr
\nr{6}&\hr{9}\cr 
\omit& &\omit& &\omit& &\omit&
&&5&&6'&&6&\cr 
\nr{8}&\hr{7}\cr
\omit& &\omit& &\omit& &\omit& &\omit&
&&6&\cr 
\nr{10}&\hr{3}\cr
 }}}
\ \ra\
{\vcenter
{\offinterlineskip
\halign{&\mystrut\vrule#&\mybox{\hss$\scriptstyle#$\hss}\cr
\hr{19}\cr 
&1&&2'&&1&&4'&&5'&&1&&2&&3&&5&\cr
\hr{19}\cr 
\omit&
&&2&&3'&&2&&\bf{5'}&&\bf2&&3&&5&&6'&\cr
\nr{2}&\hr{17}\cr 
\omit& &\omit&
&&3&&4'&&3&&\bf3&&\bf4&&6&\cr 
\nr{4}&\hr{13}\cr 
\omit& &\omit& &\omit&
&&4&&5'&&5&&5&\cr
\nr{6}&\hr{9}\cr 
\omit& &\omit& &\omit& &\omit&
&&5&&6'&&6&\cr 
\nr{8}&\hr{7}\cr
\omit& &\omit& &\omit& &\omit& &\omit&
&&6&\cr 
\nr{10}&\hr{3}\cr
 }}} 
\eeq

Then we deal with the two $6'$s by applying $\theta_{6'}$ to give
\beq
{\vcenter
{\offinterlineskip
\halign{&\mystrut\vrule#&\mybox{\hss$\scriptstyle#$\hss}\cr
\hr{19}\cr 
&1&&2'&&1&&4'&&5'&&\bf1&&\bf2&&\bf3&&\bf5&\cr
\hr{19}\cr 
\omit&
&&2&&3'&&2&&5'&&2&&3&&5&&\bf{6'}&\cr
\nr{2}&\hr{17}\cr 
\omit& &\omit&
&&3&&4'&&3&&3&&4&&6&\cr 
\nr{4}&\hr{13}\cr 
\omit& &\omit& &\omit&
&&4&&5'&&5&&5&\cr
\nr{6}&\hr{9}\cr 
\omit& &\omit& &\omit& &\omit&
&&5&&6'&&6&\cr 
\nr{8}&\hr{7}\cr
\omit& &\omit& &\omit& &\omit& &\omit&
&&6&\cr 
\nr{10}&\hr{3}\cr
 }}}
\ \ra\
{\vcenter
{\offinterlineskip
\halign{&\mystrut\vrule#&\mybox{\hss$\scriptstyle#$\hss}\cr
\hr{19}\cr 
&1&&2'&&1&&4'&&5'&&\bf{6'}&&\bf1&&\bf2&&\bf3&\cr
\hr{19}\cr 
\omit&
&&2&&3'&&2&&5'&&2&&3&&5&&\bf5&\cr
\nr{2}&\hr{17}\cr 
\omit& &\omit&
&&3&&4'&&3&&3&&4&&6&\cr 
\nr{4}&\hr{13}\cr 
\omit& &\omit& &\omit&
&&4&&5'&&5&&5&\cr
\nr{6}&\hr{9}\cr 
\omit& &\omit& &\omit& &\omit&
&&5&&6'&&6&\cr 
\nr{8}&\hr{7}\cr
\omit& &\omit& &\omit& &\omit& &\omit&
&&6&\cr 
\nr{10}&\hr{3}\cr
 }}}
\ \ = \ \
{\vcenter
{\offinterlineskip
\halign{&\mystrut\vrule#&\mybox{\hss$\scriptstyle#$\hss}\cr
\hr{19}\cr 
&1&&2'&&1&&4'&&5'&&6'&&1&&2&&3&\cr
\hr{19}\cr 
\omit&
&&2&&3'&&2&&5'&&2&&3&&5&&5&\cr
\nr{2}&\hr{17}\cr 
\omit& &\omit&
&&3&&4'&&3&&3&&4&&6&\cr 
\nr{4}&\hr{13}\cr 
\omit& &\omit& &\omit&
&&4&&5'&&\bf5&&5&\cr
\nr{6}&\hr{9}\cr 
\omit& &\omit& &\omit& &\omit&
&&5&&\bf{6'}&&6&\cr 
\nr{8}&\hr{7}\cr
\omit& &\omit& &\omit& &\omit& &\omit&
&&6&\cr 
\nr{10}&\hr{3}\cr
 }}}
\ra
{\vcenter
{\offinterlineskip
\halign{&\mystrut\vrule#&\mybox{\hss$\scriptstyle#$\hss}\cr
\hr{19}\cr 
&1&&2'&&1&&4'&&5'&&6'&&1&&2&&3&\cr
\hr{19}\cr 
\omit&
&&2&&3'&&2&&5'&&2&&3&&5&&5&\cr
\nr{2}&\hr{17}\cr 
\omit& &\omit&
&&3&&4'&&3&&3&&4&&6&\cr 
\nr{4}&\hr{13}\cr 
\omit& &\omit& &\omit&
&&4&&5'&&\bf{6'}&&5&\cr
\nr{6}&\hr{9}\cr 
\omit& &\omit& &\omit& &\omit&
&&5&&\bf5&&6&\cr 
\nr{8}&\hr{7}\cr
\omit& &\omit& &\omit& &\omit& &\omit&
&&6&\cr 
\nr{10}&\hr{3}\cr
 }}}
\eeq

This results in the juxtaposition
of $PD$ from (\ref{Eq-PD}) and $T$ from (\ref{Eq-T}) as claimed:
\beq
{\vcenter
{\offinterlineskip
\halign{&\mystrut\vrule#&\mybox{\hss$\scriptstyle#$\hss}\cr
\hr{19}\cr 
&1&&2'&&1&&4'&&5'&&6'&&1&&2&&3&\cr
\hr{19}\cr 
\omit&
&&2&&3'&&2&&5'&&2&&3&&5&&5&\cr
\nr{2}&\hr{17}\cr 
\omit& &\omit&
&&3&&4'&&3&&3&&4&&6&\cr 
\nr{4}&\hr{13}\cr 
\omit& &\omit& &\omit&
&&4&&5'&&6'&&5&\cr
\nr{6}&\hr{9}\cr 
\omit& &\omit& &\omit& &\omit&
&&5&&5&&6&\cr 
\nr{8}&\hr{7}\cr
\omit& &\omit& &\omit& &\omit& &\omit&
&&6&\cr 
\nr{10}&\hr{3}\cr
 }}}
\ \ \equiv \ \
{\vcenter
{\offinterlineskip
\halign{&\mystrut\vrule#&\mybox{\hss$\scriptstyle#$\hss}\cr
\hr{13}\cr 
&1&&2'&&1&&4'&&5'&&6'&\cr 
\hr{13}\cr 
\omit& &&2&&3'&&2&&5'&&2&\cr 
\nr{2}&\hr{11}\cr 
\omit& &\omit&
&&3&&4'&&3&&3&\cr 
\nr{4}&\hr{9}\cr 
\omit& &\omit& &\omit&
&&4&&5'&&6'&\cr
\nr{6}&\hr{7}\cr 
\omit& &\omit& &\omit& &\omit&
&&5&&5&\cr 
\nr{8}&\hr{5}\cr
\omit& &\omit& &\omit& &\omit& &\omit&
&&6&\cr 
\nr{10}&\hr{3}\cr
 }}}
\quad\cdot\quad
{\vcenter
{\offinterlineskip
\halign{&\mystrut\vrule#&\mybox{\hss$\scriptstyle#$\hss}\cr
\hr{7}\cr 
&1&&2&&3&\cr 
\hr{7}\cr 
&3&&5&&5&\cr 
\hr{7}\cr 
&4&&6&\cr 
\hr{5}\cr 
&5&\cr
\hr{3}\cr 
&6&\cr 
\hr{3}\cr
\omit& &\omit\cr
\hr{0}\cr
 }}}
\eeq

\subsection{Corollaries}
\label{subsec:gl-corollaries}

By associating $x_k$ and $y_k$ to each entry $k$ and
$k'$, respectively, in the various tableaux $PST$, $QST$,
$PD$, $QD$ and $T$ appearing in Theorem~\ref{The-PQ}
we immediately have the following corollary.

\begin{Corollary}
Let $\mu=\lambda+\delta$ be a strict partition of length $\ell(\mu)=n$,
with $\lambda$ a partition of length $\ell(\lambda)\leq n$ and 
$\delta=(n,n-1,\ldots,1)$. 
\beq
\begin{array}{rcl}
 \sum_{PST\in{\cal PST}^\mu(n)}\ (\x/\y)^{\wgt(PST)}
 &=& \sum_{PD\in{\cal PD}^\delta(n)}\ (\x/\y)^{\wgt(PD)} 
    \ \sum_{T\in{\cal T}^\lambda(n)}\ \x^{\wgt(T)}\,; \\ \\
  \sum_{QST\in{\cal QST}^\mu(n)}\ (\x/\y)^{\wgt(QST)}
 &=& \sum_{QD\in{\cal QD}^\delta(n)}\ (\x/\y)^{\wgt(QD)} 
    \ \sum_{T\in{\cal T}^\lambda(n)}\ \x^{\wgt(T)}\,. \\ \\
\end{array}
\label{Eq-PST-QST-PD-T}
\eeq
\label{Cor-PQ}
\end{Corollary}

Thanks to the definitions of $P(\x/\y)$ and $Q(\x/\y)$
given in (\ref{Eq-PQxy}), the identities (\ref{Eq-PDprod})
and (\ref{Eq-QDprod}) and the combinatorial definition
of $s_\lambda(\x)$ given in (\ref{Eq-Sf}), the above 
result is nothing other than our first main result,
Proposition~\ref{Prop-PQ}.

Other corollaries follow as special cases of these results.
Setting $\lambda=0$ we obtain 
\begin{equation}
P_\delta(\x/\y)=s_{1^{n}}(\x) \prod_{1\leq i < j \leq n} (x_i + y_j)
\quad\hbox{and}\quad
Q_{\delta}(\x/\y)=\prod_{1\leq i\leq j \leq n} (x_i + y_j).
\end{equation}
Further specialisation to the case $\y=\x$ leads to a result 
given by Macdonald~\cite{M95}(Sec.III.8,Ex.3 p259):
\begin{equation}
 P_\delta(\x)=s_\delta(\x) \quad\hbox{and}\quad 
 Q_\delta(\x)=2^n\ s_\delta(\x),
\end{equation} 
where use has been made of the fact that 
\begin{equation}
\prod_{i=1}^n x_i\ \prod_{1\leq i<j\leq n}(x_i+x_j)=
s_{1^n}(\x)\  \prod_{1\leq i<j\leq
  n}\frac{(x_i^2-x_j^2)}{(x_i-x_j)}
=s_{1^n}(\x)\ s_{\delta/1^n}(\x)=s_\delta(\x),
\label{Eq-s1ndelta}
\end{equation}
where the last step is true when,
as here, $\x$ has $n$ components $x_1,x_2,\ldots,x_n$.

More generally, if $\mu=\lambda+\delta$ for any partition
$\lambda$ of length $\ell(\lambda)\leq n$, 
but $\y=\x$ we have
another result due to Macdonald~\cite{M95}(Sec.III.8,Ex.2 p259):
\begin{equation}
 P_{\lambda+\delta}(\x)= s_\delta(\x)\ s_\lambda(\x).
\end{equation}
where (\ref{Eq-s1ndelta}) and (\ref{Eq-QxPx}) have been applied
 directly to the $\y=\x$ case of (\ref{Eq-PQ}).

On the other hand the case
$\y=t\x=(tx_1,tx_2,\ldots,tx_n)$ of (\ref{Eq-PQ}) is equivalent to
(\ref{Eq-Tokuyama}), 
the $t$-deformation of Weyl's denominator formula for the
Lie algebra $gl(n)$ due to Tokuyama \cite{T88}:

\begin{Corollary}
\begin{equation}
\prod_{i=1}^n x_i\ \prod_{1\leq i < j \leq n} (x_i+tx_j)\ s_{\lambda}(\x) =
\sum_{ST\in {\cal ST}^{\mu}(n)} t^{\hgt(ST)} (1+t)^{\str(ST)-n}\ \x^{\wgt(ST)},
\end{equation}
\end{Corollary}

\noindent{\bf Proof}\ \
While there is a combinatorial proof of this result due to
Okada~\cite{O90}, it follows immediately from Theorem~\ref{The-PQ}
by setting $y_k=tx_k$ for all $k=1,2,\ldots,n$, noting that deleting 
primes from the entries $k'$ in each $PST\in{\cal PST}^\mu(n)$
gives a shifted tableaux $ST\in{\cal ST}^\mu(n)$ with a
factor of $t$ arising from each primed entry of $PST$, and observing 
that these must occur in precisely those boxes contributing to
$\hgt(ST)$ and are optional, thereby giving rise to a factor
of $(1+t)$ in those $\str(ST)-n$ boxes at the lower left hand end
of all continuous strips of identical entries other than those
starting on the main diagonal.

The remaining corollaries mentioned in the Introduction
are the formulae (\ref{Eq-RR}) of Robbins and Rumsey~\cite{RR86}
and (\ref{Eq-Chapman}) of Chapman~\cite{C01}. These require for their 
elucidation a link with alternating sign matrices. This is
provided in Section \ref{sec:asm}.

\section{The $sp(2n)$ bijection}
\label{sec:sp-bijection}

\subsection{Main Result}
\label{Main Result}

The analogue in the symplectic case of Theorem~\ref{The-PQ} is the
following:

\begin{Theorem}
Let $\mu=\lambda+\delta$ be a strict partition of length $\ell(\mu)=n$,
with $\lambda$ a partition of length $\ell(\lambda)\leq n$ and 
$\delta=(n,n-1,\ldots,1)$. There exists a weight and barred weight
preserving, bijective map 
$\Phi$ from ${\cal QST}^\mu(n,{\ov n})$ to 
$\big({\cal QD}^\delta(n,{\ov n}),{\cal T}^\lambda(n,{\ov n})\big)$
such that for all $QST\in{\cal QST}^\mu(n,{\ov n})$ 
\beq
   \Phi:\ QST\ \mapsto\ (QD,T)  
\quad\hbox{with}\quad\left\{
\begin{array}{l}
    \wgt(QST)=\wgt(QD)+\wgt(T) \\ \\
    \bar(QST)=\bar(QD)+\bar(T) \\
\end{array}
\right.
\label{Eq-Phi}
\eeq
and $QD\in{\cal QD}^\delta(n,{\ov n})$ and
$T\in{\cal T}^\lambda(n,{\ov n})$. 
\label{The-Qbar}
\end{Theorem}

\noindent{\bf Proof}\ \

The Theorem is proved by the identification of a suitable map $\Phi$  
that it is both weight preserving and bijective. The underlying procedure is the same as 
before, in that the jeu de taquin is applied successively to all primed entries 
of $QST$ dealing in sequence with all entries ${\ov k}'$ and then $k'$
for $k=1,2\ldots,n$. In the case of the ${\ov k}'$s there is no impediment to
moving all these entries to the $k$th column by means of the jeu de taquin, but 
something slightly more subtle is required in the case of the $k'$s. As well, 
the transformations applied to entries in the $k$th column now include not
only permutations but an additional weight preserving transformation.

The structure of $\Phi$ is such that 
$\Phi=\phi_{n'}\circ\phi_{\ov{n}'}\circ\cdots
\circ\phi_{2'}\circ\phi_{\ov2'}\circ\phi_{1'}\circ\phi_{\ov1'}$.
Here $\phi_{\ov{k}'}$ differs from $\theta_{k'}$ only in that the 
jeu de taquin is played with ${\ov k}'$s rather than the $k'$s,
while $\phi_{k'}=\chi_{k'}\circ\psi_{k'}$ 
where $\psi_{k'}$ differs from $\theta_{k'}$ only if
 the final step of the path of $k'$ into the $i$th row of the 
$k$th column is blocked by an entry $\ov{k}'$. 
In such a situation the horizontal pair of 
entries $\ov{k}'k'$ in the $i$th row is replaced by the 
horizontal pair $i\,\ov{i}$. 
Having moved all the $k'$s into the $k$th column
or annihilated them as above, there may remain in the $k$th column
vertical pairs $\ov{i}\,i$. It is then necessary to
invoke $\chi_{k'}$. This replaces the lowest such
pair, for which $i$ is necessarily in the $i$th row, 
by a vertical pair $k'\ov{k}'$, moves the resulting $k'$ 
and $\ov{k}'$ north as far as possible whilst still
satisfyng T$\ov3$, and then acts in the same way
on the next lowest vertical pair $\ov{j}\,j$, replacing 
them by another vertical pair $k'\ov{k}'$, and so on. 
Having removed all unprimed vertical pairs in this way
any remaining unbarred entries $i$ or $\ov{i}$, but not both, lie 
in their own $i$th row, as required for consistency with QD$\ov1$.

To demonstrate that these various maps are well defined, we exhibit 
the relevant individual steps as below, first for the $\ov{k}'$ case.
The starting point is the 
tableau $T_0$ in which the entries $a,b,\ldots,h$ satisfy the semistandardness
conditions T1 and T2, and the symplectic constraint T$\ov3$. 
This tableau $T_0$ is to be thought of as the
subtableau surrounding a particular $\ov{k}'$ after the jeu de taquin has
been applied to all $\ov{k}'$s appearing initially above the $\ov{k}'$ in
question, moving them into the $k$th column. Some further
steps of the jeu de taquin may have already been applied to the central $\ov{k}'$
and the diagram is intended to indicate under what conditions its next
move is north or west, given that it currently lies in the $i$th row. 

If $d\leq b$ then  $\ov{k}'$ is to be interchanged
with $b$ and if $d>b$ then $\ov{k}'$ is to be interchanged with
$d$. In the first case $\ov{k}'$ moves north and the resulting tableau $T_N$
satisfies  both  $T1$  and $T2$ since $d\leq b\leq c<e$, while in the
second $\ov{k}'$ moves west and the resulting tableau $T_W$
satisfies  both  $T1$  and $T2$ since $b<d<f\leq g$. \
What is particularly important to note here is that by virtue of T$\ov3$
applied to $T_0$ it is known that $d\geq\ov{i}$. It follows that
if $d\leq b$ and the map is from $T_0$ to $T_N$ then $b\geq\ov{i}$
so that $T_N$ also satisfies T$\ov3$. Conversely if $b<\ov{i}$ then
we must have $b<d$ and the map must be from $T_0$ to $T_W$ in which case
T$\ov3$ is still satisfied.
\beq
\phi_{\ov{k}'}:\ \
\begin{array}{c}
T_0=\ {\vcenter
{\offinterlineskip
\halign{&\mystrut\vrule#&\mybox{\hss$\scriptstyle#$\hss}\cr
\hr{7}\cr 
&a&&b&&c&\cr
\hr{7}\cr 
&d&&\ov{k}'&&e&\cr
\hr{7}\cr
&f&&g&&h&\cr 
\hr{7}\cr
 }}}\\
\hbox{$\ov{k}'$ in $i$th row}\\
\ov{i}\leq d\\
\end{array}
\ \ \longrightarrow\ \ 
\left\{
\begin{array}{ll} 
T_N=\ {\vcenter
{\offinterlineskip
\halign{&\mystrut\vrule#&\mybox{\hss$\scriptstyle#$\hss}\cr
\hr{7}\cr 
&a&&\ov{k}'&&c&\cr
\hr{7}\cr 
&d&&b&&e&\cr
\hr{7}\cr
&f&&g&&h&\cr 
\hr{7}\cr
 }}}
&\quad\hbox{if}\quad \ov{i}\leq d\leq b
\\  \\
T_W=\ {\vcenter
{\offinterlineskip
\halign{&\mystrut\vrule#&\mybox{\hss$\scriptstyle#$\hss}\cr
\hr{7}\cr 
&a&&b&&c&\cr
\hr{7}\cr 
&\ov{k}'&&d&&e&\cr
\hr{7}\cr
&f&&g&&h&\cr 
\hr{7}\cr
 }}}
&\quad\hbox{if}\quad b<\ov{i}\leq d
\\
\end{array}
\right.
\label{Eq-phi-NW}
\eeq
If $\ov{k}'$ is already in the topmost row or leftmost column the
following diagrams illustrate the allowed 
moves of $\ov{k}'$ 
west and north respectively:
\beq
\phi_{\ov{k}'}:\ \
T_0=\ {\vcenter
{\offinterlineskip
\halign{&\mystrut\vrule#&\mybox{\hss$\scriptstyle#$\hss}\cr
\hr{7}\cr 
&d&&\ov{k}'&&e&\cr
\hr{7}\cr
&f&&g&&h&\cr 
\hr{7}\cr
 }}}
\quad
\ \ \longrightarrow\ \ 
\quad
T_W=\ {\vcenter
{\offinterlineskip
\halign{&\mystrut\vrule#&\mybox{\hss$\scriptstyle#$\hss}\cr
\hr{7}\cr 
&\ov{k}'&&d&&e&\cr
\hr{7}\cr
&f&&g&&h&\cr 
\hr{7}\cr
 }}}
\label{Eq-phi-W}
\eeq
and
\beq
\phi_{\ov{k}'}:\ \
\begin{array}{c}
T_0=\ {\vcenter
{\offinterlineskip
\halign{&\mystrut\vrule#&\mybox{\hss$\scriptstyle#$\hss}\cr
\hr{5}\cr 
&b&&c&\cr
\hr{5}\cr 
&\ov{k}'&&e&\cr
\hr{5}\cr
&g&&h&\cr 
\hr{5}\cr
 }}}\\
\hbox{$\ov{k}'$ in $i$th row}\\
\end{array}
\quad
\ \ \longrightarrow\ \
\quad
T_N=\ {\vcenter
{\offinterlineskip
\halign{&\mystrut\vrule#&\mybox{\hss$\scriptstyle#$\hss}\cr
\hr{5}\cr 
&\ov{k}'&&c&\cr
\hr{5}\cr 
&b&&e&\cr
\hr{5}\cr
&g&&h&\cr 
\hr{5}\cr
 }}}
\quad\hbox{if}\quad \ov{i}\leq b.
\label{Eq-phi-N}
\eeq
Notice that in order to satisfy T$\ov3$ the map from $T_0$ to $T_N$ 
only takes place if $b\geq\ov{i}$, otherwise $T_0$ is unaltered
and $\ov{k}'$ occupies the site in the $i$th row of the $k$th column,
until disturbed by any incoming $k'$s, as we shall see later. In the 
meantime it is to be noted that $T_W$ satisfies T$\ov3$, since
by hypothesis $T_0$ does so. By the same token, 
since the unprimed entries of $T_0$ satisfy both T1 and T2,
these conditions are also satisfied by both $T_W$ and $T_N$ 
by virtue of the fact that we still have $d<f\leq g$ and $b\leq c<e$.

It still remains to be shown that by means of the above moves of 
$\ov{k}'$ in a north-westerly direction it actually reaches the $k$th
column. This time the only possible impediments to this is the 
existence of another $\ov{k}'$ to its immediate left. However, by the
same argument as that used in the $\theta_{k'}$ case, this 
cannot occur because the path followed by $\ov{k}'$ always remains 
strictly south of the path followed by all preceding $\ov{k}'$s.

Having completed the jeu de taquin moves for all $\ov{k}'$s and moved them
as far north as possible in the $k$th column, it remains to deal with 
any $k'$s in $QST$ using $\phi_{k'}=\chi_{k'}\circ\psi_{k'}$. 
The action of $\psi_{k'}$ is carried out in the same way as before
with the diagrams of (\ref{Eq-phi-NW}), (\ref{Eq-phi-W}) and (\ref{Eq-phi-N})
just altered by changing $\ov{k}'$ to $k'$, provided that $d$ is unprimed.
If $d$ is primed then $d\neq k'$ since $d=k'$ would give a horizontal pair $k'k'$,
and this cannot occur since the path of the second $k'$ must stay
strictly south of the path of the first $k'$. To deal with the case 
$d=\ov{k}'$ the following map is required if $a$ is unprimed
\beq
\phi_{{k}'}:\ \
\begin{array}{c}
T_0=\ {\vcenter
{\offinterlineskip
\halign{&\mystrut\vrule#&\mybox{\hss$\scriptstyle#$\hss}\cr
\hr{7}\cr 
&a&&b&&c&\cr
\hr{7}\cr 
&\ov{k}'&&k'&&e&\cr
\hr{7}\cr
&f&&g&&h&\cr 
\hr{7}\cr
 }}}\\ 
\hbox{$k'$ in $i$th row}\\
a<\ov{i}\\
\end{array}
\ \ \longrightarrow\ \ 
\left\{
\begin{array}{ll} 
T_N=\ {\vcenter
{\offinterlineskip
\halign{&\mystrut\vrule#&\mybox{\hss$\scriptstyle#$\hss}\cr
\hr{7}\cr 
&a&&k'&&c&\cr
\hr{7}\cr 
&\ov{k}'&&b&&e&\cr
\hr{7}\cr
&f&&g&&h&\cr 
\hr{7}\cr
 }}}
&\quad\hbox{if}\quad a<\ov{i}\leq b
\\  \\
T_W=\ {\vcenter
{\offinterlineskip
\halign{&\mystrut\vrule#&\mybox{\hss$\scriptstyle#$\hss}\cr
\hr{7}\cr 
&a&&b&&c&\cr
\hr{7}\cr 
&i&&\ov{i}&&e&\cr
\hr{7}\cr
&f&&g&&h&\cr 
\hr{7}\cr
 }}}
&\quad\hbox{if}\quad a\leq b<\ov{i}
\\
\end{array}
\right.
\label{Eq-phi-NWkk}
\eeq
where advantage has been taken of the fact that $a<\ov{i}<i$, since
if this were not the case the action of $\phi_{\ov{k}'}$ would
have required that $a$ and $\ov{k}'$ be interchanged in $T_0$. 
It follows that $T_W$ satisfies T$\ov3$.
The fact that the pair $i\,\ov{i}$ violates $T1$ is not a problem
because the $i$ lies in the $k$th column of what will become $QD$
and the $\ov{i}$ lies either in the $(k+1)$th column of $QD$
or the first column of $T$. In neither case does the condition 
T$1$ apply to both $i$ and $\ov{i}$. If they are both in $QD$
they automatically satisfy the condition QD1=PD1, and if one is 
in $QD$ and the other in $T$ then $i$ and $\ov{i}$ automatically
satisfy QST$\ov5$ and the combination of ST1 and ST$\ov4$, 
respectively. In the case that $\phi_{k'}$ maps $T_0$ to
$T_N$, then one reverts to the use of $\psi_{k'}=\theta_{k'}$ for the
next step (or equivalently $\phi_{\ov{k}'}$ with $\ov{k}'$ replaced
by $k'$ in (\ref{Eq-phi-NW})-(\ref{Eq-phi-N})).

Finally, if $d=\ov{k}'$ and $a$ is primed then for $a=\ov{k}'$ the action
of $\psi_{k'}$ gives
\beq
\psi_{k'}:\ \
\begin{array}{c}
T_0=\ {\vcenter
{\offinterlineskip
\halign{&\mystrut\vrule#&\mybox{\hss$\scriptstyle#$\hss}\cr
\hr{7}\cr 
&\ov{k}'&&b&&c&\cr
\hr{7}\cr 
&\ov{k}'&&k'&&e&\cr
\hr{7}\cr
&f&&g&&h&\cr 
\hr{7}\cr
 }}}\\ 
\hbox{$k'$ in $i$th row}\\
\end{array}
\ \ \longrightarrow\ \ 
\left\{
\begin{array}{ll} 
T_N=\ {\vcenter
{\offinterlineskip
\halign{&\mystrut\vrule#&\mybox{\hss$\scriptstyle#$\hss}\cr
\hr{7}\cr 
&\ov{k}'&&k'&&c&\cr
\hr{7}\cr 
&\ov{k}'&&b&&e&\cr
\hr{7}\cr
&f&&g&&h&\cr 
\hr{7}\cr
 }}}
&\quad\hbox{if}\quad \ov{i}\leq b
\\  \\
T_W=\ {\vcenter
{\offinterlineskip
\halign{&\mystrut\vrule#&\mybox{\hss$\scriptstyle#$\hss}\cr
\hr{7}\cr 
&\ov{k}'&&b&&c&\cr
\hr{7}\cr 
&i&&\ov{i}&&e&\cr
\hr{7}\cr
&f&&g&&h&\cr 
\hr{7}\cr
 }}}
&\quad\hbox{if}\quad b<\ov{i}
\\
\end{array}
\right.
\label{Eq-psi-NWkkk}
\eeq
while for $ a=k'$ the action is
\beq
\psi_{k'}:\ \
\begin{array}{c}
T_0=\ {\vcenter
{\offinterlineskip
\halign{&\mystrut\vrule#&\mybox{\hss$\scriptstyle#$\hss}\cr
\hr{7}\cr 
&k'&&b&&c&\cr
\hr{7}\cr 
&\ov{k}'&&k'&&e&\cr
\hr{7}\cr
&f&&g&&h&\cr 
\hr{7}\cr
 }}}\\ 
\hbox{$\ov{k}'$ in $i$th row}\\
\end{array}
\quad
\begin{array}{c}
\ b<\ov{i}\\ 
\ \ \longrightarrow\ \ \\
\end{array}
\quad
T_W=\ {\vcenter
{\offinterlineskip
\halign{&\mystrut\vrule#&\mybox{\hss$\scriptstyle#$\hss}\cr
\hr{7}\cr 
&k'&&b&&c&\cr
\hr{7}\cr 
&i&&\ov{i}&&e&\cr
\hr{7}\cr
&f&&g&&h&\cr 
\hr{7}\cr
 }}}
\label{Eq-psi-Wkkk}
\eeq
since the uppermost $k'$ must have arrived at its present position by
displacing $b$, and $b$ could only have been in that position
above $\ov{k}'$ in the $i$th row if it satisfied $b<\ov{i}$. This
prevents the lower $k'$ from moving north by interchanging with $b$,
since the latter move would violate T$\ov3$.

If $k'$ is already in the topmost row or leftmost column the
following diagrams (\ref{Eq-psi-W})and (\ref{Eq-psi-N}) illustrate the allowed moves of $k'$ 
west and north respectively:
\beq
\psi_{k'}:\ \
T_0=\ {\vcenter
{\offinterlineskip
\halign{&\mystrut\vrule#&\mybox{\hss$\scriptstyle#$\hss}\cr
\hr{7}\cr 
&d&&k'&&e&\cr
\hr{7}\cr
&f&&g&&h&\cr 
\hr{7}\cr
 }}}
\qquad
\ \ \longrightarrow\ \
T_W=\ {\vcenter
{\offinterlineskip
\halign{&\mystrut\vrule#&\mybox{\hss$\scriptstyle#$\hss}\cr
\hr{7}\cr 
&k'&&d&&e&\cr
\hr{7}\cr
&f&&g&&h&\cr 
\hr{7}\cr
 }}}
\label{Eq-psi-W}
\eeq
and
\beq
\psi_{k'}:\ \
\begin{array}{c}
T_0=\ {\vcenter
{\offinterlineskip
\halign{&\mystrut\vrule#&\mybox{\hss$\scriptstyle#$\hss}\cr
\hr{5}\cr 
&b&&c&\cr
\hr{5}\cr 
&\ov{k}'&&e&\cr
\hr{5}\cr
&g&&h&\cr 
\hr{5}\cr
 }}}\\
\hbox{$k'$ in $i$th row}\\
\end{array}
\quad
\begin{array}{c}
\ \ov{i}\leq b\\ 
\ \ \longrightarrow\ \ \\
\end{array}
\quad
T_N=\ {\vcenter
{\offinterlineskip
\halign{&\mystrut\vrule#&\mybox{\hss$\scriptstyle#$\hss}\cr
\hr{5}\cr 
&k'&&c&\cr
\hr{5}\cr 
&b&&e&\cr
\hr{5}\cr
&g&&h&\cr 
\hr{5}\cr
 }}}
\label{Eq-psi-N}
\eeq

Following all these moves of $k'$ and thereby completing the action 
of $\psi_{k'}$ one may still
have one or more unprimed pairs $i$ and $\ov{i}$ in the same column.
If so then $\chi_{k'}$ acts on the lowest such pair. This pair must
be such that $\ov{i}$ and $i$ lie in the $(i-1)$th
and $i$th rows, respectively. To see this let them
lie in the $j-1$th and $j$th rows with $j\leq i$ by
virtue of T$\ov3$. Then following
the action of $\psi_{k'}$ all unprimed entries below $i$ are unpaired
and must lie in their own row, whether they are barred or unbarred,
by virtue of the argument already given in the $\theta_{k'}$ case. 
Thus the entry immediately below $i$ is either $j+1$, $\ov{j+1}$,
$k'$ or $\ov{k}'$. In all four cases the condition T2 
and the rule for shifting $k'$ and $\ov{k}'$ as far north as 
is consistent with T$\ov3$ imply that $i\leq j$. Since $j\leq i$ it 
follows that $j=i$, as claimed. The action of
$\chi_{k'}$ is then to map this vertical pair
$\ov{i}\,i$ with $i$ in the $i$th row to $k'\ov{k}'$. 
\beq
\chi_{k'}:\ \
\begin{array}{c}
T_0=\ {\vcenter
{\offinterlineskip
\halign{&\mystrut\vrule#&\mybox{\hss$\scriptstyle#$\hss}\cr
\hr{5}\cr 
&\ov{i}&&b&\cr
\hr{5}\cr
&i&&d&\cr 
\hr{5}\cr
 }}}\\
\hbox{$i$ in $i$th row}\\
\ov{i}\leq b<d\\
\end{array}
\quad
\ \ \longrightarrow\ \ 
\quad
T=\ {\vcenter
{\offinterlineskip
\halign{&\mystrut\vrule#&\mybox{\hss$\scriptstyle#$\hss}\cr
\hr{5}\cr 
&k'&&b&\cr
\hr{5}\cr
&\ov{k}'&&d&\cr 
\hr{5}\cr
 }}}
\label{Eq-chi-iikk}
\eeq
Then moving the new $k'$ and $\ov{k}'$ 
as far north as permitted by QD1=PD1, the process is repeated for any
other vertical pair $\ov{j}j$ in the $k$th column with
$j$ in the $j$th row, until all such pairs are eliminated. The
result of this is that every unprimed entry $i$ or $\ov{i}$
in the $k$th column lies in the $i$th row, with all
other entries in this column equal to $k'$ or $\ov{k}'$.

Repeating this procedure for all $k=1,2,\ldots,n$ results, as 
required, in the juxtaposition of a primed tableau
$QD\in{\cal PD}^\delta(n,\ov{n})$ and an umprimed tableau
$T\in{\cal T}^\lambda(n,\ov{n})$. It should be noted that
each indvidual step of the maps constituting $\Phi$
is either a simple interchange or a map from a horizontal
pair $\ov{k}'k'$ to a horizontal pair $i\,\ov{i}$, or a map from
a vertical pair $\ov{i}\,i$ to a vertical pair $k'\ov{k}'$.
In each case the weight of each pair is zero, and the number
of barred entries is one. It follows that such steps are all
weight and barred weight preserving, so that (\ref{Eq-Phi})
is always satisfied.

To show that the map $\Phi$ from $QST\in{\cal QST}^\mu(n,\ov{n})$
to $(QD,T)$ with $QD\in{\cal PD}^\delta(n,\ov{n})$ and 
$T\in{\cal T}^\lambda(n,\ov{n})$ is bijective it is sufficient to note 
that each step of the map $\Phi$ may be reversed and that if $\Phi^{-1}$
is defined by such a reversal of each step, then the action
of $\Phi^{-1}$ on the juxtaposition of any $QD\in{\cal PD}^\delta(n,\ov{n})$ 
and any $T\in{\cal T}^\lambda(n,\ov{n})$ always leads to some 
$QST\in{\cal QST}^\mu(n,\ov{n})$.

To see this note that since 
$\Phi=\phi_{n'}\circ\phi_{\ov{n}'}\circ\cdots
\circ\phi_{2'}\circ\phi_{\ov2'}\circ\phi_{1'}\circ\phi_{\ov1'}$
with $\phi_{k'}=\chi_{k'}\circ\psi_{k'}$,
then $\Phi^{-1}=\phi^{-1}_{\ov1'}\circ\phi^{-1}_{1'}\circ
\phi^{-1}_{\ov2'}\circ\phi^{-1}_{2'}\circ\cdots\circ
\phi^{-1}_{\ov{n}'}\circ\phi^{-1}_{n'}$ with
$\phi^{-1}_{k'}=\psi^{-1}_{k'}\circ\chi^{-1}_{k'}$.
The action of the $\phi^{-1}_{\ov{k}'}$ is defined
by the action of $\theta^{-1}_{k'}$ illustrated in (\ref{Eq-jeu-ES})
with $k'$ replaced by $\ov{k}'$, while that of $\psi^{-1}_{k'}$ 
coincides precisely with that of $\theta^{-1}_{k'}$ together
with additional transformations of the type
\beq
\psi^{-1}_{k'}:\ \ 
\begin{array}{c}
T_0=\ {\vcenter
{\offinterlineskip
\halign{&\mystrut\vrule#&\mybox{\hss$\scriptstyle#$\hss}\cr
\hr{7}\cr 
&a&&b&&c&\cr
\hr{7}\cr 
&i&&\ov{i}&&e&\cr
\hr{7}\cr
&f&&g&&h&\cr 
\hr{7}\cr
 }}}\\
\hbox{$i$ in the $i$th row}\\
\end{array}
\quad
\begin{array}{c} 
\ b<\ov{i}\\
\ \ \longrightarrow\ \ \\ 
\end{array} 
T_E=\ {\vcenter
{\offinterlineskip
\halign{&\mystrut\vrule#&\mybox{\hss$\scriptstyle#$\hss}\cr
\hr{7}\cr 
&a&&b&&c&\cr
\hr{7}\cr 
&\ov{k}'&&k'&&e&\cr
\hr{7}\cr
&f&&g&&h&\cr 
\hr{7}\cr
 }}}
\label{Eq-psi-Ekk}
\eeq
where, as indicated, the pair $i\,\ov{i}$ lies in the $i$th 
row of the tableau.
The fact that we necessarily have $b<\ov{i}$ ensures that
this is the exact inverse of the map $\phi_{k'}$
that takes $T_0$ to $T_W$ in (\ref{Eq-phi-NWkk}). 
If $a=\ov{k}'$ then the above gives the inverse
of the map $\phi_{k'}$ that takes $T_0$ to $T_W$
in (\ref{Eq-phi-W}), while if $a=k'$ then 
the inverse of (\ref{Eq-psi-Wkkk}) takes the form
\beq
\psi^{-1}_{k'}:\ \
\begin{array}{c}
T_0=\ {\vcenter
{\offinterlineskip
\halign{&\mystrut\vrule#&\mybox{\hss$\scriptstyle#$\hss}\cr
\hr{7}\cr 
&k'&&b&&c&\cr
\hr{7}\cr 
&i&&\ov{i}&&e&\cr
\hr{7}\cr
&f&&g&&h&\cr 
\hr{7}\cr
 }}}\\
\hbox{$i$ in the $i$th row}\\
\end{array}
\quad
\begin{array}{c} \
\ b<\ov{i}\\ 
\ \ \longrightarrow\ \ \\
\end{array}
\quad 
T_E=\ {\vcenter
{\offinterlineskip
\halign{&\mystrut\vrule#&\mybox{\hss$\scriptstyle#$\hss}\cr
\hr{7}\cr 
&k'&&b&&c&\cr
\hr{7}\cr 
&\ov{k}'&&k'&&e&\cr
\hr{7}\cr
&f&&g&&h&\cr 
\hr{7}\cr
 }}}
\label{Eq-psi-Ekkk}
\eeq
where it is important to note that when applying $\psi^{-1}_{k'}$
the map of any $i\,\ov{i}$ to 
$\ov{k}'k'$ in the $i$th row must be carried out before moving 
any $k'$s that may appear in the tableau higher than the $i$th row.

Similarly, the basic action of the inverse of $\chi_{k'}$, 
which must be applied before $\psi^{-1}_{k'}$ comes into play, 
lowers $k'$s in the $k$th column in accordance with
\beq
\chi^{-1}_{k'}:\ \ \
T_0=\ {\vcenter
{\offinterlineskip
\halign{&\mystrut\vrule#&\mybox{\hss$\scriptstyle#$\hss}\cr
\hr{5}\cr 
&k'&&b&\cr
\hr{5}\cr
&a&&d&\cr 
\hr{5}\cr
 }}}
\quad \ \ 
\longrightarrow\ \ 
\quad 
T=\ {\vcenter
{\offinterlineskip
\halign{&\mystrut\vrule#&\mybox{\hss$\scriptstyle#$\hss}\cr
\hr{5}\cr 
&a&&b&\cr
\hr{5}\cr
&k'&&d&\cr 
\hr{5}\cr
 }}}
\label{Eq-chi-k}
\eeq
if $a$ is unprimed and $a\leq b$. 
If $a=k'$ then the action of $\chi^{-1}_{k'}$
produces no change, while if $a=\ov{k}'$
then 
\beq
\chi^{-1}_{k'}:\ \
\begin{array}{c}
T_0=\ {\vcenter
{\offinterlineskip
\halign{&\mystrut\vrule#&\mybox{\hss$\scriptstyle#$\hss}\cr
\hr{5}\cr 
&k'&&b&\cr
\hr{5}\cr
&\ov{k}'&&d&\cr 
\hr{5}\cr
 }}}\\
\hbox{$\ov{k}'$ in the $i$th row}\\
\ov{i}\leq b<d
\end{array}
\quad
\ \ \longrightarrow\ \ 
\quad 
T=\ {\vcenter
{\offinterlineskip
\halign{&\mystrut\vrule#&\mybox{\hss$\scriptstyle#$\hss}\cr
\hr{5}\cr 
&\ov{i}&&b&\cr
\hr{5}\cr
&i&&d&\cr 
\hr{5}\cr
 }}}
\label{Eq-chi-kkii}
\eeq
where the constraint $\ov{i}\leq b<d$
is required so as to ensure that $T_0$ satisfies 
the semistandardness
conditions T1 and T2. Of course if $b<\ov{i}$ then
$\chi^{-1}_{k'}$ cannot map as above, so it simply leaves
$T_0$ unchanged but then subject to $\phi^{-1}_{k'}$
which acts as follows
\beq
\phi^{-1}_{k'}:\ \ \
T=\ {\vcenter
{\offinterlineskip
\halign{&\mystrut\vrule#&\mybox{\hss$\scriptstyle#$\hss}\cr
\hr{5}\cr 
&k'&&b&\cr
\hr{5}\cr
&\ov{k}'&&d&\cr 
\hr{5}\cr
 }}}
\quad 
\ \ \longrightarrow\ \ 
\quad 
T_E=\ {\vcenter
{\offinterlineskip
\halign{&\mystrut\vrule#&\mybox{\hss$\scriptstyle#$\hss}\cr
\hr{5}\cr 
&b&&k'&\cr
\hr{5}\cr
&\ov{k}'&&d&\cr 
\hr{5}\cr
 }}}
\label{Eq-phi-kE}
\eeq
in accordance with the inverse of $\theta_{k'}$
applied in the case $b<\ov{i}$.
In applying $\chi^{-1}_{k'}$
as in (\ref{Eq-chi-k}) and (\ref{Eq-chi-kkii}) one
starts the action with the most northerly $k'$
and the most northerly vertical pairs $k'\ov{k}'$,
respectively, working down the $k$th column. 
This is in contrast to the subsequent action of $\psi^{1}_{k'}$
which acts on the most southerly $k'$ first.

Finally, it is necessary to show that all $k'$s and $\ov{k}'$s are mapped
by $\phi^{-1}_{k'}$ and $\phi^{-1}_{\ov{k}'}$ to endpoints
consistent with the conditions QST1-4 and QST$\ov5$
on all primed tableaux
$QST\in{\cal QST}^\mu(n,\ov{n})$. First it should be noted that 
QST3 and QST$\ov5$ are satisfied throughout the application of 
$\Phi^{-1}$. Furthermore, a violation of QST1 or QST2 by 
any $k'$ or $\ov{k}'$ simply means that application of 
$\phi^{-1}_{k'}$ or $\phi^{-1}_{\ov{k}'}$, respectively, has not been 
completed. The argument used in Section~\ref{sec:gl-bijection}
regarding reverse paths staying strictly north of one another in any
given column, is sufficient to ensure that no two identical primed 
entries may appear in the same row, thereby ensuring
that the final condition QST3 is always satisfied. Thus
the image of $\Phi^{-1}$ of any pair $(QD,T)$ with
$QD\in{\cal PD}^\delta(n,\ov{n})$ and 
$T\in{\cal T}^\lambda(n,\ov{n})$ 
is some $QST\in{\cal QST}^\mu(n,\ov{n})$.

This implies that the original map $\Phi$ is bijective. Since it 
is also both weight and barred weight preserving, this completes
the proof of Theorem~\ref{The-Qbar}.

It should be pointed out that, unlike the $gl(n)$ case, a corresponding 
result does not apply to $PST\in{\cal PST}^\mu(n,\ov{n})$ in 
the $sp(2n)$ case because of the necessity of using (\ref{Eq-chi-iikk}).
If used, as will sometimes be necessary, in the case $i=k$ it results
in a primed entry $\ov{k}'$ appearing on the main diagonal, in direct
violation of PD3.

\subsection{Example}
\label{subsec:sp-example}

The bijection $\Phi$ is illustrated by the following map for $n=5$ and 
$\mu=(9,7,6,2,1)$:
\beq
PST=\
{\vcenter
 {\offinterlineskip
 \halign{&\mystrut\vrule#&\mybox{\hss$\scriptstyle#$\hss}\cr
  \hr{19}\cr
  &\ov1&&1'&&\ov2'&&2'&&\ov3'&&3&&\ov4&&\ov4&&5&\cr
  \hr{19}\cr
  \omit& &&\ov2'&&\ov2&&\ov3&&\ov3&&4'&&4&&4&\cr
  \nr{2}&\hr{15}\cr
  \omit& &\omit& &&\ov3&&3&&\ov4&&5'&&5&&5&\cr
  \nr{4}&\hr{13}\cr
  \omit& &\omit& &\omit& &&4&&4&\cr
  \nr{6}&\hr{5}\cr
  \omit& &\omit& &\omit& &\omit& &&5'&\cr
  \nr{8}&\hr{3}\cr
 }}}
\ \ \longrightarrow\ \ 
PD\ =\ \
{\vcenter
 {\offinterlineskip
 \halign{&\mystrut\vrule#&\mybox{\hss$\scriptstyle#$\hss}\cr
  \hr{11}\cr
  &1'&&1&&\ov3'&&4'&&\ov1&\cr
  \hr{11}\cr
  \omit& &&\ov2'&&\ov2&&4'&&5'&\cr
  \nr{2}&\hr{9}\cr
  \omit& &\omit& &&\ov3&&\ov4'&&\ov3&\cr
  \nr{4}&\hr{7}\cr
  \omit& &\omit& &\omit& &&4&&4&\cr
  \nr{6}&\hr{5}\cr
  \omit& &\omit& &\omit& &\omit& &&5'&\cr
  \nr{8}&\hr{3}\cr
 }}}
\ \ ,\ \
T\ =\ \ 
{\vcenter
 {\offinterlineskip
 \halign{&\mystrut\vrule#&\mybox{\hss$\scriptstyle#$\hss}\cr
  \hr{9}\cr
  &\ov1&&\ov4&&\ov4&&5&\cr
  \hr{9}\cr
  &3&&4&&4&\cr
  \hr{7}\cr
  &\ov4&&5&&5&\cr
  \hr{7}\cr
  }}}
\label{Eq-Ex-sp}
\eeq

Once again we indicate by means of boldface entries both the paths traced
out by elements $\ov{k}'$ and $k'$ as they move to the $k$th column 
under the action of $\phi_{\ov{k}'}$ and $\psi_k'$, respectively, as
well as annihilations and creation of $\ov{k}'k'$ pairs under
$\psi_{k'}$ and $\chi_{k'}$, respectively.

There are no $\ov{1}'$s, so we first move the single $1'$ under the action of 
$\psi_{1'}$ as shown:
\beq
{\vcenter
 {\offinterlineskip
 \halign{&\mystrut\vrule#&\mybox{\hss$\scriptstyle#$\hss}\cr
  \hr{19}\cr
  &\ov{\bf1}&&\bf{1}'&&\ov2'&&2'&&\ov3'&&3&&\ov4&&\ov4&&5&\cr
  \hr{19}\cr
  \omit& &&\ov2'&&\ov2&&\ov3&&\ov3&&4'&&4&&4&\cr
  \nr{2}&\hr{15}\cr
  \omit& &\omit& &&\ov3&&3&&\ov4&&5'&&5&&5&\cr
  \nr{4}&\hr{13}\cr
  \omit& &\omit& &\omit& &&4&&4&\cr
  \nr{6}&\hr{5}\cr
  \omit& &\omit& &\omit& &\omit& &&5'&\cr
  \nr{8}&\hr{3}\cr
 }}}
\ \ \longrightarrow\ \ 
{\vcenter
 {\offinterlineskip
 \halign{&\mystrut\vrule#&\mybox{\hss$\scriptstyle#$\hss}\cr
  \hr{19}\cr
  &\bf{1}'&&\ov{\bf1}&&\ov2'&&2'&&\ov3'&&3&&\ov4&&\ov4&&5&\cr
  \hr{19}\cr
  \omit& &&\ov2'&&\ov2&&\ov3&&\ov3&&4'&&4&&4&\cr
  \nr{2}&\hr{15}\cr
  \omit& &\omit& &&\ov3&&3&&\ov4&&5'&&5&&5&\cr
  \nr{4}&\hr{13}\cr
  \omit& &\omit& &\omit& &&4&&4&\cr
  \nr{6}&\hr{5}\cr
  \omit& &\omit& &\omit& &\omit& &&5'&\cr
  \nr{8}&\hr{3}\cr
 }}}
\eeq
The application of $\phi_{\ov2}'$ then gives
\beq
{\vcenter
 {\offinterlineskip
 \halign{&\mystrut\vrule#&\mybox{\hss$\scriptstyle#$\hss}\cr
  \hr{19}\cr
  &1'&&\ov{\bf1}&&\ov{\bf2}'&&2'&&\ov3'&&3&&\ov4&&\ov4&&5&\cr
  \hr{19}\cr
  \omit& &&\ov2'&&\ov2&&\ov3&&\ov3&&4'&&4&&4&\cr
  \nr{2}&\hr{15}\cr
  \omit& &\omit& &&\ov3&&3&&\ov4&&5'&&5&&5&\cr
  \nr{4}&\hr{13}\cr
  \omit& &\omit& &\omit& &&4&&4&\cr
  \nr{6}&\hr{5}\cr
  \omit& &\omit& &\omit& &\omit& &&5'&\cr
  \nr{8}&\hr{3}\cr
 }}}
\ \ \longrightarrow\ \ 
{\vcenter
 {\offinterlineskip
 \halign{&\mystrut\vrule#&\mybox{\hss$\scriptstyle#$\hss}\cr
  \hr{19}\cr
  &1'&&\ov{\bf2}'&&\ov{\bf1}&&2'&&\ov3'&&3&&\ov4&&\ov4&&5&\cr
  \hr{19}\cr
  \omit& &&\ov2'&&\ov2&&\ov3&&\ov3&&4'&&4&&4&\cr
  \nr{2}&\hr{15}\cr
  \omit& &\omit& &&\ov3&&3&&\ov4&&5'&&5&&5&\cr
  \nr{4}&\hr{13}\cr
  \omit& &\omit& &\omit& &&4&&4&\cr
  \nr{6}&\hr{5}\cr
  \omit& &\omit& &\omit& &\omit& &&5'&\cr
  \nr{8}&\hr{3}\cr
 }}}
\eeq
where there is no possibility of moving the lower $\ov{2}'$. Then
the application of $\psi_{2'}$ on the only $2'$ involves first a 
transposition and then the replacement of the
resulting horizontal pair $\ov2'\,2'$ in the first row by $1\,\ov1$:
\beq 
{\vcenter
 {\offinterlineskip
 \halign{&\mystrut\vrule#&\mybox{\hss$\scriptstyle#$\hss}\cr
  \hr{19}\cr
  &1'&&\ov{2}'&&\ov{\bf1}&&\bf{2}'&&\ov3'&&3&&\ov4&&\ov4&&5&\cr
  \hr{19}\cr
  \omit& &&\ov2'&&\ov2&&\ov3&&\ov3&&4'&&4&&4&\cr
  \nr{2}&\hr{15}\cr
  \omit& &\omit& &&\ov3&&3&&\ov4&&5'&&5&&5&\cr
  \nr{4}&\hr{13}\cr
  \omit& &\omit& &\omit& &&4&&4&\cr
  \nr{6}&\hr{5}\cr
  \omit& &\omit& &\omit& &\omit& &&5'&\cr
  \nr{8}&\hr{3}\cr
 }}}
\ \ \longrightarrow\ \ 
{\vcenter
 {\offinterlineskip
 \halign{&\mystrut\vrule#&\mybox{\hss$\scriptstyle#$\hss}\cr
  \hr{19}\cr
  &1'&&\ov{2}'&&\bf{2}'&&\ov{\bf1}&&\ov3'&&3&&\ov4&&\ov4&&5&\cr
  \hr{19}\cr
  \omit& &&\ov2'&&\ov2&&\ov3&&\ov3&&4'&&4&&4&\cr
  \nr{2}&\hr{15}\cr
  \omit& &\omit& &&\ov3&&3&&\ov4&&5'&&5&&5&\cr
  \nr{4}&\hr{13}\cr
  \omit& &\omit& &\omit& &&4&&4&\cr
  \nr{6}&\hr{5}\cr
  \omit& &\omit& &\omit& &\omit& &&5'&\cr
  \nr{8}&\hr{3}\cr
 }}}
\ \ =\ \ 
{\vcenter
 {\offinterlineskip
 \halign{&\mystrut\vrule#&\mybox{\hss$\scriptstyle#$\hss}\cr
  \hr{19}\cr
  &1'&&\ov{\bf2}'&&\bf{2}'&&\ov1&&\ov3'&&3&&\ov4&&\ov4&&5&\cr
  \hr{19}\cr
  \omit& &&\ov2'&&\ov2&&\ov3&&\ov3&&4'&&4&&4&\cr
  \nr{2}&\hr{15}\cr
  \omit& &\omit& &&\ov3&&3&&\ov4&&5'&&5&&5&\cr
  \nr{4}&\hr{13}\cr
  \omit& &\omit& &\omit& &&4&&4&\cr
  \nr{6}&\hr{5}\cr
  \omit& &\omit& &\omit& &\omit& &&5'&\cr
  \nr{8}&\hr{3}\cr
 }}}
\ \ \longrightarrow\ \ 
{\vcenter
 {\offinterlineskip
 \halign{&\mystrut\vrule#&\mybox{\hss$\scriptstyle#$\hss}\cr
  \hr{19}\cr
  &1'&&\bf1&&\ov{\bf1}&&\ov1&&\ov3'&&3&&\ov4&&\ov4&&5&\cr
  \hr{19}\cr
  \omit& &&\ov2'&&\ov2&&\ov3&&\ov3&&4'&&4&&4&\cr
  \nr{2}&\hr{15}\cr
  \omit& &\omit& &&\ov3&&3&&\ov4&&5'&&5&&5&\cr
  \nr{4}&\hr{13}\cr
  \omit& &\omit& &\omit& &&4&&4&\cr
  \nr{6}&\hr{5}\cr
  \omit& &\omit& &\omit& &\omit& &&5'&\cr
  \nr{8}&\hr{3}\cr
 }}}
\eeq
The single $\ov3'$ is moved as shown to the $3$rd
column under the action of $\psi_{3'}$:
\beq
{\vcenter
 {\offinterlineskip
 \halign{&\mystrut\vrule#&\mybox{\hss$\scriptstyle#$\hss}\cr
  \hr{19}\cr
  &1'&&1&&\ov{\bf1}&&\ov{\bf1}&&\bf\ov3'&&3&&\ov4&&\ov4&&5&\cr
  \hr{19}\cr
  \omit& &&\ov2'&&\ov2&&\ov3&&\ov3&&4'&&4&&4&\cr
  \nr{2}&\hr{15}\cr
  \omit& &\omit& &&\ov3&&3&&\ov4&&5'&&5&&5&\cr
  \nr{4}&\hr{13}\cr
  \omit& &\omit& &\omit& &&4&&4&\cr
  \nr{6}&\hr{5}\cr
  \omit& &\omit& &\omit& &\omit& &&5'&\cr
  \nr{8}&\hr{3}\cr
 }}}
\ \ \longrightarrow\ \
{\vcenter
 {\offinterlineskip
 \halign{&\mystrut\vrule#&\mybox{\hss$\scriptstyle#$\hss}\cr
  \hr{19}\cr
  &1'&&1&&\bf\ov3'&&\ov{\bf1}&&\ov{\bf1}&&3&&\ov4&&\ov4&&5&\cr
  \hr{19}\cr
  \omit& &&\ov2'&&\ov2&&\ov3&&\ov3&&4'&&4&&4&\cr
  \nr{2}&\hr{15}\cr
  \omit& &\omit& &&\ov3&&3&&\ov4&&5'&&5&&5&\cr
  \nr{4}&\hr{13}\cr
  \omit& &\omit& &\omit& &&4&&4&\cr
  \nr{6}&\hr{5}\cr
  \omit& &\omit& &\omit& &\omit& &&5'&\cr
  \nr{8}&\hr{3}\cr
 }}}
\eeq
Next, the single $4'$ is moved under $\phi_{{4}'}$ as follows:
\beq
{\vcenter
 {\offinterlineskip
 \halign{&\mystrut\vrule#&\mybox{\hss$\scriptstyle#$\hss}\cr
  \hr{19}\cr
  &1'&&1&&\ov3'&&\ov{\bf1}&&\ov{\bf1}&&\bf3&&\ov4&&\ov4&&5&\cr
  \hr{19}\cr
  \omit& &&\ov2'&&\ov2&&\ov3&&\ov3&&\bf4'&&4&&4&\cr
  \nr{2}&\hr{15}\cr
  \omit& &\omit& &&\ov3&&3&&\ov4&&5'&&5&&5&\cr
  \nr{4}&\hr{13}\cr
  \omit& &\omit& &\omit& &&4&&4&\cr
  \nr{6}&\hr{5}\cr
  \omit& &\omit& &\omit& &\omit& &&5'&\cr
  \nr{8}&\hr{3}\cr
 }}}
\ \ \longrightarrow\ \
{\vcenter
 {\offinterlineskip
 \halign{&\mystrut\vrule#&\mybox{\hss$\scriptstyle#$\hss}\cr
  \hr{19}\cr
  &1'&&1&&\ov3'&&\bf4'&&\ov{\bf1}&&\ov{\bf1}&&\ov4&&\ov4&&5&\cr
  \hr{19}\cr
  \omit& &&\ov2'&&\ov2&&\ov3&&\ov3&&\bf3&&4&&4&\cr
  \nr{2}&\hr{15}\cr
  \omit& &\omit& &&\ov3&&3&&\ov4&&5'&&5&&5&\cr
  \nr{4}&\hr{13}\cr
  \omit& &\omit& &\omit& &&4&&4&\cr
  \nr{6}&\hr{5}\cr
  \omit& &\omit& &\omit& &\omit& &&5'&\cr
  \nr{8}&\hr{3}\cr
 }}}
\eeq
There are no $\ov4'$s. However, the $4$th column contains the pair
$\ov3\,3$ which must be replaced by $4'\,\ov4'$ under the 
action of $\chi_{4'}$:
\beq
{\vcenter
 {\offinterlineskip
 \halign{&\mystrut\vrule#&\mybox{\hss$\scriptstyle#$\hss}\cr
  \hr{19}\cr
  &1'&&1&&\ov3'&&4'&&\ov1&&\ov1&&\ov4&&\ov4&&5&\cr
  \hr{19}\cr
  \omit& &&\ov2'&&\ov2&&\ov{\bf3}&&\ov3&&3&&4&&4&\cr
  \nr{2}&\hr{15}\cr
  \omit& &\omit& &&\ov3&&\bf3&&\ov4&&5'&&5&&5&\cr
  \nr{4}&\hr{13}\cr
  \omit& &\omit& &\omit& &&4&&4&\cr
  \nr{6}&\hr{5}\cr
  \omit& &\omit& &\omit& &\omit& &&5'&\cr
  \nr{8}&\hr{3}\cr
 }}}
\ \ \longrightarrow\ \
{\vcenter
 {\offinterlineskip
 \halign{&\mystrut\vrule#&\mybox{\hss$\scriptstyle#$\hss}\cr
  \hr{19}\cr
  &1'&&1&&\ov3'&&4'&&\ov1&&\ov1&&\ov4&&\ov4&&5&\cr
  \hr{19}\cr
  \omit& &&\ov2'&&\ov2&&\bf4'&&\ov3&&3&&4&&4&\cr
  \nr{2}&\hr{15}\cr
  \omit& &\omit& &&\ov3&&\ov{\bf4}'&&\ov4&&5'&&5&&5&\cr
  \nr{4}&\hr{13}\cr
  \omit& &\omit& &\omit& &&4&&4&\cr
  \nr{6}&\hr{5}\cr
  \omit& &\omit& &\omit& &\omit& &&5'&\cr
  \nr{8}&\hr{3}\cr
 }}}
\eeq
There are no $\ov{5}'$s and it is then important to notice that one 
does not replace the vertical pair $\ov4\,4$ in the 
$5$th column by $5'\,\ov{5}'$ under the action of $\chi_{5'}$
because one must first apply $\psi_{5'}$ to the two $5'$s. In any
case the premature action of $\chi_{5'}$ would lead to two $5'$s in
the same row, which is forbidden. Instead, the algorithm dictates
that one first acts on the two $5'$s with $\psi_{5'}$. The uppermost 
$5'$ is moved into the $5$th column and up that column until it is 
just below the entry $\ov1$ which cannot be moved into the second
row. This leaves both the $\ov3$ and $4$ in their own rows in the
$5$th column, and no $\ov4\,4$ pair. 
\beq
 {\vcenter
 {\offinterlineskip
 \halign{&\mystrut\vrule#&\mybox{\hss$\scriptstyle#$\hss}\cr
  \hr{19}\cr
  &1'&&1&&\ov3'&&4'&&\ov1&&\ov1&&\ov4&&\ov4&&5&\cr
  \hr{19}\cr
  \omit& &&\ov2'&&\ov2&&4'&&\ov{\bf3}&&3&&4&&4&\cr
  \nr{2}&\hr{15}\cr
  \omit& &\omit& &&\ov3&&\ov4'&&\ov{\bf4}&&\bf{5}'&&5&&5&\cr
  \nr{4}&\hr{13}\cr
  \omit& &\omit& &\omit& &&4&&4&\cr
  \nr{6}&\hr{5}\cr
  \omit& &\omit& &\omit& &\omit& &&5'&\cr
  \nr{8}&\hr{3}\cr
 }}}
\ \ \longrightarrow\ \
{\vcenter
 {\offinterlineskip
 \halign{&\mystrut\vrule#&\mybox{\hss$\scriptstyle#$\hss}\cr
  \hr{19}\cr
  &1'&&1&&\ov3'&&4'&&\ov1&&\ov1&&\ov4&&\ov4&&5&\cr
  \hr{19}\cr
  \omit& &&\ov2'&&\ov2&&4'&&\bf5'&&3&&4&&4&\cr
  \nr{2}&\hr{15}\cr
  \omit& &\omit& &&\ov3&&\ov4'&&\ov{\bf3}&&\ov{\bf4}&&5&&5&\cr
  \nr{4}&\hr{13}\cr
  \omit& &\omit& &\omit& &&4&&4&\cr
  \nr{6}&\hr{5}\cr
  \omit& &\omit& &\omit& &\omit& &&5'&\cr
  \nr{8}&\hr{3}\cr
 }}}
\eeq 
The second $5'$ does not move under the action of $\psi_{5'}$ 
since the $4$ immediately above it cannot move into the $5$th row.
The final result can then be seen to be, as claimed, the juxtaposition 
of a primed tableaux $PD\in{\cal PD}^{54321}(5,\ov{5})$
and an unprimed tableaux $T\in{\cal T}^{433}(5,\ov{5})$:
\beq
{\vcenter
 {\offinterlineskip
 \halign{&\mystrut\vrule#&\mybox{\hss$\scriptstyle#$\hss}\cr
  \hr{19}\cr
  &1'&&1&&\ov3'&&4'&&\ov1&&\ov1&&\ov4&&\ov4&&5&\cr
  \hr{19}\cr
  \omit& &&\ov2'&&\ov2&&4'&&5'&&3&&4&&4&\cr
  \nr{2}&\hr{15}\cr
  \omit& &\omit& &&\ov3&&\ov4'&&\ov3&&\ov4&&5&&5&\cr
  \nr{4}&\hr{13}\cr
  \omit& &\omit& &\omit& &&4&&4&\cr
  \nr{6}&\hr{5}\cr
  \omit& &\omit& &\omit& &\omit& &&5'&\cr
  \nr{8}&\hr{3}\cr
 }}}
\ \ \equiv\ \
{\vcenter
 {\offinterlineskip
 \halign{&\mystrut\vrule#&\mybox{\hss$\scriptstyle#$\hss}\cr
  \hr{11}\cr
  &1'&&1&&\ov3'&&4'&&\ov1&\cr
  \hr{11}\cr
  \omit& &&\ov2'&&\ov2&&4'&&5'&\cr
  \nr{2}&\hr{9}\cr
  \omit& &\omit& &&\ov3&&\ov4'&&\ov3&\cr
  \nr{4}&\hr{7}\cr
  \omit& &\omit& &\omit& &&4&&4&\cr
  \nr{6}&\hr{5}\cr
  \omit& &\omit& &\omit& &\omit& &&5'&\cr
  \nr{8}&\hr{3}\cr
 }}}
\ \ \cdot\ \ 
{\vcenter
 {\offinterlineskip
 \halign{&\mystrut\vrule#&\mybox{\hss$\scriptstyle#$\hss}\cr
  \hr{9}\cr
  &\ov1&&\ov4&&\ov4&&5&\cr
  \hr{9}\cr
  &3&&4&&4&\cr
  \hr{7}\cr
  &\ov4&&5&&5&\cr
  \hr{7}\cr
  }}}
\eeq

\subsection{Corollaries}
\label{subsec:sp-corollaries}

By associating $x_k$, $t^2\,\ov{x}_k$, $y_k$ and $t^2\,\ov{y}_k$ to each entry 
$k$, $\ov{k}$, $k'$ and $\ov{k}'$, respectively, in the various 
tableaux $QST$, $QD$ and $T$ appearing in Theorem~\ref{The-Qbar}
we immediately have the following corollary.

\begin{Corollary}
Let $\mu=\lambda+\delta$ be a strict partition of length $\ell(\mu)=n$,
with $\lambda$ a partition of length $\ell(\lambda)\leq n$ and 
$\delta=(n,n-1,\ldots,1)$. 
\beq
   \sum_{QST\in{\cal QST}^\mu(n)}\ t^{2\bar(QST)}\ (\x/\y)^{\wgt(QST)}
 \ =\ \sum_{QD\in{\cal QD}^\delta(n)}\ t^{2\bar(QD)}\ (\x/\y)^{\wgt(QD)} 
    \ \sum_{T\in{\cal T}^\lambda(n)}\ t^{2\bar(T)}\ \x^{\wgt(T)}\,. 
\label{Eq-QST-QD-T}
\eeq
\label{Cor-PQbar}
\end{Corollary}

Thanks to the definition of $Q(\x/\y;t)$
given in (\ref{Eq-PQx}), the identity (\ref{Eq-QDprodbar}) 
and the combinatorial definition of the $t$-deformation 
$sp_\lambda(\x;t)$ of symplectic characters given in 
(\ref{Eq-Spt}), the above 
result is precisely our second main result
Proposition~\ref{Prop-PQsp}. 

Other corollaries follow as special cases of these results.
Setting $\lambda=0$ we obtain 
\begin{equation}
Q_{\delta}(\x/\y;t)=\prod_{1\leq i\leq j \leq n} 
(x_i + t^2\,\ov{x}_i+y_j+t^2\,\ov{y}_j).
\end{equation}

On the other hand the case
$\y=t\x=(tx_1,tx_2,\ldots,tx_n)$ of (\ref{Eq-MainResult-sp}) is equivalent to
the $t$-deformation of Weyl's denominator formula (\ref{Eq-HK}) for the
Lie algebra $sp(2n)$ derived elsewhere~\cite{HK02,HK03} by much more 
circuitous means.

\begin{Corollary}
\begin{equation}
\begin{array}{l}
\prod_{i=1}^n\ (x_i+t\,\ov{x}_i) \ 
\prod_{1\leq i<j \leq n} (x_i+t^2\,\ov{x}_i+t\,x_j+t\,\ov{x}_j)\ 
sp_{\lambda}(\x;t)\\ \\
\qquad\qquad =\ \ 
\sum_{ST\in {\cal ST}^\mu(n\ov{n})} 
\ t^{\var(ST)+\bar(ST)} (1+t)^{\str(ST)-n}\ \x^{\wgt(ST)}.\\
\end{array}
\label{Eq-Cor-sp}
\end{equation}
\end{Corollary}

\noindent{\bf Proof}\ \ 
First it should be noted that
\beq
\prod_{1\leq i\leq j \leq n} (x_i+t^2\,\ov{x}_i+t\,x_j+t\,\ov{x}_j)\ 
=\prod_{i=1}^n\ (1+t)\,(x_i+t\,\ov{x}_i)\ 
\prod_{1\leq i < j \leq n} (1+t\,\ov{x}_i\,\ov{x}_j)
  (1+t\,\ov{x}_i\,\ov{x}_j)\ ,
\label{Eq-prods}
\eeq
which includes a factor $(1+t)^n$.

Then it suffices to recognise that deleting primes from the
entries $k'$ and $\ov{k}'$ in $QST\in{\cal QST}^\mu(n,\ov{n})$
gives a symplectic shifted tableau $ST\in{\cal ST}^\mu(n,\ov{n})$
with a factor of $t$ arising from each primed entry since $y_k=tx_k$
and $t^2y_k^{-1}=tx_k^{-1}$. Additional factors of $t^2$ arise from
each barred entry $k$ since these are associated with $t^2x_k^{-1}$.
Compulsorily primed entries in each primed shifted tableau $QST$ 
corresponding to a fixed shifted tableau $ST$ appear once and only 
once in each row of each connected strip of identical entries, 
whether barred or unbarrred, while the leftmost lowest entry of 
each connected strip may each be primed or unprimed. To summarise,
each connected component of a $k$-strip gives rise to a factor
$(1+t)\,t^{\row_k-1}$, where $\row_k$ is the number of rows
of the $k$-strip component, and each component of a $\ov{k}$-strip
gives rise to a factor of $(1+t^{-1})\,t^{2\bar_{\ov{k}}-\row_{\ov{k}}+1}
=(1+t)\,t^{\bar_{\ov{k}}+\col_{\ov{k}}-1}$ where $\bar_{\ov{k}}$ 
is the length of 
the component of the $\ov{k}$-strip, while $\row_{\ov{k}}$
and $\col_{\ov{k}}$ are the numbers of rows and columns, respectively,
that it occupies, with $\bar_{\ov{k}}=\row_{\ov{k}}+\col_{\ov{k}}-1$.
Combining all these factors for $k=1,2,\ldots,n$ gives 
$(1+t)^{\str(ST)}\,t^{\var(ST)+\bar(ST)}$, as required, since 
as defined earlier $\var(ST)=\sum (\row_k+\col_{\ov{k}}-1)$
where the sum is over all connected components of all $k$ and $\ov{k}$
strips, for all $k$, and $\bar(ST)=\sum \bar_{\ov{k}}$ is the total
number of barred entries in $ST$.

Another significant corollary involves a link with 
U-turn alternating sign matrices. This is
provided in Section 5.

\section{Connection to Alternating Sign Matrices}
\label{sec:asm}

\subsection{$gl(n)$ case}\ \
\label{subsec:gl-case}

In this section we show how to map from shifted tableaux,
$ST$, to alternating sign matrices.
Using the analogous relationship for primed shifted tableaux, $PST$, 
a result of Chapman \cite{C01} is a straightforward
consequence of Theorem~\ref{The-PQ}.

An alternating sign matrix (ASM) is an $n \times n$ matrix filled with 
$0$'s, $1$'s, and $-1$'s such that the first and last nonzero entries of each
row and column are $1$'s and the nonzero entries within a row or column
alternate in sign.  There is a famous formula, conjectured by Mills,
Robbins, and Rumsey \cite{MRR83} and proved by Zeilberger \cite{Z96}, that
counts the number of ASM of size $n$ as 
$\prod_{j=0}^{n-1}\ (3j+1)! \big/ (n+j)!$.
See also Bressoud \cite{B99}.

We work with a generalisation of ASM called $\mu$--ASM introduced by 
Okada~\cite{O93} that can be associated with semistandard shifted tableaux.  
 The new feature here is that the alternating sign matrix is no longer square. 
Its row sums are all $1$ but its column sums are $1$ or $0$ according
as the column label is or is not a part of some partition $\mu$.
To be more precise,
given a partition $\mu$ with distinct parts
and such that $\ell(\mu)=n$ and $\mu_1 \leq m$ for some 
$m\geq n$, the set ${\cal A}^\mu(n)$ of $\mu$--alternating
sign matrices is the set of $n\times m$ matrices 
$A=(a_{ij})_{1\leq i,j\leq n}$ that satisfy
the following conditions:

\begin{tabular}{rl}
ASM1& $a_{iq} \in \{ -1, 0, 1\}$ for $1 \leq i \leq n, 1 \leq q\leq m$;\\
ASM2& $\sum_{q=p}^{m} a_{iq} \in \{ 0,1\}$ for $ 1 \leq i \leq n, 1\leq p
\leq m$;\\
ASM3& $\sum_{i=j}^{n} a_{iq} \in \{ 0,1\}$ for $1\leq j \leq n, 1 \leq q
\leq m$ \\
ASM4& $\sum_{q=1}^{m} a_{iq} = 1$ for $1 \leq i \leq n$;\\
ASM5& $\sum_{i=1}^{n} a_{iq} =1$ if $q=\mu_j$ for some $j$; or
$\sum_{i=1}^{n} a_{iq} = 0$ otherwise; for $1 \leq q \leq m$.
\end{tabular}

In what follows we also require U-turn alternating sign matrices, UASMs,
and their generalisation $\mu$--UASMs that are associated with
semistandard shifted symplectic tableaux~\cite{HK02,HK03}. In fact
the bijection between semistandard shifted tableaux $ST\in{\cal ST}^\mu(n)$ and
$\mu$--ASMs, $A\in{\cal A}^\mu(n)$, is a special case of a bijection 
between semistandard shifted symplectic tableaux $ST\in{\cal ST}^\mu(n,\ov{n})$
and $\mu$--UASMs $A\in{\cal A}^\mu(n,\ov{n})$~\cite{HK03}.  

Briefly, in the $\mu$--ASM case, we
associate to each semistandard shifted tableaux $ST\in{\cal ST}^\mu(n)$
of shape $\mu$ with $\ell(\mu)=n$ and $\mu_1=m$ an $n \times m$ 
matrix $M(ST)$ filled with the entries from $ST$ together with zeros such
that if there is an $i$ on diagonal $j$ of $ST$ (where the
main diagonal is diagonal 1 and the last box in the first row is in
diagonal $\mu_1=m$), then there is an $i$  in row $i$,
column $j$ of the matrix.  All other entries are zero.

For example, in the case $n=6$ and $\mu=(9,8,6,4,3,1)$ a given semistandard 
shifted tableau $ST$ of shape $\mu$ yields a $6\times9$ matrix, $M(ST)$, as shown:
\beq
ST=\ {\vcenter
{\offinterlineskip
\halign{&\mystrut\vrule#&\mybox{\hss$\scriptstyle#$\hss}\cr
\hr{19}\cr 
&1&&1&&1&&2&&3&&3&&4&&4&&4&\cr
\hr{19}\cr 
\omit&
&&2&&2&&2&&3&&4&&5&&5&&5&\cr
\nr{2}&\hr{17}\cr 
\omit& &\omit&
&&3&&4&&4&&4&&5&&6&\cr 
\nr{4}&\hr{13}\cr 
\omit& &\omit& &\omit&
&&4&&5&&5&&6&\cr
\nr{6}&\hr{9}\cr 
\omit& &\omit& &\omit& &\omit&
&&5&&6&&6&\cr 
\nr{8}&\hr{7}\cr
\omit& &\omit& &\omit& &\omit& &\omit&
&&6&\cr 
\nr{10}&\hr{3}\cr
 }}}\ 
\ \Longrightarrow\ 
M(ST) = \left[\begin{array}{lllllllll}
1 & 1 & 1 & 0 & 0 & 0 & 0 & 0 & 0 \\
2 & 2 & 2 & 2 & 0 & 0 & 0 & 0 & 0 \\
3 & 0 & 0 & 3 & 3 & 3 & 0 & 0 & 0 \\
4 & 4 & 4 & 4 & 4 & 0 & 4 & 4 & 4 \\
5 & 5 & 5 & 0 & 5 & 5 & 5 & 5 & 0 \\
6 & 6 & 6 & 6 & 0 & 6 & 0 & 0 & 0 \\
\end{array}\right]
\eeq

A primed semistandard shifted tableau $PST\in{\cal PST}^\mu(N)$ 
yields a similar matrix $M(PST)$ in the same way:
\beq
PST=\ {\vcenter
{\offinterlineskip
\halign{&\mystrut\vrule#&\mybox{\hss$\scriptstyle#$\hss}\cr
\hr{19}\cr 
&1&&1&&1&&2'&&3'&&3&&4&&4&&4&\cr
\hr{19}\cr 
\omit&
&&2&&2&&2&&3'&&4'&&5'&&5&&5&\cr
\nr{2}&\hr{17}\cr 
\omit& &\omit&
&&3&&4'&&4&&4&&5&&6&\cr 
\nr{4}&\hr{13}\cr 
\omit& &\omit& &\omit&
&&4&&5'&&5&&6'&\cr
\nr{6}&\hr{9}\cr 
\omit& &\omit& &\omit& &\omit&
&&5&&6'&&6&\cr 
\nr{8}&\hr{7}\cr
\omit& &\omit& &\omit& &\omit& &\omit&
&&6&\cr 
\nr{10}&\hr{3}\cr
 }}}\ 
\Longrightarrow
M(PST) = \left[\begin{array}{lllllllll}
1 & 1 & 1 & 0 & 0 & 0 & 0 & 0 & 0 \\
2 & 2 & 2 & 2' & 0 & 0 & 0 & 0 & 0 \\
3 & 0 & 0 & 3' & 3' & 3 & 0 & 0 & 0 \\
4 & 4' & 4 & 4 & 4' & 0 & 4 & 4 & 4 \\
5 & 5' & 5 & 0 & 5 & 5' & 5 & 5 & 0 \\
6 & 6' & 6 & 6' & 0 & 6 & 0 & 0 & 0 \\
\end{array}\right]
\eeq
where as we shall see it is possible to
distinguish various types of entry $0$
as characterised by their set of nearest
non-vanishing neighbours. 

Each of these matrices can be converted into a 
$\mu$--alternating sign matrix by replacing the rightmost 
entry of each continuous sequence of nonzero entries
by a $1$ and each zero immediately to the left of a
nonzero entry by $-1$, leaving all other entries $0$.
In the case of the above example we obtain in this way
\beq
A=
\left[\begin{array}{rrrrrrrrr}
0 & 0 & 1 & 0 & 0 & 0 & 0 & 0 & 0 \\
0 & 0 & 0 & 1 & 0 & 0 & 0 & 0 & 0 \\
1 & 0 &-1 & 0 & 0 & 1 & 0 & 0 & 0 \\
0 & 0 & 0 & 0 & 1 &-1 & 0 & 0 & 1 \\
0 & 0 & 1 &-1 & 0 & 0 & 0 & 1 & 0 \\
0 & 0 & 0 & 1 &-1 & 1 & 0 &0 & 0 \\
\end{array}\right]
\ \in {\cal A}^{986431}(6)
\label{Ex_Amu}
\eeq


Square ice provides a further refinement of the relationship between 
shifted tableaux and $\mu$--ASM.
Square ice is a directed graph that models the orientation of oxygen and
hydrogen molecules in frozen water. The vertices are laid out in an $n\times n$
grid and each vertex has two incoming and two outgoing edges in a
north, south, east, west orientation. The square ice graph
corresponding to (\ref{Ex_Amu}) appears in Figure \ref{newfig}.

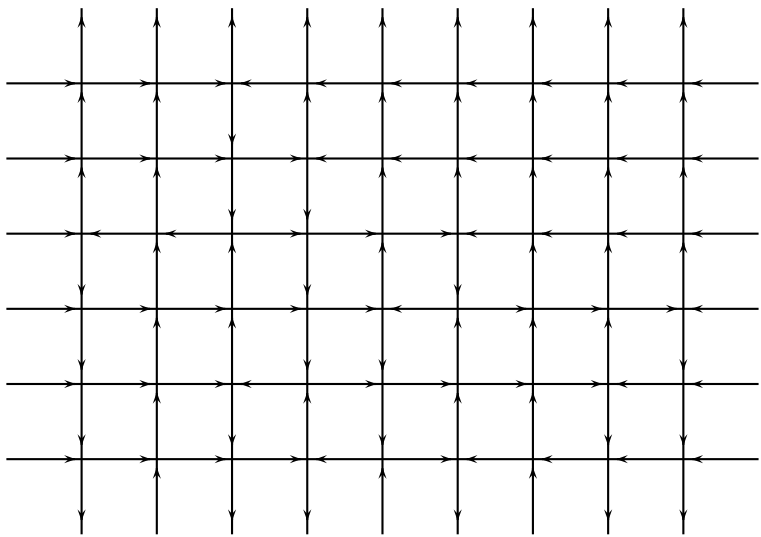
\begin{fullfigure}{newfig}{Square ice for equation (\ref{Ex_Amu})}  


\pspicture(3,0)(7,7)

\psline{->}(.91,1)(.92,1)
\psline{->}(1.91,1)(1.92,1) \psline{->}(2.91,1)(2.92,1)
\psline{->}(3.91,1)(3.92,1) 
\psline{->}(4.12,1)(4.11,1) \psline{->}(5.91,1)(5.92,1)
\psline{->}(6.12,1)(6.11,1) 
\psline{->}(7.12,1)(7.11,1) \psline{->}(8.12,1)(8.11,1)
\psline{->}(9.12,1)(9.11,1) 

\psline{->}(.91,2)(.92,2)
\psline{->}(1.91,2)(1.92,2) \psline{->}(2.91,2)(2.92,2)
\psline{->}(3.12,2)(3.11,2) 
\psline{->}(4.91,2)(4.92,2) \psline{->}(5.91,2)(5.92,2)
\psline{->}(6.91,2)(6.92,2) 
\psline{->}(7.91,2)(7.92,2) \psline{->}(8.12,2)(8.11,2)
\psline{->}(9.12,2)(9.11,2) 

\psline{->}(.91,3)(.92,3)
\psline{->}(1.91,3)(1.92,3) \psline{->}(2.91,3)(2.92,3)
\psline{->}(3.91,3)(3.92,3) 
\psline{->}(4.91,3)(4.92,3) \psline{->}(5.12,3)(5.11,3)
\psline{->}(6.91,3)(6.92,3) 
\psline{->}(7.91,3)(7.92,3) \psline{->}(8.91,3)(8.92,3)
\psline{->}(9.12,3)(9.11,3) 

\psline{->}(.91,4)(.92,4)
\psline{->}(1.12,4)(1.11,4) \psline{->}(2.12,4)(2.11,4)
\psline{->}(3.91,4)(3.92,4) 
\psline{->}(4.91,4)(4.92,4) \psline{->}(5.91,4)(5.92,4)
\psline{->}(6.12,4)(6.11,4) 
\psline{->}(7.12,4)(7.11,4) \psline{->}(8.12,4)(8.11,4)
\psline{->}(9.12,4)(9.11,4) 

\psline{->}(.91,5)(.92,5)
\psline{->}(1.91,5)(1.92,5) \psline{->}(2.91,5)(2.92,5)
\psline{->}(3.91,5)(3.92,5) 
\psline{->}(4.12,5)(4.11,5) \psline{->}(5.12,5)(5.11,5)
\psline{->}(6.12,5)(6.11,5) 
\psline{->}(7.12,5)(7.11,5) \psline{->}(8.12,5)(8.11,5)
\psline{->}(9.12,5)(9.11,5) 

\psline{->}(.91,6)(.92,6)
\psline{->}(1.91,6)(1.92,6) \psline{->}(2.91,6)(2.92,6)
\psline{->}(3.12,6)(3.11,6) 
\psline{->}(4.12,6)(4.11,6) \psline{->}(5.12,6)(5.11,6)
\psline{->}(6.12,6)(6.11,6) 
\psline{->}(7.12,6)(7.11,6) \psline{->}(8.12,6)(8.11,6)
\psline{->}(9.12,6)(9.11,6)

\psline{->}(1,.19)(1,.18)
\psline{->}(1,1.19)(1,1.18)\psline{->}(1,2.19)(1,2.18)
\psline{->}(1,3.19)(1,3.18)\psline{->}(1,4.88)(1,4.89)\psline{->}(1,5.88)(1,5.89)
\psline{->}(1,6.88)(1,6.89)

\psline{->}(2,.88)(2,.89)
\psline{->}(2,1.88)(2,1.89)\psline{->}(2,2.88)(2,2.89)
\psline{->}(2,3.88)(2,3.89)\psline{->}(2,4.88)(2,4.89)\psline{->}(2,5.88)(2,5.89)
\psline{->}(2,6.88)(2,6.89)

\psline{->}(3,.19)(3,.18)
\psline{->}(3,1.19)(3,1.18)\psline{->}(3,2.88)(3,2.89)
\psline{->}(3,3.88)(3,3.89)\psline{->}(3,4.19)(3,4.18)\psline{->}(3,5.19)(3,5.18)
\psline{->}(3,6.88)(3,6.89)

\psline{->}(4,.19)(4,.18)
\psline{->}(4,1.88)(4,1.89)\psline{->}(4,2.19)(4,2.18)
\psline{->}(4,3.19)(4,3.18)\psline{->}(4,4.19)(4,4.18)\psline{->}(4,5.88)(4,5.89)
\psline{->}(4,6.88)(4,6.89)

\psline{->}(5,.88)(5,.89)
\psline{->}(5,1.19)(5,1.18)\psline{->}(5,2.19)(5,2.18)
\psline{->}(5,3.88)(5,3.89)\psline{->}(5,4.88)(5,4.89)\psline{->}(5,5.88)(5,5.89)
\psline{->}(5,6.88)(5,6.89)

\psline{->}(6,.19)(6,.18)
\psline{->}(6,1.88)(6,1.89)\psline{->}(6,2.88)(6,2.89)
\psline{->}(6,3.19)(6,3.18)\psline{->}(6,4.88)(6,4.89)\psline{->}(6,5.88)(6,5.89)
\psline{->}(6,6.88)(6,6.89)

\psline{->}(7,.88)(7,.89)
\psline{->}(7,1.88)(7,1.89)\psline{->}(7,2.88)(7,2.89)
\psline{->}(7,3.88)(7,3.89)\psline{->}(7,4.88)(7,4.89)\psline{->}(7,5.88)(7,5.89)
\psline{->}(7,6.88)(7,6.89)

\psline{->}(8,.19)(8,.18)
\psline{->}(8,1.19)(8,1.18)\psline{->}(8,2.88)(8,2.89)
\psline{->}(8,3.88)(8,3.89)\psline{->}(8,4.88)(8,4.89)\psline{->}(8,5.88)(8,5.89)
\psline{->}(8,6.88)(8,6.89)

\psline{->}(9,.19)(9,.18)
\psline{->}(9,1.19)(9,1.18)\psline{->}(9,2.19)(9,2.18)
\psline{->}(9,3.88)(9,3.89)\psline{->}(9,4.88)(9,4.89)\psline{->}(9,5.88)(9,5.89)
\psline{->}(9,6.88)(9,6.89)


\psline(0,1)(10,1) \psline(0,2)(10,2)
\psline(0,3)(10,3) \psline(0,4)(10,4)
\psline(0,5)(10,5) \psline(0,6)(10,6)

\psline(1,0)(1,7) \psline(2,0)(2,7)
\psline(3,0)(3,7) \psline(4,0)(4,7)
\psline(5,0)(5,7) \psline(6,0)(6,7)
\psline(7,0)(7,7) \psline(8,0)(8,7)
\psline(9,0)(9,7)

\endpspicture
\end{fullfigure}

At each vertex there are six possible orientations of the four
directed edges. These orientations may be specified by the 
pairs of compass points giving the directions of the incoming 
edges. In this way the above square ice graph is
specified by a corresponding ``compass points'' matrix:   

\beq
CM\ =\
\left[\begin{array}{rrrrrrrrr}
\ssc{NE} &\ssc{NE} &\ssc{WE} &\ssc{NW} &\ssc{NW} &\ssc{NW} &\ssc{NW} &\ssc{NW} &\ssc{NW} \\
\ssc{NE} &\ssc{NE} &\ssc{SE} &\ssc{WE} &\ssc{NW} &\ssc{NW} &\ssc{NW} &\ssc{NW} &\ssc{NW} \\
\ssc{WE} &\ssc{NW} &\ssc{NS} &\ssc{SE} &\ssc{NE} &\ssc{WE} &\ssc{NW} &\ssc{NW} &\ssc{NW} \\
\ssc{SE} &\ssc{NE} &\ssc{NE} &\ssc{SE} &\ssc{WE} &\ssc{NS} &\ssc{NE} &\ssc{NE} &\ssc{WE} \\
\ssc{SE} &\ssc{NE} &\ssc{WE} &\ssc{NS} &\ssc{SE} &\ssc{NE} &\ssc{NE} &\ssc{WE} &\ssc{SW} \\
\ssc{SE} &\ssc{NE} &\ssc{SE} &\ssc{WE} &\ssc{NS} &\ssc{WE} &\ssc{NW} &\ssc{SW} &\ssc{SW} \\
\end{array}\right]
\eeq

The bijection between compass point matrices, square ice 
graphs and $\mu$-ASMs is provided by the following correspondences:

\begin{fullfigure}{V.1}
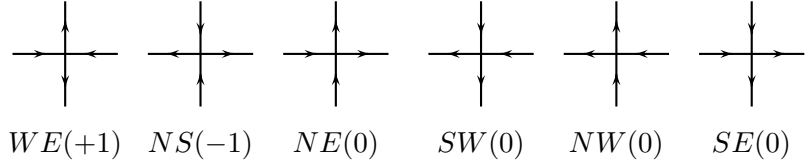
    {Square ice and corresponding compass points and ASM entries. Figure adapted from Kuperberg \cite{K02}}
\verta{WE (+1)}\hspace{.4cm}\vertb{NS (-1)}\hspace{.4cm}\vertc{NE (0)}\hspace{.4cm}
\vertd{SW (0)}\hspace{.4cm}\verte{NW (0)}\hspace{.4cm}\vertf{SE (0)}\eatline
\end{fullfigure}

The horizontal orientation (with both horizontal edges directed in),
$WE$, corresponds
to each entry $+1$ in $A$, and the vertical orientation (with both vertical edges
directed in), $NS$, 
corresponds to each entry $-1$ in $A$; the other four orientations, $NE$, $SW$, $NW$
and $SE$ correspond to the entries $0$ in $A$.
Accordingly there are northwest zeros (with edges pointing in the north
and west directions), southwest zeros, northeast zeros, and southeast zeros.
Northwest zeros are those whose nearest nonzero neighbour to the right, if it
has one,  is $-1$, and whose nearest nonzero neighbour below is $1$.
Southwest zeros are those whose nearest nonzero neighbour to the right, if
it has one, is $-1$, and whose nearest nonzero neighbour below, if it has one,
is $-1$.
Northeast zeros are those whose nearest nonzero neighbour to the right
is $1$, and whose nearest nonzero neighbour below is $1$.
Southeast zeros are those whose nearest nonzero neighbour to the right
is $1$, and whose nearest nonzero neighbour below, if it has one, is $-1$.

The compass points matrices $CM$ can then be associated 
to the set of all primed shifted tableaux $PST$ that may be obtained 
by adding primes to the entries of the unprimed tableau $ST$.  For example,
the entries $NE$ in the $k$th row are 
associated with an entry $k$ in $PST$
and correspondingly to a weight factor $x_k$.
The entries $SE$ in the $k$th row
are associated with an entry $k'$ in $PST$
and correspondingly to a weight factor $y_k$.
The entries $NS$ in the $k$th row 
are associated with the two possible labels $k$ and
$k'$ of the first
box of each connected component of 
$\str_k(PST)$ other than the one starting on the main diagonal.
Correspondingly each $NS$ in row $k$ is associated 
with a weight factor $(x_k+y_k)$. 
It should 
be pointed out that the main diagonal is not included at all in the compass points
matrix so that the first column corresponds to the second diagonal and indeed in general,
column $k$ of $CM$ corresponds to diagonal $k+1$ of $PST$.  This implies that the above
weighting excludes the weight $x_1x_2\cdots x_n$
arising from the entries $1,2,\ldots,n$ on the main 
diagonal of each $PST$.

Combining the weight factors we have a total weight
associated with each $A\in{\cal A}^\mu(n)$ given by
\beq
\prod_{k=1}^n\ 
x_k^{NE_k(A)}\ y_k^{SE_k(A)}\ (x_k+y_k)^{NS_k(A)}
\eeq
where $NE_k(A)$, $SE_k(A)$ and $NS_k(A)$ 
are the numbers of entries $NE$, $SE$ and $NS$ 
in the $k$th row of the compass matrix $CM(A)$
corresponding to $A$.

Thanks to the connection already made between $PST$s and 
weighted $ST$s, the following is then an immediate corollary
of Proposition~\ref{Prop-PQ}:

\begin{Corollary}
Let $\mu=\lambda+\delta$ be a strict partition of length $\ell(\mu)=n$,
with $\lambda$ a partition of length $\ell(\lambda)<n$ and \
$\delta=(n,n-1,\ldots,1)$. Then for all $\x=(x_1,x_2,\ldots,x_n)$ and 
$\y=(y_1,y_2,\dots,y_n)$ we have
\beq
 \prod_{1\leq i < j\leq n}\ (x_i + y_j)\ s_{\lambda}(\x) =
\sum_{A\in {\cal A}^{\mu}(n)} \prod_{k=1}^{n}
x_{k}^{NE_{k}(A)} y_{k}^{SE_{k}(A)} (x_k + y_k)^{NS_{k}(A)}.
\label{Eq-5.1}
\eeq
\label{Cor-5.1}
\end{Corollary}

This generalises a result of Chapman \cite{C01}.  In his original
paper he weights by column instead of row so the parameters in his
paper correspond to the transpose matrix. Now setting $\lambda=0$
so that $\mu=\delta$, and noting that
${\cal A}^\delta(n)={\cal A}(n)$, the set of all $n\times n$
ASMs, we have

\begin{Corollary}[Chapman \cite{C01}]
\beq
\prod_{1\leq i < j\leq n} (x_i + y_j) 
= \sum_{A\in {\cal A}(n)} \prod_{k=1}^{n}
x_{k}^{NE_{k}(A)} y_{k}^{SE_{k}(A)} (x_k + y_k)^{NS_{k}(A)}.
\eeq
\end{Corollary}

Corollary~\ref{Cor-5.1} has a further consequence:

\begin{Corollary}
Let $\mu=\lambda+\delta$ be a strict partition of length $\ell(\mu)=n$,
with $\lambda$ a partition of length $\ell(\lambda)\leq n$ and \
$\delta=(n,n-1,\ldots,1)$. For any $m$ for which $m>n$ and $\mu_1\leq m$,
let ${\cal A}(n,m,\mu)\in{\cal A}(m)$ be the subset consisting
of those ASMs, $C$, whose top $n$ rows constitute
an ASM, $A$, in ${\cal A}^\mu(n)$.
Then for all $(x_1,\ldots,x_n,x_{n+1},\ldots,x_m)$ and 
$(y_1,\dots,y_n,y_{n+1},\ldots,y_m)$ we have
\begin{eqnarray}
 &&\prod_{1\leq i < j\leq n}(x_i + y_j)\ \ s_{\lambda}(x_1,\ldots,x_n)\ 
 \prod_{n+1\leq i < j\leq m}(x_i + y_j)\ \ s_{\kappa}(y_{n+1},\ldots,y_{m})\cr
&&\qquad=
\sum_{C\in {\cal A}(n,m,\mu)} \prod_{k=1}^{m}
x_{k}^{NE_{k}(C)} y_{k}^{SE_{k}(C)} (x_k + y_k)^{NS_{k}(C)},
\label{Eq-5.new}
\end{eqnarray}
where $\kappa$ is the conjugate of the complement of $\lambda$ with 
respect to the rectangular partition $((m-n)^n)$, that is 
$\kappa=((m-n)^n/\lambda)'$.
\label{Cor-5.new}
\end{Corollary}

\noindent{\bf Proof}\ \
Let the top $n$ rows of $C$ and the bottom $(m-n)$ rows of $C$, 
reversed in order, form the matrices $A$ and $B$, respectively.
Then the application of Corollary~\ref{Cor-5.1} to $A$
gives the contribution made by the top $n$ rows of each $C$ 
on the right hand side of (\ref{Eq-5.new}) in the form of the first 
two factors on the left hand side.
Similarly, the remaining two factors on the left hand side arise
from the contribution of the bottom $(m-n)$ rows of each $C$.
To see this one applies Corollary~\ref{Cor-5.1} to $B$, but
this time with $\mu=\lambda+\delta$ replaced by $\nu=\kappa+\epsilon$,
where $\epsilon=(m-n,m-n-1,\ldots,1)$, and with $(x_1,\ldots,x_n)$ 
and $(y_1,\dots,y_n)$ replaced by $(y_m,,\dots,y_{n+1})$
and $(x_m,,\ldots,x_{n+1})$, respectively. It only remains to 
relate $\lambda$ and $\kappa$. Since the parts of $\mu$ and $\nu$ 
specify those columns of $A$ and $B$, respectively, whose column 
sums are $1$, and $A$ and $B$ are constructed from an ASM $C$,
all these parts must be distinct and together constitute 
$(m,m-1,\ldots,1)$. It follows that the union of  
$\{\lambda_i+n-i+1\,|\,i=1,\ldots,n\}$ and
$\{\kappa_j+(m-n)-j+1\,|\,j=1,\ldots,m-n\}$ must be 
$\{1,\ldots,m\}$. However, it is well known~\cite{M95}p3 that
the complement of $\{\lambda_i+n-i+1\,|\,i=1,\ldots,n\}$ in
$\{1,\ldots,m\}$ is $\{n+k-\lambda'_k\,|\,k=1,\ldots,m-n\}$. 
By setting $k=m-n-j+1$ it can then be seen that 
$\kappa_j=n-\lambda'_{m-n-j+1}$ for $j=1,\ldots,m-n$, so 
that $\kappa=(n^{m-n}/\lambda'=((m-n)^n/\lambda)'$, as claimed.  

Alternatively, Corollary~\ref{Cor-5.new}, may be proved bijectively
by taking each primed shifted tableau $PST$ specified by some 
$C\in{\cal A}(n,m,\mu)$ and using the jeu de taquin, first as described
above, to move all entries $k'$ with $1\leq k\leq n$ north-west
to the $k$th column, and then in an analogous manner, to move all
entries $k$ with $n+1\leq k\leq m$ south-east to the $k$th row. The result
is a pair of triangular subtableaux, both of type $PD$ but with all
entries $k$ and $k'$ such that $k\leq n$ in one case and $k>n$ in the other,
together with a pair of semistandard tableaux, one of shape $\lambda$ 
with unprimed entries subject to the order relation $1<2<\cdots<n$ and 
the other of shape $\kappa$ with primed entries subject to the order 
relation $m'<(m-1)'<\cdots<(n+1)'$.  

By way of an example, let $m=6$ and consider the following $6\times6$ ASM
\beq
C=
\left[\begin{array}{rrrrrr}
0 & 0 & 1 & 0 & 0 & 0 \\
0 & 1 &-1 & 1 & 0 & 0 \\
1 &-1 & 1 & 0 & 0 & 0 \\
0 & 0 & 0 & 0 & 1 & 0 \\
0 & 1 & 0 & 0 &-1 & 1 \\
0 & 0 & 0 & 0 & 1 & 0 \\
\end{array}\right]\ \ \in{\cal A}(6)\,.
\label{Ex-A-new}
\eeq
Taking $n=2$, the top two rows of $C$ constitute $A$, and the bottom four 
rows of $C$, reversed in order, constitute $B$, where:
\beq
A=\ \left[\begin{array}{rrrrrr}
0 & 0 & 1 & 0 & 0 & 0 \\
0 & 1 &-1 & 1 & 0 & 0 \\
\end{array}\right]\ \ \in{\cal A}^{4,2}(2)
\quad\hbox{and}\quad
B=\ \left[\begin{array}{rrrrrr}
0 & 0 & 0 & 0 & 1 & 0 \\
0 & 1 & 0 & 0 &-1 & 1 \\
0 & 0 & 0 & 0 & 1 & 0 \\
1 &-1 & 1 & 0 & 0 & 0 \\
\end{array}\right]\ \ \in{\cal A}^{6,5,3,1}(4)\,.
\label{Ex-AA-new}
\eeq
The superscripts on ${\cal A}^{4,2}(2)$ and ${\cal A}^{6,5,3,1}(4)$
specify those columns of $A$ and $B$, respectively, having column sums $1$. 
They indicate that $A\in{\cal A}(2,4,\mu)$ with $\mu=(4,2)$ so that
$\lambda=(2,1)$, while $\nu=(6,5,3,1,)$ so that $\kappa=(2,2,1)$. This
is in accordance with the formula $\kappa=(4^2/\lambda)'=(3,2)'=(2,2,1)$.

The compass point matrix corresponding to $C$ takes the form:
\beq
CM\ =\
\left[\begin{array}{rrrrrrrrr}
\ssc{NE} &\ssc{NE} &\ssc{WE} &\ssc{NW} &\ssc{NW} &\ssc{NW} \\
\ssc{NE} &\ssc{WE} &\ssc{NS} &\ssc{WE} &\ssc{NW} &\ssc{NW} \\
\ssc{WE} &\ssc{NS} &\ssc{WE} &\ssc{SW} &\ssc{NW} &\ssc{NW} \\
\ssc{SE} &\ssc{NE} &\ssc{SE} &\ssc{SE} &\ssc{WE} &\ssc{NW} \\
\ssc{SE} &\ssc{WE} &\ssc{SW} &\ssc{SW} &\ssc{NS} &\ssc{WE} \\
\ssc{SE} &\ssc{SE} &\ssc{SE} &\ssc{SE} &\ssc{WE} &\ssc{SW} \\
\end{array}\right]
\eeq
so that the contribution of $C$ to the right hand side of (\ref{Eq-5.new})
is
\beq
x_1^2\,x_2\,(x_2+y_2)\,(x_3+y_3)\,x_4\,y_4^3\,y_5\,(x_5+y_5)\,y_6^4.
\label{Eq-xy-new}
\eeq

The three factors $(x_k+y_k)$ arise from three entries $NS$ in $CM$,
that themselves arise from the three entries $-1$ in $C$. 
There must be therefore be precisely $2^3$ primed shifted tableaux 
$PST$ corresponding to $C$. Choosing just one of these for illustrative 
purposes, the use of the jeu de taquin to move all $k'$s with 
$k\leq 2$ north-west and all $k$s with $k>2$ south-east, gives
the following bijective map:
\beq
PST\ \ =\ \
{\vcenter
{\offinterlineskip
\halign{&\mystrut\vrule#&\mybox{\hss$\scriptstyle#$\hss}\cr
\hr{13}\cr 
&1&&1&&1&&2'&&4'&&5&\cr
\hr{13}\cr 
\omit&
&&2&&2&&3&&4'&&6'&\cr
\nr{2}&\hr{11}\cr 
\omit& &\omit&
&&3&&4'&&4&&6'&\cr 
\nr{4}&\hr{9}\cr 
\omit& &\omit& &\omit&
&&4&&5'&&6'&\cr
\nr{6}&\hr{7}\cr 
\omit& &\omit& &\omit& &\omit&
&&5&&6'&\cr 
\nr{8}&\hr{5}\cr
\omit& &\omit& &\omit& &\omit& &\omit&
&&6&\cr 
\nr{10}&\hr{3}\cr
 }}}
\ \ \longleftrightarrow\ \ 
{\vcenter
{\offinterlineskip
\halign{&\mystrut\vrule#&\mybox{\hss$\scriptstyle#$\hss}\cr
\hr{13}\cr 
&1&&2'&&1&&1&&4'&&6'&\cr 
\hr{13}\cr 
\omit& &&2&&2&&4'&&5'&&6'&\cr 
\nr{2}&\hr{11}\cr 
\omit& &\omit&
&&3&&4'&&3&&6'&\cr 
\nr{4}&\hr{9}\cr 
\omit& &\omit& &\omit&
&&4&&4&&6'&\cr
\nr{6}&\hr{7}\cr 
\omit& &\omit& &\omit& &\omit&
&&5&&5&\cr 
\nr{8}&\hr{5}\cr
\omit& &\omit& &\omit& &\omit& &\omit&
&&6&\cr 
\nr{10}&\hr{3}\cr
 }}}
\label{Eq-PST-new}
\eeq
The corresponding contribution to the left hand side
of (\ref{Eq-5.new}) is then given by
\beq
x_1^2\,x_2\,y_2\,x_3\,x_4\,y_4^3\,x_5\,y_5\,y_6^4
=
y_2 \cdot x_1^2\,x_2 \cdot x_3\,x_4\,y_4\,x_5\,y_6^2 \cdot y_4^2\,y_5\,y_6^2,
\label{Eq-xyfact-new}
\eeq
where the arrangement of the terms on the right exhibits the
contributions to each of the four factors constituting the
left hand side of (\ref{Eq-5.new}). Both tableaux in (\ref{Eq-PST-new}) 
may be displayed, as shown below, in terms of suitably re-oriented 
subtableaux involving entries $k$ and $k'$, with all $k\leq2$ in one case, 
and all $k>2$ on the other.
\beq
{\vcenter
{\offinterlineskip
\halign{&\mystrut\vrule#&\mybox{\hss$\scriptstyle#$\hss}\cr
\hr{9}\cr 
&1&&1&&1&&2'&\cr
\hr{9}\cr 
\omit&
&&2&&2&\cr
\nr{2}&\hr{5}\cr 
}}}
\ \ \cdot\ \
{\vcenter
{\offinterlineskip
\halign{&\mystrut\vrule#&\mybox{\hss$\scriptstyle#$\hss}\cr
\hr{13}\cr 
&6&&6'&&6'&&6'&&6'&&5&\cr
\hr{13}\cr 
\omit&
&&5&&5'&&4&&4'&&4'&\cr
\nr{2}&\hr{11}\cr 
\omit& &\omit&
&&4&&4'&&3&\cr 
\nr{4}&\hr{7}\cr 
\omit& &\omit& &\omit&
&&3&\cr
\nr{6}&\hr{3}\cr 
 }}}
\ \ \longleftrightarrow \ \
{\vcenter
{\offinterlineskip
\halign{&\mystrut\vrule#&\mybox{\hss$\scriptstyle#$\hss}\cr
\hr{5}\cr 
&1&&2'&\cr 
\hr{5}\cr 
\omit& &&2&\cr 
\nr{2}&\hr{3}\cr 
}}}
\ \ \cdot\ \ 
{\vcenter
{\offinterlineskip
\halign{&\mystrut\vrule#&\mybox{\hss$\scriptstyle#$\hss}\cr
\hr{5}\cr
&1&&1&\cr 
\hr{5}\cr 
&2&\cr
\hr{3}\cr 
}}}
\ \ \cdot\ \ 
{\vcenter
{\offinterlineskip
\halign{&\mystrut\vrule#&\mybox{\hss$\scriptstyle#$\hss}\cr
\hr{9}\cr 
&6&&5&&6'&&6'&\cr
\hr{9}\cr 
\omit&
&&5&&4&&3&\cr
\nr{2}&\hr{7}\cr 
\omit& &\omit& 
&&4&&4'&\cr 
\nr{4}&\hr{5}\cr
\omit& &\omit& &\omit&
&&3&\cr 
\nr{6}&\hr{3}\cr
 }}}
\ \ \cdot\ \
{\vcenter
{\offinterlineskip
\halign{&\mystrut\vrule#&\mybox{\hss$\scriptstyle#$\hss}\cr
\hr{5}\cr 
&6'&&6'&\cr 
\hr{5}\cr 
&5'&&4'&\cr 
\hr{5}\cr 
&4'&\cr
\hr{3}\cr
}}}
\label{Eq-PST-PD-new}
\eeq
This illustrates the outcome of applying the jeu de taquin to
primed shifted tableaux corresponding to the submatrices $A$ 
and $B$ of the ASM $C$. The resulting contribution of the four final tableaux to
the left hand side of (\ref{Eq-5.new}), is then confirmed to be as given 
on the right hand side of (\ref{Eq-xyfact-new}).

\subsection{$sp(2n)$ case}\ \
\label{subsec:sp-case}

The symplectic case involves a modified alternating sign matrix called
a U--turn $\mu$--ASM or $\mu$--UASM.  
Informally, the U--turn condition means that two consecutive rows and
the U--turn between them must follow 
the $\mu$--ASM summation rules, ASM2--5; that is, the cumulative sum
must be zero or one, and the total sum 
must be one.  These $\mu$--UASM were first defined in Hamel and 
King~\cite{HK02} 
where they were called 
sp($2n$)--generalised alternating sign matrics.  
They are discussed at length in Hamel and King \cite{HK03}.  
A formal definition is as follows:

Let $\mu$ be a partition of length $\ell(\mu)=n$, all of whose parts 
are distinct, and for which $\mu_1\leq m$. Then the matrix 
$UA=(a_{iq})_{1\leq i\leq 2n,1\leq q\leq m}$ is said to belong to the 
set ${\cal UA}^\mu(2n)$ of 
$\mu$--alternating sign matrices with a U--turn boundary if it is
a $2n\times m$ matrix whose elements $a_{iq}$ satisfy the conditions: 

\begin{tabular}{rll}
 UA1 &$a_{iq}\in\{-1,0,1\} $
&\hbox{for $1\leq i\leq 2n$, $1\leq q\leq m$};\\
UA2&$\sum_{q=p}^m a_{iq}\in\{0,1\}$ 
&\hbox{for $1\leq i\leq 2n$, $1\leq p\leq m$}; \\
UA3&$\sum_{i=j}^{2n} a_{iq}\in\{0,1\}$ 
&\hbox{for $1\leq j\leq 2n$, $1\leq q\leq m$}. \\
 UA4 &$\sum_{q=1}^m (a_{2i-1,q}+a_{2i,q})=1 $
&\hbox{for $1\leq i\leq n$};\\
 UA5 &$\sum_{i=1}^{2n} a_{iq}=
  \bigg\{ \begin{array}{ll}
 1& \mbox{if } q=\mu_k \mbox{ for some $k$}\\ 
 0&\mbox{otherwise}
\end{array} $
 &\hbox{for $1\leq q\leq m$, $1\leq k\leq n$}.\cr
\end{tabular}

In the case $\mu=\delta=(n,n-1,\ldots,1)$ and $m=n$, for 
which UA5 becomes $\sum_{i=1}^{2n}a_{iq}=1$ for $1\leq q\leq n$, 
this definition is such that the set of $\mu$--UASM coincides 
with the set of U--turn alternating sign matrices, UASMs, defined by 
Kuperberg~\cite{K02}. 

As noted above, Hamel and King \cite{HK03} established a bijection 
between $\mu$--UASM and semistandard shifted symplectic tableaux.  
An example of this association is illustrated below
in the case $n=5$ and $\mu=(9,7,6,2,1)$:

\beq 
{\vcenter {\offinterlineskip
\halign{&\mystrut\vrule#&\mybox{\hss$\scriptstyle#$\hss}\cr \hr{19}\cr
&\ov1&&1&&\ov2&&2&&\ov3&&\ov3&&\ov4&&4&&5&\cr
\hr{19}\cr
\omit& &&\ov2&&\ov2&&2&&3&&\ov4&&\ov4&&4&\cr
\nr{2}&\hr{15}\cr \omit& &\omit& &&3&&\ov4&&4&&4&&4&&4&\cr 
\nr{4}&\hr{13}\cr 
\omit& &\omit& &\omit& &&4&&4&\cr
\nr{6}&\hr{5}\cr
\omit& &\omit& &\omit& &\omit&&&\ov5&\cr 
\nr{8}&\hr{3}\cr }}}
\Longrightarrow \left[\begin{array}{rrrrrrrrr}
\ov1 & 0 & 0 & 0 & 0 & 0 & 0 & 0 & 0 \\
0 & 1 & 0 & 0 & 0 & 0 & 0 & 0 & 0 \\
\ov2 & \ov2 & \ov2 & 0 & 0 & 0 & 0 & 0 & 0 \\
0 & 0 & 2 & 2 & 0 & 0 & 0 & 0 & 0 \\
0 & 0 & 0 & 0 & \ov3 & \ov3 & 0 & 0 & 0 \\
3 & 0 & 0 & 3 & 0 & 0 & 0 & 0 & 0 \\
0 & \ov4 & 0 & 0 & \ov4 & \ov4 & \ov4 & 0 & 0 \\
4 & 4 & 4 & 4 & 4 & 4 & 4 & 4 & 0 \\
\ov5 & 0 & 0 & 0 & 0 & 0 & 0 & 0 & 0 \\
0 & 0 & 0 & 0 & 0 & 0 & 0 & 0 & 5 \\
\end{array}\right]
\Longrightarrow \left[\begin{array}{rrrrrrrrr}
1 & 0 & 0 & 0 & 0 & 0 & 0 & 0 & 0 \\
-1 & 1 & 0 & 0 & 0 & 0 & 0 & 0 & 0 \\
0 & 0 & 1 & 0 & 0 & 0 & 0 & 0 & 0 \\
0 & -1 & 0 & 1 & 0 & 0 & 0 & 0 & 0 \\
0 & 0 & 0 & -1 & 0 & 1 & 0 & 0 & 0 \\
1 & 0 & -1 & 1 & 0 & 0 & 0 & 0 & 0 \\
-1 & 1 & 0 & -1 & 0 & 0 & 1 & 0 & 0 \\
0 & 0 & 0 & 0 & 0 & 0 & 0 & 1 & 0 \\
1 & 0 & 0 & 0 & 0 & 0 & 0 & 0 & 0 \\
0 & 0 & 0 & 0 & 0 & 0 & 0 & -1 & 1\\
\end{array}\right]
\label{IV.1}
\eeq 
where the columns are labeled from left 
to right $1,2,\ldots,9=m=\mu_1$, and the rows from
top to bottom $\ov1,1,\ov2,2,\ldots,\ov5,5=n$.

The translation to square ice is also natural and just requires a
modification of the left boundary by the insertion of a U--turn. 
The square ice graph in Figure~\ref{V.1b} corresponds to the above 
$\mu$-UASM matrix. The same example appeared in Hamel and King \cite{HK03}.

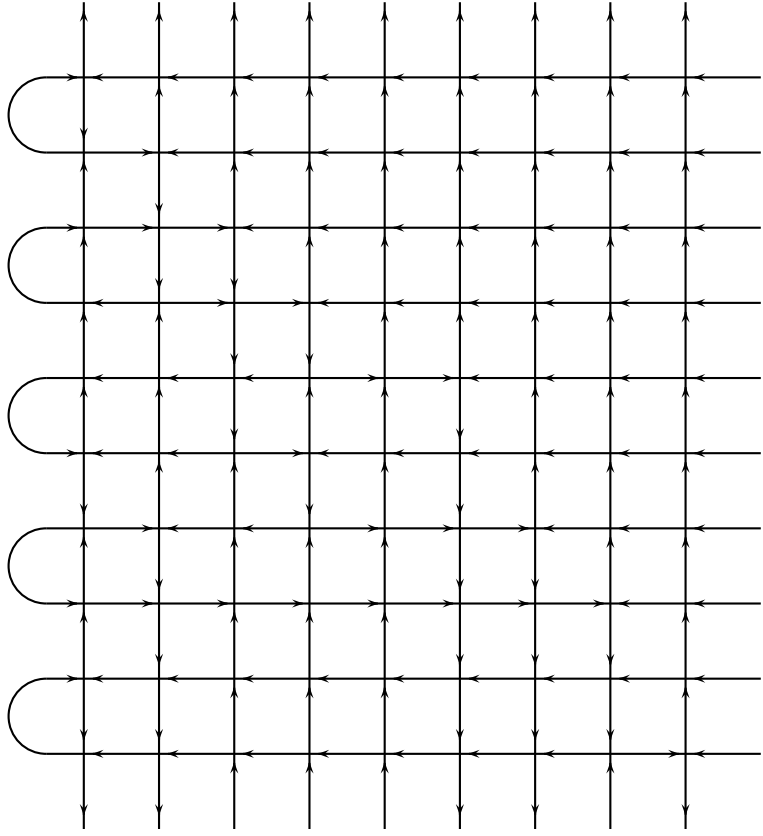
\begin{fullfigure}{V.1b}{Square ice with U--turn boundary} \pspicture(-.5,-.5)(11,12)  
\psline{->}(1.12,1)(1.11,1) \psline{->}(2.12,1)(2.11,1) \psline{->}(3.12,1)(3.11,1)  
\psline{->}(4.12,1)(4.11,1) \psline{->}(5.12,1)(5.11,1) \psline{->}(6.12,1)(6.11,1)  
\psline{->}(7.12,1)(7.11,1) \psline{->}(8.91,1)(8.92,1) \psline{->}(9.12,1)(9.11,1) 
  
\psline{->}(1.12,2)(1.11,2) \psline{->}(2.12,2)(2.11,2) \psline{->}(3.12,2)(3.11,2) 
 \psline{->}(4.12,2)(4.11,2) \psline{->}(5.12,2)(5.11,2) \psline{->}(6.12,2)(6.11,2)  
\psline{->}(7.12,2)(7.11,2) \psline{->}(8.12,2)(8.11,2) \psline{->}(9.12,2)(9.11,2)
   
\psline{->}(1.91,3)(1.92,3) \psline{->}(2.91,3)(2.92,3) \psline{->}(3.91,3)(3.92,3) 
 \psline{->}(4.91,3)(4.92,3) \psline{->}(5.91,3)(5.92,3) \psline{->}(6.91,3)(6.92,3)  
\psline{->}(7.91,3)(7.92,3) \psline{->}(8.12,3)(8.11,3) \psline{->}(9.12,3)(9.11,3) 
  
\psline{->}(1.91,4)(1.92,4) \psline{->}(2.12,4)(2.11,4) \psline{->}(3.12,4)(3.11,4) 
 \psline{->}(4.91,4)(4.92,4) \psline{->}(5.91,4)(5.92,4) \psline{->}(6.91,4)(6.92,4) 
 \psline{->}(7.12,4)(7.11,4) \psline{->}(8.12,4)(8.11,4) \psline{->}(9.12,4)(9.11,4) 
 
 \psline{->}(1.12,5)(1.11,5) \psline{->}(2.12,5)(2.11,5) \psline{->}(3.91,5)(3.92,5)  
\psline{->}(4.12,5)(4.11,5) \psline{->}(5.12,5)(5.11,5) \psline{->}(6.12,5)(6.11,5) 
 \psline{->}(7.12,5)(7.11,5) \psline{->}(8.12,5)(8.11,5) \psline{->}(9.12,5)(9.11,5) 
 
 \psline{->}(1.12,6)(1.11,6) \psline{->}(2.12,6)(2.11,6) \psline{->}(3.12,6)(3.11,6)  
\psline{->}(4.91,6)(4.92,6) \psline{->}(5.91,6)(5.92,6) \psline{->}(6.12,6)(6.11,6)  
\psline{->}(7.12,6)(7.11,6) \psline{->}(8.12,6)(8.11,6) \psline{->}(9.12,6)(9.11,6)
  
 \psline{->}(1.12,7)(1.11,7) \psline{->}(2.91,7)(2.92,7) \psline{->}(3.91,7)(3.92,7)  
\psline{->}(4.12,7)(4.11,7) \psline{->}(5.12,7)(5.11,7) \psline{->}(6.12,7)(6.11,7)  
\psline{->}(7.12,7)(7.11,7) \psline{->}(8.12,7)(8.11,7) \psline{->}(9.12,7)(9.11,7)
  
 \psline{->}(1.91,8)(1.92,8) \psline{->}(2.91,8)(2.92,8) \psline{->}(3.12,8)(3.11,8) 
 \psline{->}(4.12,8)(4.11,8) \psline{->}(5.12,8)(5.11,8) \psline{->}(6.12,8)(6.11,8)  
\psline{->}(7.12,8)(7.11,8) \psline{->}(8.12,8)(8.11,8) \psline{->}(9.12,8)(9.11,8) 
  
\psline{->}(1.91,9)(1.92,9) \psline{->}(2.12,9)(2.11,9) \psline{->}(3.12,9)(3.11,9)  
\psline{->}(4.12,9)(4.11,9) \psline{->}(5.12,9)(5.11,9) \psline{->}(6.12,9)(6.11,9) 
 \psline{->}(7.12,9)(7.11,9) \psline{->}(8.12,9)(8.11,9) \psline{->}(9.12,9)(9.11,9) 
 
 \psline{->}(1.12,10)(1.11,10) \psline{->}(2.12,10)(2.11,10) \psline{->}(3.12,10)(3.11,10) 
 \psline{->}(4.12,10)(4.11,10) \psline{->}(5.12,10)(5.11,10) \psline{->}(6.12,10)(6.11,10) 
 \psline{->}(7.12,10)(7.11,10) \psline{->}(8.12,10)(8.11,10) \psline{->}(9.12,10)(9.11,10)
   
  
\psline{->}(.91,2)(.92,2) \psline{->}(.91,3)(.92,3) \psline{->}(.91,5)(.92,5) 
\psline{->}(.91,8)(.92,8) \psline{->}(.91,10)(.92,10) 

 \psline{->}(1,.19)(1,.18)

 \psline{->}(1,1.19)(1,1.18)\psline{->}(1,2.88)(1,2.89) \psline{->}(1,3.88)(1,3.89)
\psline{->}(1,4.19)(1,4.18)\psline{->}(1,5.88)(1,5.89) \psline{->}(1,6.88)(1,6.89)
\psline{->}(1,7.88)(1,7.89)\psline{->}(1,8.88)(1,8.89) \psline{->}(1,9.19)(1,9.18)

\psline{->}(1,10.88)(1,10.89)  

\psline{->}(2,.19)(2,.18) \psline{->}(2,1.19)(2,1.18)
\psline{->}(2,2.19)(2,2.18) \psline{->}(2,3.19)(2,3.18)\psline{->}(2,4.88)(2,4.89)
\psline{->}(2,5.88)(2,5.89) \psline{->}(2,6.88)(2,6.89)\psline{->}(2,7.19)(2,7.18)
\psline{->}(2,8.19)(2,8.18) \psline{->}(2,9.88)(2,9.89)\psline{->}(2,10.88)(2,10.89)
 
 \psline{->}(3,.88)(3,.89) \psline{->}(3,1.88)(3,1.89)\psline{->}(3,2.88)(3,2.89)
 \psline{->}(3,3.88)(3,3.89)\psline{->}(3,4.88)(3,4.89)\psline{->}(3,5.19)(3,5.18)
 \psline{->}(3,6.19)(3,6.18)\psline{->}(3,7.19)(3,7.18)\psline{->}(3,8.88)(3,8.89)

 \psline{->}(3,9.88)(3,9.89)\psline{->}(3,10.88)(3,10.89) 

 \psline{->}(4,.88)(4,.89) 
\psline{->}(4,1.88)(4,1.89)\psline{->}(4,2.88)(4,2.89) \psline{->}(4,3.88)(4,3.89)
\psline{->}(4,4.19)(4,4.18)\psline{->}(4,5.88)(4,5.89) \psline{->}(4,6.19)(4,6.18)
\psline{->}(4,7.88)(4,7.89)\psline{->}(4,8.88)(4,8.89) \psline{->}(4,9.88)(4,9.89)
\psline{->}(4,10.88)(4,10.89)  

\psline{->}(5,.88)(5,.89) \psline{->}(5,1.88)(5,1.89)
\psline{->}(5,2.88)(5,2.89) \psline{->}(5,3.88)(5,3.89)\psline{->}(5,4.88)(5,4.89)
\psline{->}(5,5.88)(5,5.89) \psline{->}(5,6.88)(5,6.89)\psline{->}(5,7.88)(5,7.89)
\psline{->}(5,8.88)(5,8.89) \psline{->}(5,9.88)(5,9.89)\psline{->}(5,10.88)(5,10.89)

\psline{->}(6,.19)(6,.18) \psline{->}(6,1.19)(6,1.18)\psline{->}(6,2.19)(6,2.18) 
\psline{->}(6,3.19)(6,3.18)\psline{->}(6,4.19)(6,4.18)\psline{->}(6,5.19)(6,5.18) 
\psline{->}(6,6.88)(6,6.89)\psline{->}(6,7.88)(6,7.89)\psline{->}(6,8.88)(6,8.89) 
\psline{->}(6,9.88)(6,9.89)\psline{->}(6,10.88)(6,10.89)  

\psline{->}(7,.19)(7,.18) 
\psline{->}(7,1.19)(7,1.18)\psline{->}(7,2.19)(7,2.18) \psline{->}(7,3.19)(7,3.18)
\psline{->}(7,4.88)(7,4.89)\psline{->}(7,5.88)(7,5.89) \psline{->}(7,6.88)(7,6.89)
\psline{->}(7,7.88)(7,7.89)\psline{->}(7,8.88)(7,8.89) \psline{->}(7,9.88)(7,9.89)
\psline{->}(7,10.88)(7,10.89) 

 \psline{->}(8,.88)(8,.89) \psline{->}(8,1.19)(8,1.18)
\psline{->}(8,2.19)(8,2.18) \psline{->}(8,3.88)(8,3.89)\psline{->}(8,4.88)(8,4.89)
\psline{->}(8,5.88)(8,5.89) \psline{->}(8,6.88)(8,6.89)\psline{->}(8,7.88)(8,7.89)
\psline{->}(8,8.88)(8,8.89) \psline{->}(8,9.88)(8,9.89)\psline{->}(8,10.88)(8,10.89)

 \psline{->}(9,.19)(9,.18) \psline{->}(9,1.88)(9,1.89)\psline{->}(9,2.88)(9,2.89)
 \psline{->}(9,3.88)(9,3.89)\psline{->}(9,4.88)(9,4.89)\psline{->}(9,5.88)(9,5.89)
 \psline{->}(9,6.88)(9,6.89)\psline{->}(9,7.88)(9,7.89)\psline{->}(9,8.88)(9,8.89) 
\psline{->}(9,9.88)(9,9.89)\psline{->}(9,10.88)(9,10.89) 


\psline(.5,1)(10,1) \psline(.5,2)(10,2) \psline(.5,3)(10,3)
 \psline(.5,4)(10,4) \psline(.5,5)(10,5) \psline(.5,6)(10,6) \psline(.5,7)(10,7)
 \psline(.5,8)(10,8) \psline(.5,9)(10,9) \psline(.5,10)(10,10)  \psline(1,0)(1,11) 
\psline(2,0)(2,11) \psline(3,0)(3,11) \psline(4,0)(4,11) \psline(5,0)(5,11) \psline(6,0)(6,11) 
\psline(7,0)(7,11) \psline(8,0)(8,11) \psline(9,0)(9,11)    \psarc(.5,1.5){.5}{90}{270} 
\psarc(.5,3.5){.5}{90}{270} \psarc(.5,5.5){.5}{90}{270} \psarc(.5,7.5){.5}{90}{270} 
\psarc(.5,9.5){.5}{90}{270} 
\endpspicture 
 \end{fullfigure}

In this symplectic case the bijection from the U--turn $\mu$-ASMs to
U--turn square ice graphs is precisely as before, with entries $+1$
and $-1$ mapped to $WE$ and $NS$ vertex orientations, and 
$NE$, $SW$, $NW$ and $SE$ entries $0$ distinguished by their
nearest non-zero neighbouring entries. This map is encoded in the 
corresponding compass points matrix. For the above example, this takes
the form:

\beq
CM=\
\left[
\begin{array}{ccccccccc}
\ssc{WE}&\ssc{NW}&\ssc{NW}&\ssc{NW}&\ssc{NW}&\ssc{NW}&\ssc{NW}&\ssc{NW}&\ssc{NW}\\
\ssc{NS}&\ssc{WE}&\ssc{NW}&\ssc{NW}&\ssc{NW}&\ssc{NW}&\ssc{NW}&\ssc{NW}&\ssc{NW}\\
\ssc{NE}&\ssc{SE}&\ssc{WE}&\ssc{NW}&\ssc{NW}&\ssc{NW}&\ssc{NW}&\ssc{NW}&\ssc{NW}\\
\ssc{NW}&\ssc{NS}&\ssc{SE}&\ssc{WE}&\ssc{NW}&\ssc{NW}&\ssc{NW}&\ssc{NW}&\ssc{NW}\\
\ssc{NW}&\ssc{NW}&\ssc{SW}&\ssc{NS}&\ssc{NE}&\ssc{WE}&\ssc{NW}&\ssc{NW}&\ssc{NW}\\
\ssc{WE}&\ssc{NW}&\ssc{NS}&\ssc{WE}&\ssc{NW}&\ssc{SW}&\ssc{NW}&\ssc{NW}&\ssc{NW}\\
\ssc{NS}&\ssc{WE}&\ssc{NW}&\ssc{NS}&\ssc{NE}&\ssc{SE}&\ssc{WE}&\ssc{NW}&\ssc{NW}\\
\ssc{NE}&\ssc{SE}&\ssc{NE}&\ssc{NE}&\ssc{NE}&\ssc{SE}&\ssc{SE}&\ssc{WE}&\ssc{NW}\\
\ssc{WE}&\ssc{SW}&\ssc{NW}&\ssc{NW}&\ssc{NW}&\ssc{SW}&\ssc{SW}&\ssc{SW}&\ssc{NW}\\
\ssc{SW}&\ssc{SW}&\ssc{NW}&\ssc{NW}&\ssc{NW}&\ssc{SW}&\ssc{SW}&\ssc{NS}&\ssc{WE}\\

\end{array}
\right]
\label{V.2}
\eeq

Then we can generate a weighting in the same manner as for the $gl(n)$
case, with the $k$th column of $CM$ corresponding to the $(k+1)$th
diagonal of $QST$.  
In this case we have unbarred entries corresponding to odd rows and barred entries 
corresponding to even rows.  An entry $NE$ in the $k$th row is
associated to an entry $k$ in $QST$ and correspondingly to a weight
factor $x_k$.  An entry $SE$ in row $k$ is associated to an entry $k'$
in $QST$ and correspondingly to a weight factor $y_k$.  
An entry $NE$ in row $\ov{k}$ is associated to an entry $\ov{k}$ in
$QST$ and correspondingly to a weight factor $t^2x_{k}^{-1}$. 
An entry $SE$ in row $\ov{k}$ is associated to an entry $\ov{K}'$ in
$QST$ and correspondingly to $t^2y_{k}^{-1}$.  
An entry $NS$ in the $k$th row is associated with the four possible
labels $k, k', \ov{k}, \ov{k}'$ of the first box of each connected
component of $str_k(QST)$ (other than the one starting on the main
diagonal) and correspondingly to a weight $(x_i + y_i + t^2 \ov{x}_{i} + t^2\ov{y}_{i})$.
Combining the weight factors we have
\[
\prod_{k=1}^{n} x_{k}^{NE_{k}(A)}(t^2 \ov{x}_{k})^{NE_{\ov{k}}(A)}
y_{k}^{SE_{k}(A)} (t^2\ov{y}_{k})^{SE_{\ov{k}}(A)}
(x_k + y_k)^{NS_{k}(A)} (t^2 \ov{x}_{k} + t^2 \ov{y}_{k})^{NS_{\ov{k}}(A)}
\]
where $SE_{k}(A),\; NE_{k}(A),\; NS_k(A)$ (resp.\ $SE_{\ov{k}}(A),\; 
NE_{\ov{k}}(A),\; NS_{\ov{k}}(A)$) are the numbers of entries 
$SE,\; NE,\; NS$ in row $k$ (resp.\ $\ov{k}$) of the compass matrix $CM(A)$.

We then have the following immediate corollary of Proposition \ref{Prop-PQsp}:
\begin{Corollary}
\begin{eqnarray}
\prod_{1\leq i < j \leq n} (x_i + t^2 \ov{x}_{i} + y_j + t^2 \ov{y}_{j}) sp_{\lambda} (x;t) &=&
\sum_{A} \prod_{k=1}^{n} x_{k}^{NE_{k}(A)}(t^2
\ov{x}_{k})^{NE_{\ov{k}}(A)} y_{k}^{SE_{k}(A)} (t^2\ov{y}_{k})^{SE_{\ov{k}}(A)}\\
 & & \times\ (x_k + y_k)^{NS_{k}(A)} (t^2 \ov{x}_{k} + t^2 \ov{y}_{k})^{NS_{\ov{k}}(A)}
\end{eqnarray}
\label{cor5.3}
\end{Corollary}

This Corollary is a generalisation of Theorem 6.4 of Hamel and King \cite{HK03}.
This Theorem 6.4 may be recovered from Corollary  \ref{cor5.3} by setting $\y =t\x$,
exploiting the bijection between compass point matrices $CM(A)$ and the
U--turn $\mu$-ASM's $A$, and noting that the number of entries $NS$ and
$WE$ in any row of $CM$ are either the same or differ by one 
according to
the nature, barred or unbarred, of the corresponding entry on the main
diagonal of the associated semistandard shifted symplectic tableau $ST$.
Note also that Theorem 6.4 includes the weighting for the main diagonal on 
each side of the equation, whereas Corollary \ref{cor5.3} does not.

\vspace{0.5cm}

\noindent{\bf Acknowledgements}

\bigskip
 The existence of Corollary~\ref{Cor-5.new} was pointed out to us
by an anonymous referee of an earlier version of this paper. 
The first author (AMH) acknowledges the support of a Discovery Grant
from the Natural Sciences and Engineering Research Council of Canada,
and the second (RCK) the support of a Leverhulme Emeritus Fellowship.

\end{document}